\documentclass{amsart}
\usepackage{amscd,amssymb,amsmath,amsthm}
\usepackage[dvips]{graphicx}
\usepackage[all]{xy}
\theoremstyle{plain}
\newtheorem{definition}{Definition}
\newtheorem{proposition}{Proposition}
\newtheorem{theorem}[proposition]{Theorem}

\newtheorem{lemma}[proposition]{Lemma}
\newtheorem{remark}[proposition]{Remark}
\newtheorem*{proposition*}{Proposition}
\newtheorem*{theorem*}{Theorem}
\newtheorem*{corollary*}{Corollary}
\newtheorem*{lemma*}{Lemma}
\newtheorem*{remark*}{Remark}
\newtheorem*{example*}{Example}

\newcommand{\Z}{\mathbb{Z}}
\newcommand{\Q}{\mathbb{Q}}
\newcommand{\R}{\mathbb{R}}
\newcommand{\C}{\mathbb{C}}

\begin{document}

\title{Open Gromov-Witten invariants and Boundary States}

\author{Vito Iacovino}

%\address{Max-Planck-Institut f\"ur Mathematik, Bonn}

\email{vito.iacovino@gmail.com}

\date{version: \today}

%%%%%%%%%%%%%%%%%%%%%%%%%%%%%%%%%%%%%%%%%%%%%%%%%%%%%%%%%%%
%%%%%%%%%%%%%%%%%%%%%%%%%%%%%%%%%%%%%%%%%%%%%%%%%%%%%%%%%%%
%%%%%%%%%%%%%%%%%%%%%%%%%%%%%%%%%%%%%%%%%%%%%%%%%%%%%%%%%%%

\begin{abstract} 
%We define Open Gromov Witten invariants for non-compact lagrangian submanifolds as elements of a vector space associated to the lagrangian submanifold. We show how the quantum Kontsevich-Soibelman algebra arise naturally  in this context.

In this note we show how Kontsevich-Soibelman algebra arise naturally in Open Gromov-Witten theory for not compact geometries.

\end{abstract}

\maketitle

\section{Introduction}

Let $X$ be a Calabi-Yau simplectic six-manifold  and let $L$ be a Maslov index zero Lagrangian submanifold of $X$. Assume that $[L]=0 \in H_3(X, \Z)$ and pick a four chain $K$ such that $\partial K =L$.
% and pick a four chain $K $ such that $\partial K =L$.
In \cite{OGW3} we developed the general  Open Gromov-Witten theory associated to $(X,L)$. We introduced the multi-curve chain complex and , for each $\beta \in H(X,L)$ and $\chi \in \Z$,  we constructed a $MC$-cycle $Z_{\beta,\chi}$ from the moduli-space of  pseudo-holomorphic (multi-)curves . $Z_{\beta,\chi}$ is not uniquely determined but it well defined up to isotopy.   Our theory is the natural mathematical counterpart of the partition function of the A-model Open Topological String introduced by Witten in \cite{W}, and it makes mathematically precise the relation with the perturbative Chern Simons claimed in \cite{W}.  $MC$-cycles (or more precisely what we called nice-$MC$-cycles) can be considered as the mathematical definition of the configuration space of points in the  point splitting perturbative Chern-Simons theory.

In this paper we consider more closely the structures arising in not compact geometries. 
%We assume standard  convexity properties at infinity in order to assure  that the moduli space of  pseudo-holomorphic curves positive area is compact.  In particular we  assume that $L$ has cylindrical ends, that is, we assume that outside a compact set $L$ is diffeomorphic to  $\Sigma \times \R$ , for some riemmanian surface $\Sigma$.
Under standard assumptions on the geometry at infinity, the construction of \cite{OGW3} applies essentially the same way. The main difference consists in setting boundary conditions at infinity  for the perturbation of  area zero curves. This introduces new structures in the problem.   
%The open-closed Gromov-Witten $MC$-cycle depends on the choice of the boundary condition.  
%This problem is present also in closed Gromov-Witten invariants. 
%where changing the boundary condition multiply the for a power series in $g_s$. 
%In the open-closed case there is more beyond this scalar ambiguity. 
%In particular, an important role is played by the area zero annulus with one marked boundary  point.  

%In  \cite{OGW3} we associate to the pair $(X,L)$ together with a four chain $K$ with $\partial K =L$, and $\beta \in H_2(X,L)$ and $$
%   

% In  \cite{OGW3} we associate to $\beta \in H_2(X,L)$ and $\chi \in \Z$ a $MC$-cycle $Z_{\beta,\chi}$ defined up to isotopy.
 In this paper we consider the collection of $MC$-cycles $Z=(Z_{\beta,\chi})_{\beta,\chi}$ for different $\beta,\chi$.  In the not compact case, up to isotopy, the vector $Z$  is determinate by the choice of a boundary condition $\mathcal{B}$ on the perturbation of area zero curves. 
 %We interpret $Z$ as an element of an infinity dimensional vector space $\mathfrak{H}$ .
 %on the field $\Q((g_s ,T))$,  where $g_s$ is the  formal parameter used as weight for the Euler Characteristic of the multi-curves and $T$ is the formal variable of the Novikov ring, such that  $Z \in \mathfrak{H}$.
 %The construction of \cite{OGW3} associate an element $Z \in \mathfrak{H}$ from the moduli space of pseudo-holomorphic curves (that depends on the choice of a boundary condition and a spin structure of $L$).   
%This vector space is related to the perturbative quantization of abelian Chern-Simons theory on $\Sigma \times \R$.
 We interpret $Z$ as an element of an infinity dimensional vector space $\mathfrak{H}_{  \mathfrak{c}  }$, which depends on an element  $\mathfrak{c} \in H_1(L, \Z)$.  The class $\mathfrak{c} $ depends on the boundary conditions, and we refer to it  as the Chern Class.

This result can be considered as the analogous in open-closed Gromov-Witten theory of the main claim  of  Topological Field Theories,  where the partition function on a manifold  with boundary $M$ defines a vector on the space of boundary states, which is defined quantizing the theory on $\partial M \times \R$.  In our case, $M=L$ and the relevant topological field theory is the perturbative abelian Chern-Simons theory.

Let $\Sigma = \partial L$.
Consider the quantum Kontsevich-Soibelman algebra $\hat{\mathfrak{g}}$ associated to 
the abelian group  $ H_1(\Sigma,\Z) $ equipped with the standard intersection  form $\langle \bullet , \bullet \rangle$:
$$  \hat { \mathfrak{g} } = \bigoplus_{\gamma \in  H_1(\Sigma,\Z) } \Q [q^{\frac{1}{2}}, q^{-\frac{1}{2}}] \hat{e}_{\gamma} ,$$
$$\hat{e}_{\gamma_1}  \hat{e}_{\gamma_2} =  q^{ \frac{\langle \gamma_1 , \gamma_2 \rangle}{2}}  \hat{e}_{\gamma_1 + \gamma_2} \text{   for    }  \gamma_1, \gamma_2 \in H_1(\Sigma,\Z) ,$$
%We show that there exists a natural representation on $\mathfrak{H}$ of the quantum KS Kontsevich-Soibelman algebra $\mathfrak{g}$ with
with
$$     q^{\frac{1}{2}} = - e^{ \frac{g_s}{2}}       .$$
%where $g_s$ is the string coupling  used as formal

%The vector space $\mathfrak{H}$ depends on the Chern Class $\mathfrak{c}$, which in turn depends on the choice of a frame $\mathfrak{fr}$ of $\Sigma \times \R$, which is invariant by translation on the $\R$-direction. The Chern Class $\mathfrak{c}$ is integer if and only if the spin structure defined by $\mathfrak{fr}$ can be extended to $L$. 

In this paper we consider an extension   $ \mathfrak{g}^{\spadesuit}$ of the $KS$-algebra, which is associated to $(X,L)$. Namely the generators of $ \mathfrak{g}^{\spadesuit}$ are labeled by pairs $(\gamma, \beta) \in  H_1(\Sigma,\Z) \times   H_2(X, L,\Z)  $  with    the image of $\gamma$ on  $ H_1(L,\Z) $ agreeing with $\partial \beta$. There is a natural action of $ \mathfrak{g}^{\spadesuit}$ on $\mathfrak{H}$.
%$\hat{e}_{\gamma, \beta} $

We consider another version $KS$-algebra $\hat { \mathfrak{g} } $ with $ q^{\frac{1}{2}} = e^{ \frac{g_s}{2}}$ and $ \gamma \in \frac{1}{2} H_1(\Sigma, \Z) $ . This algebra will play an important rule in the study of the boundary conditions, and it will be related with the perturbation of the area zero annulus. We denote its generators with  $ (\hat{e}_{\gamma}^{ann0})_{\gamma}$.
%The vector space $\mathfrak{H}_{  \mathfrak{c}  }$ depends on an element  $\mathfrak{c} \in H_1(L, \Z)$, that we call Chern class. $\mathfrak{c} $ depends on the boundary conditions. 
The element
 $ \hat{e}_{\gamma}^{ann0}$ maps $\mathfrak{H}_{  \mathfrak{c}  }$ on $\mathfrak{H}_{  \mathfrak{c} + Q(\gamma) }$, where $Q(\gamma)$ is the image of $\gamma $ on $H_1(L, \Z)$.

%We define an action of $\hat{ \mathfrak{g} } $ on $\mathfrak{H}$. Denote by $ \hat{e}_{\gamma}^{an}$ the image of $ \hat{e}_{\gamma}$. 
%Up to isotopy, the vector $Z \in \mathfrak{H}$  is determinate by the choice of a boundary condition $\mathcal{B}$ on the perturbation of area zero curves. The frame $\mathbf{fr}$ on $\R \times \Sigma$ should be considered as part of $\mathcal{B}$.  

Assume that $X$ has cylindrical ends , i.e., outside a compact region $X$ can be identified with the symplectization  $ \R^+ \times Y$  of a contact five-manifold $Y$.
In closed Gromov-Witten theory (which is the particular case of open-closed theory when $L = \emptyset$ and $K$ is a closed four chain) a choice of a boundary condition consists on the choice of a perturbation of the moduli space of area zero curves on  $ \R^+ \times Y$ which is invariant under the translation on the $\R$-direction. 
%Note that this makes sense since there is an obvious action of $ \R^+ $ on RHS of (\ref{obstruction-bundle0}), which it is trivial on the second factor. 
%Note that in the closed case the boundary condition affect only the invariant of degree zero.  
%\begin{remark} \label{boundary-closed}
% For example ,  for area zero spheres with three marked points mapped on $K$, it defines the intersection pairing $K \cap K \cap K$, which depends on the choice on the choice of a small perturbation of the diagonal of $Y^3$ transversal to $K_{\infty}^3$. 
%\end{remark}

The generalization to open-closed Gromov-Witten theory of the  boundary condition  is not straightforward. Since in this case we need to work with the moduli of multi-curves $  \overline{\mathcal{M}}_{G,m}$ the  set of area zero components (of a multi-curve)  mapped on $ \R^+ \times Y$ changes when we move inside  the moduli space. It is necessary to set the invariance in the $\R$ direction consistently in a suitable way.   
 Observe that ,  since   there can be  components of the multi-curve of   area zero  in each degree   , the choice of the boundary condition affects  the invariants of not zero degree also . This is in contrast with the closed case, where the ambiguity associated to the choice of the boundary condition affect only the invariants of degree zero. 
% all the components of the moduli of multi-curves of positive area live in a compact region of $X$. However the moduli of the area zero components are not compact. 
%Therefore we need  to fix boundary conditions on the perturbation on the moduli of area zero pseudo-holomorphic curves. 

 We consider  a particular type of boundary conditions which in particular gives  a frame $\mathbf{fr}$ of $\R^+ \times \Sigma $ invariant by translation on the $\R$-direction. 
 %The Chern Class $\mathfrak{c}$ is even if and only if the spin structure defined by $\mathfrak{fr}$ can be extended to $L$. 

If $\mathcal{B}_1,\mathcal{B}_2$ are two such boundary conditions, the associate $MC$-cycles $Z^{\mathcal{B}_1},Z^{\mathcal{B}_2}$ of $\mathfrak{H}$, are related by the action of the quantum KS-algebra: 
\begin{equation} \label{boundary-shift}
Z^{\mathcal{B}_1}= a \hat{e}_{\gamma}^{ann0} Z^{\mathcal{B}_2} 
\end{equation} 
for some $a \in \Q((g_s))$ and $\gamma \in \frac{1}{2} H_1(\Sigma,\Z)$.  $\gamma$ is determined by the  frames associated to $\mathcal{B}_1$ and  $\mathcal{B}_2$.
%The fact that $\gamma$ is half-integer (instead than integer) is of critical importance. 

The scalar ambiguity $a $ of (\ref{boundary-shift}) can be considered the analogous of the one appearing in closed Gromov-Witten theory. For example, in the case of area zero spheres with three marked points, it is  the classical ambiguity appearing in  the definition of $K \cap K \cap K$. 
%(in the closed case $\partial K = \emptyset$).

%--The factor  $\hat{e}_{\gamma}$ in (\ref{boundary-shift}) is more interesting. We shall constrain partially the boundary condition in order to fix $\hat{e}_{\gamma}$ by some topological data. In particular we use a frame of $T_*L|_{\Sigma}$ to refine the definition of four chain. There is a natural action of  $H_1(\Sigma, \Z)$ on such frames.  Even elements do not change the spin structure of $\Sigma$. 

%To partially fix the boundary conditions we refine the four chain with a frame of $T_*L|_{\Sigma}$.

%The group $H_1(\Sigma, \Z)$ acts on the boundary conditions. 
%%%$$Z  \mapsto \hat{e}_{\frac{\gamma_0}{2}} Z$$

%Two different boundary conditions are related by an element $MCH(\Sigma \times \R)$. This can be written as $a \hat{e}_{\gamma}$, for $a \in \Q((g_s))$. 
%$ \hat{e}_{\gamma}$ is determined by the frame of $\Sigma$.

%In this introduction for simplicity we assume that the symplectic form of $M$ is exact, this condition is removed in the paper. The case if exact Calabi-Yau is however important in mirror symmetry (see \cite{Z} for a discussion related to this paper). Under this assumption we can identify the topological charges $\Gamma$ of holomorphic curves with $H_1(L,\Z)$. 
%We shall be interested to configuration of (multi-) curves with topological charge in the lattice 

We also consider compactification frames of $L$, which are defined in terms of  maximal isotropic sub-lattice $F \subset H_1(\Sigma , \Z)$  ( this notion is a straightforward extension of \cite{frame}). Using the multi-linking number defined in \cite{OGW3} we can define a wave function $ \Psi_F(x)$. Alternative this can be defined using the Chern-Simons Propagator (see \cite{BoundaryStates2}).

%For each $\gamma \in \Gamma_{flavor}$, $\mathcal{H}^{coh}(\gamma)$ is  a one dimensional vector space. The choice of an isomorphism $\mathcal{H}^{coh}(\gamma) \cong \Q$ involves the choice of  a \emph{frame} of $L$,  that is maximal isotropic sub-lattice $F \subset H_1(\Sigma , \Z)$  ( this notion is a straightforward extension of \cite{frame}). In particular this defines rational numbers invariants using a variant of linking number between pairs of curve that depends on $F$ . The wave function: $$ \Psi_F(x) = \sum_{\gamma \in \Gamma_{flavor}} \text{Link}_F (Z_{\gamma}) x_{\gamma}.$$
%We stress that the wave function $\Psi_F$ depends on the choice of a frame, instead the vector  $\Psi$ in (\ref{psi}) does \emph{not} depend on it.

\subsection{Quantum shift of the coordinates}

The frame $\mathbf{fr}$ is strictly related to the quantum shift of the coordinated appearing in the physical literature which is critical in order to match the B-model picture where it appears in
the argument of the quantum dilogarithm. Hence, it is critical to consider this extra data in order to have a dictonary between A-model and B-model. We know explain this relation.

The definition of $MCH(M,\gamma| \mathfrak{c})$ depends also on the Chern class $\mathfrak{c} \in  H_1(L, \Z)$
%$ \mathfrak{c} \in \text{Im} \{ H_1(L, \Z) \rightarrow H_1(L, \Q)  \}$ 
which is determinate in a purely topological way from the four chain $K$. In the compact case  $\mathfrak{c}$ is always even, i.e.,
 $ \mathfrak{c} \in 2 H_1(L, \Z) $. In the not compact case this is not always the case, $ \mathfrak{c}$
 is even if the spin structure on $\Sigma$ associate to $\mathbf{fr}$ can be extended to $L$.

The Chern class needs to be even if we want  the wave function be a periodic function of 
%$x \in \text{Im} \{ H_1(\Sigma, \R) \rightarrow H_1(L, \R)  \}$ 
$x \in   \frac{H^1(L, \R)}{H^1(L, \Sigma, \R)}   $:
$$x \mapsto x + 2 \pi H^1(L , \Z) .$$

%$$x \in   \frac{H^1(L, \R)}{H^1_{comp}(L, \R)}   $$

%It is easy to show that the last condition is equivalent to require that the spin structure on $\Sigma$ associate to $\mathbf{fr}$ can be extended to $L$. 
%Hence it is natural to constrain   the frame $\mathbf{fr}$ is compatible with the spin structure of $L$.
 
 It is necessary to fix a spin structure of $L$ also in order to fix an orientation of the moduli space of (multi-)curves.
 If the spin structure is changed by an element $\alpha \in  H^1(L, \Z_2)$, the orientation of the moduli space of (multi-)curves of homology class $\beta$ changes by $(-1)^{   \langle \alpha , \partial \beta \rangle}$. The action on the wave function can be obtained by the shift of  the coordinates  
 $$x \rightsquigarrow x + i \pi \alpha .$$ 
 
If we restrict our attention to the frames $\mathbf{fr}$ compatible with the spin structure of $L$, 
 the change of the spin structure by  $\alpha \in  H^1(L, \Z_2)$ needs to be accompanied by a twist of the frame by an element $\gamma \in H^1(\Sigma , \Z)$ such that
$$\alpha =  \gamma  \text{   mod  }   2  H^1(\Sigma , \Z) .$$
  Hence, if  $\gamma \in Ker (H_1(\Sigma, \Z) \rightarrow H_1(L, \Z) )$, we obtain a change of coordinates 
\begin{equation} \label{shift-coordinates}
x \rightsquigarrow x + (i \pi +  \frac{g_s}{2} )  \gamma  \text{   mod  }   2 \pi H^1(L , \Z)  .
\end{equation}
%where the term $\pi$ take in account the shift of the spin structure.

%The $MC$-cycle $Z = (Z_{\beta})_{\beta}$ obtained from the moduli space of multi-curves depends in particular  by the choice of a spin structure of $L$, since it determines the orientation of the moduli space of curves.

%If we change the change the spin structure 
%$$z_{\beta} \rightsquigarrow z_{\beta} + i \pi$$

%If we want that the wave function has to be periodic on $x$, i.e. , invariant by $x \mapsto x + 2\pi i$, the Chern Class $\mathfrak{ch}$ has to be even. Because of ...,  this condition equivalent to ask that the spin structure on $\Sigma$ associate to $\mathbf{fr}$ can be extended to $L$.  In order to archive this, it is natural to require  that the frame $\mathbf{fr}$ is compatible with the spin structure . 

%Thus a twist of $\mathbf{fr}$  by $\gamma \in H_1(\Sigma , \Z)$ induces  a shift of the the spin structure by the image of $\gamma$ on $H_1(\Sigma , \Z)/ 2 \times H_1(\Sigma , \Z)= H_1(\Sigma ,\frac{1}{2} \Z)$. 

The quantum shift appearing in (\ref{shift-coordinates}) is related to the metaplectic correction  when the wave function is interpreted as a perturbative state for the geometric quantization of abelian flat connections (see \cite{BoundaryStates2}).   

%The shift of the spin structure is necessary in particular when a singularity on the moduli space of Lagrangians is crossed. Since $Ker (Q)$ differs on the two sides, the condition that the spin structure has to be extended on $L$ forces the shift of the spin structure. To compare the wave function on the two sides it will be necessary to perform a  Konstant-Rosemberger type homomorphism, the shift  (\ref{shift-coordinates}) is part of it. This point is discussed more in XXX.  The most simple example arise in the case of Topological Vertex, where this conjecture matches with the Fourier transformation property of the quantum dilogharitm. 

The shift of the spin structure is necessary in particular when a singularity on the moduli space of Lagrangians is crossed (see also Remark \ref{analitic-continuation})  
since  in the process $Ker (H_1(\Sigma, \Z) \rightarrow H_1(L, \Z) )$ changes.
%the condition that the spin structure has to be extended on $L$ forces the shift of the spin structure. 
 The quantum shift  (\ref{shift-coordinates}) can be seen as  part of a  Konstant-Rosemberger type homomorphism mentioned in   Remark \ref{analitic-continuation}.
 % After these considerations are taken in account, we conjecture that the wave-function is the analytical continuation  
  The most simple example of this phenomena  arises in the case of Topological Vertex, where in the B-model side corresponds to the shift  of the quantum dilogharitm appearing in its Fourier Transformation. (See for example section $5.10$ of \cite{ff1}.)

\begin{remark} \label{analitic-continuation}
%The objects defined considered in this paper are associated to a fixed pair $(X,L)$. 
% and are formal power series on $T$ and $x$.
%The considerations we made above hold for a fixed $L$ in the formal power series of the Novikov ring $T$ and $x$. 

It is natural to ask  how the objects defined in this paper change when the Lagrangian manifold $L$ moves on its moduli space,  in particular when we cross a singularity, what in  physical literature is usually referred to as open phase transition.

The boundary conditions can be moved across the singularity without any problem. In some sense the geometry at infinity does not see the singularity. 

After the boundary conditions are fixed, we would like to consider objects which, in a suitable sense, do not change when we cross the singularity, 
%It is trivial to prove that the $MCH$ element defined in this paper does not change as long as we move $L$ on the smooth part of its moduli. 
%It is straightforward to prove that the $MCH$ element associate does not change when wave function this means and prove it in the smooth part of the moduli. 
in accord with the general physical principle stating that the moduli spaces are smooth after the quantum corrections are included . 

To compare the wave function in different regions of the moduli, it is necessary to implement an analytical continuation.
The objects considered in this paper are defined as  formal power series, hence we need to make suitable converge assumptions in order to obtain an actual function and consider analytical continuation (this step may involve to consider not perturbative corrections also).  Moreover  we need to consider the wave function as a perturbative state for the geometric quantization of abelian flat connections (see \cite{BoundaryStates2}). Since $Ker (H_1(\Sigma, \Z) \rightarrow H_1(L, \Z) )$ changes when the singularity is crossed, the polarization changes and it is necessary to implement a    Konstant-Rosemberger type homomorphism in order to compare the wave functions in the two sides.  
%Conjecturaly, the perturbative state defined by the wave function in the two side of the wall  analytical continuation of each other. 
After these considerations are taken in account, we  expect that the wave function is invariant by open phase transition. 

% Moreover it is necessary to take in account the interesting problem of instant corrections in Gromov-Witten theory, which are still scarcely understood. 

This result provide the general setting to understand the quantization of moduli of $A$-branes (see also remark \ref{A-branes-remark}) .

\end{remark}

\subsection{Non Abelian Case} The $MCH$ homology we consider in this paper are associated to decorated graphs space of  fixed Euler characteristic, which is reffered  as the abelian case. We can use this construction in order to consider curves of a fixed genus and number of boundary components $(g,h)$ using the following standard procedure. 

Identify the tubular neighbourhood of $L$ with a neighbourhood of the zero section of $T^*L$. Consider a $k$-brunched covering $\tilde{L}$ of $L $ in $T^*L \rightarrow L$. The $MC$-cycle of fixed euler characteristic associated to $\tilde{L}$ is related to the moduli of fixed $(g,h)$ of $L$. If we apply this for arbitrary $k$ we recover all the informations about the curves of fixed $(g,h)$.

Alternatively it is not hard to extend the definition of $MCH$ of \cite{OGW3} to consider decorated graphs of fixed genus and boundary component, providing a more direct way to determine the refined generating function.

According the proposal of Aganagic and Vafa (\cite{AV}),  this result can be applied in order to obtain  a mathematical theory of HOMFLY polynomials of a link $\mathcal{L}$ in $S^3$. This apply in the particular case that $X$ is the conifold and $L$ is the conifold transition of the normal bundle of $\mathcal{L}$ in $T^*S^3$. Actually in this case we do not need to consider brunched covering of $L$ but it is enough to consider the conifold transition of the normal bundle of parallel copies of  $\mathcal{L}$, i.e., $\tilde{L}$ is a disjoint union of copies of $L$.

\section{Multi Curve Cycles}

Fix an oriented three manifold $M$ and a class $\gamma \in H_1(M, \Z)$. 
Consider the  objects  
\begin{equation} \label{generator-nice}
( H , ( w_h  )_{h \in H} ) 
\end{equation}
where $H$ is a finite set , $\{ w_h  \}_{h \in H} $ are closed one  dimensional  chains on $M$ on the homology class $\gamma$  which are close on the $C^0$-topology. 
Denote by $\mathfrak{Gen}(\gamma)^{\dagger}$ the set of the objects (\ref{generator-nice}) modulo the obvious equivalence relation:  $( H , ( w_h  )_{h \in H} ) \cong ( H' , ( w_h'  )_{h \in H'} )$ if there exists an identification of sets $H =H'$ such that $w_h=w_h'$. 

Define $\mathfrak{Gen}(\gamma) \subset \mathfrak{Gen}(\gamma)^{\dagger}$ as the subset obtained imposing the extra condition  
\begin{equation} \label{generator-transversality}
w_h \cap w_{h'} = \emptyset \text{    if    } h \neq h'.
\end{equation}
The vector space $\mathcal{Z}_{\gamma}$ of $\mathbf{nice \thickspace MC-cycles}$ in the homology class $\gamma$  is the formal vector space generated by $\mathfrak{Gen}(\gamma)$.

In this paper we consider $\mathcal{Z}_{\gamma}[[g_s]]$ the formal power series on the formal variable $g_s$ with coefficients  $\mathcal{Z}_{\gamma}$. 

%The vector space $\mathcal{Z}_{\gamma}$ of $\mathbf{nice \thickspace MC-cycles}$ in the homology class $\gamma$  is the space of formal power series on $g_s$ with coefficients the vector space generated by  $\mathfrak{Gen}(\gamma) $.

%\begin{remark}
%In \cite{OGW3} we defined 
%\end{remark}
%\subsection{isotopies}

%As in \cite{OGW3}, an isotopy $\tilde{Z}$ defines a one parameter family of element $\mathcal{Z}_{\gamma}[[g_s]]$. In this paper we use a slightly different notation. 

An isotopy of $MC$-cycles is defined by a formal power series 
$$\tilde{Z}= \sum_i  g_s^{k_i}   r_i ( H_i , ( \tilde{w}_{h,i}  )_{h \in H_i} , [a_i,b_i] ) ,$$
for some $r_i \in \Q$, $k_i \rightarrow \infty$.
In particular $\tilde{Z}$ defines a one parameter family  of elements of $Z_t \in \mathcal{Z}_{\gamma}$, $t \in [0,1]$, which can be discontinues on a finite number of times (we do not define $Z_t$ if $t$ is a discontinuity point). 
if for $h,h' \in H_i$ , $ \tilde{w}_{h_1,i},  \tilde{w}_{h_2,i}$ cross transversely at a time $t_0$, we require that $\tilde{Z}$ jumps according to the formula
\begin{equation} \label{jump}
 Z^{t_0^+} - Z^{t_0^-} =  \pm g_s^{k_i+1} r_i   (H_i , ( w_{h,i}^{t_0}  )_{h \in H_i \setminus \{h_1,h_2 \} }  )  .
\end{equation}
where the sign is defined by the sign of the crossing.
Denote by  $\tilde{\mathcal{Z}}_{\gamma}$ the set of isotopies of nice $MC$-cycles.

We need also to introduce the space of $MC$-cycles $\mathcal{Z}_{\gamma| \mathfrak{c}}[[g_s]]$ with not vanishing Chern Class $\mathfrak{c}$. This is done as in \cite{OGW3}.

\subsection{Graph version}

Fix $\epsilon>0$ small enough.
For $\gamma \in  H_1(M,\Z)$ define the set of pair of one dimensional closed chains
$$  \mathcal{P}_{\gamma} = \{      (w,w' ) |  [w]=[w' ] = \gamma ,  w \cap  w'=\emptyset , d(w,w') < \epsilon \}/\sim$$
where $\sim$ is the equivalence relation generated by small isotopies of the pair of curves that do not cross each other . For $\epsilon$ small enough $\mathcal{P}_{\gamma} $ does not depend on $\epsilon$.

We should think to an element of   $\mathcal{P}_{\gamma}$ as (a small perturbation of) a framed link in the homology class $\gamma$, i.e.,  $w'$ is a small translation of $w$ in the direction of the frame.

There is a $\Z$-action on $\mathcal{P}_{\gamma}$
$$ \mathbf{tw}_k : \mathcal{P}_{\gamma} \rightarrow \mathcal{P}_{\gamma} , \text{       } k \in \Z .$$
The pair $\mathbf{tw}_k(w,w' )$ is obtained  crossing $k$ times (with sign)   the pair $(w,w')$.

Let $\mathfrak{Gen}(\gamma)^{\sharp}$ the set of oBjects 
\begin{equation} \label{generator-graphs}
 ( E^{in}, E^{ex}, ( (w_e,w_e') )_{e \in E^{in}} , (  w_e )_{e \in E^{ex}} ) . 
 \end{equation}
 modulo the obvious equivalence relation.
Let $\mathcal{Z}^{\sharp}$ the vector space of formal linear combinations with coefficients in $\Q[[g_s]]$.

We can introduce the notion of isotopy in a analogous way as for XXX. We obtain a one parameter family of elements of elements of $\mathcal{Z}^{\sharp}$, with the jumping condition:
if  for an addendum (\ref{generator-graphs}) with coefficient $g_s^k r$ , for some $e_0 \in E^{in}$, $w_{e_0}$ and $w_{e_0}'$ intersects transversally,  $Z_t$ jumps of 
$$  g_s^{k} r ( E^{in} \setminus \{ e_0 \}, E^{ex},  ((w_e,w_e') )_{e \in E^{in} \setminus \{ e_0\}} , (  w_e )_{e \in E^{ex}} ) .  $$

There is a map 
$$  \mathcal{Z} \rightarrow \mathcal{Z}^{\sharp} $$
defined by
$$ (H ,  ( w_h )_h ) \mapsto \sum_{E} \frac{g_s^{E^{in}}}{Aut(E)} ( E^{in}, E^{ex}, ( (w_e,w_e') )_{e \in E^{in}} , (  w_e )_{e \in E^{ex}} ) $$
where $ (w_e,w_e')  = (w_h,w_{h'})$ for $e = \{ h,h'\} \in E^{in}$, $w_e = w_h$ for $e = \{ h \} \in E^{ex}$.

\subsection{Open Gromov-Witten $MC$-cycle}

Let $(X,L)$ be a pair given by a Calabi-Yau simplectic six-manifold $X$ and a Maslov index zero lagrangian submanifold $L$. 
%As in \cite{OGW1}, to deal with the nodes of type $E$ in sense of \cite{L}, 
We assume $[L] =0 \in H_3(X, \Z)$. 
%Fix a $[K] \in H_4(X,L)$ with $\partial [K] = [L]$. 
Fix  a four chain $K$ with $\partial K =L$.

To the four chain $K$ it is associated a Chern Class $\mathfrak{c}(K) \in H_1(L,\Z)$.

\begin{theorem} \label{main-theorem}
Let $\beta \in H_2(X,L, \Z)$. To the moduli space of pseudoholomotphic multi-curves of homology class $\beta$  it is associated a multi-curve cycle $Z_{\beta} \in \mathcal{Z}_{\partial \beta| \mathfrak{c}} [[g_s]]$ with Chern class $\mathfrak{c}$. $Z_{\beta}$ depends by the varies choices we made to define the  Kuranishi structure and its perturbation on the moduli space of multi-curves.
%( the compatible almost complex structure $J$, the representative $K$,  the various choices we made to define the  Kuranishi structure and its perturbation on the moduli space of curves). 
%Moreover  $Z_{\beta, \chi}$ does not change if add to $K$ a four-chain in $M$ homologically trivial.
Different choices lead to isotopic $MC$-cycles.
\end{theorem}

The result is an extension of the main theorem of \cite{OGW3}. In  \cite{OGW3} we considered $MC$-cycles of a fixed Euler Characteristic. In particular only a finite number of decorated graphs are involved in the construction. This point is critical  in order to  construct a perturbation of the moduli space of multi-curve. 
The extension to a full $MC$-cycle is made using the argument subsection \ref{extension-section}. Hence, using  the $MC$-cycles constructed  in \cite{OGW3} (which involves only a finite number of graphs) we can construct a full $MC$-cycles.

%--However  defined $MC$-cycles $Z^k \in \mathcal{Z}^k$ for each $k$ such that there exists  $\tilde{Z}^k$  isotopy between $Z_{k+1}$ and $Z_k$  for each $k$. 

\subsubsection{Review of the construction of \cite{OGW3}} \label{review-section}

In this subsection we review the construction of the Gromov-Witten $MC$-cycle of  \cite{OGW3}.

A Kuranishi structure on the moduli space of pseudo-holomorphic curves  in particular associate to
 each pseudo-holomorphic curve $\mathbf{p}=(\Sigma,u)$ an  obstruction bundle 
\begin{equation} \label{obstruction-bundle}
E_{\mathbf{p}} \subset   C^{\infty} (\Sigma, u^*(TX) \otimes \Lambda^{0,1}(\Sigma)).
\end{equation}
to each pseudo-holomorphic curve $\mathbf{p}=(\Sigma,u)$. 
%satisfying suitable conditions. 

The standard procedure to construct a Kuranishi structure on moduli space of stable pseudo-homorphic curves uses an argument by induction on the singular strata. For each strata, it is used a  gluing argument to extend the Kuranishi structure to a neighborhood of the substratas and after the Kuranishi structure is extended inside the strata.  The method does not lead to a unique Kuranishi structure, however it determinate it up to isotopy, i.e., given two Kuranishi structures constructed in this way there exists a Kuranishi structure on the moduli spaces times $[0,1]$ whose restrictions to $0$ and $1$  agree with the starting Kuranishi structures.
 This isotopy is constructed using an analogous inductive argument on the strata.

Denote by $  \overline{\mathcal{M}}_{(g,V, D, (H_v)_v)}(\beta) $ 
the moduli  of curves of genus $g$ , homology class $\beta$, whose boundary components are labelled by $V$, internal marked points are  labelled by $D$ and boundary  marked points are labelled by $(H_v)_v$.

In \cite{OGW3}, to a decorated graph $(G,m) \in  \mathfrak{G}_l$, it is associated the moduli space of multi-curves
\begin{equation}    \label{multi-curve-interp0}
 \overline{\mathcal{M}}_{G,m} :=    \left( \prod_{c \in Comp(G)} \overline{\mathcal{M}}_{g_c,D_c, V_c, (H_v)_{v \in V_c} }(\beta_c)  \right)  / \text{Aut}(G,m) \times \Delta^l
 \end{equation}   
where $\Delta^l$ is the standard simplex  of dimension $l$.

%The moduli space $ \overline{\mathcal{M}}_{G,m} $ is endowed with a Kuranishi structure constructed using an inductive argument similar to the one used in  the construction of the Kuranishi structures for stable pseudo-holomorphic curves. 

Using an inductive argument similar to the one used in  the construction of the Kuranishi structures for stable pseudo-holomorphic curves, we can obtain a collection of Kuranishi structures $(\overline{\mathcal{M}}_{G, m })_{(G,m) \in \mathfrak{G}(\beta, \leq \kappa   )}$  such that:
\begin{itemize}
\item  forgetful compatibily holds;
% with respect each puncture associate to all the external edges $E^{ex}(G)$;
\item the evaluation map 
$$\text{ev} : \overline{\mathcal{M}}_{G,m} \rightarrow L^{H(G) \setminus H_l} \times X^{D(G)}$$ 
is weakly submersive;
\item   the corner faces of $\overline{\mathcal{M}}_{(G,m)}$ are in bijection with  the graphs $\mathfrak{G}(G,m)$
$$  \mathfrak{G}(G,m ) \leftrightarrow  \{    \text{   corner faces of  } \overline{\mathcal{M}}_{(G,m)}    \} .$$
The corner face $ \overline{\mathcal{M}}_{G,m }( G',m',E') $ corresponding to  $(G',m',E') \in \mathfrak{G}(G,m)$ comes with an identification of Kuranishi spaces: 
\begin{equation} \label{corner-faces-multi-interp0}
 \overline{\mathcal{M}}_{G, m }(G' ,m', E') \cong   \delta_{E'}  \overline{\mathcal{M}}_{G',m'}.
  % \overline{\mathcal{M}}_{G, m }(G' ,E', F') \cong   \left( \prod_{c \in C(G)} \overline{\mathcal{M}}_{(g_c,h_c), H_c, m_c}(\beta_c) \times_{L^{H(G)} \times X^{n_G} \times \partial_F \Delta^{C(G)}} ( L^{E(G)} \times L^{n_G} \times \partial_F \Delta ) \right) / \text{Aut}(G)
\end{equation}
%where  $m'$ is the labeling of $E(G')$ corresponding to $\partial_{F'} m$.
\end{itemize}
The Kuranishi spaces $\delta_{E'}  \overline{\mathcal{M}}_{G',m'}$ are defined as fiber products. See \cite{OGW3}.
%Note that in point $(3)$ it is used point $(2)$ to define the Kuranishi structure on $\delta_{E'}  \overline{\mathcal{M}}_{G',m'}$ as fiber product. 
%in (\ref{bedge}).

%For each pseudo-holomorphic curve $(\Sigma,u)$ the obstruction bundle is a sub-vector space
%$$  \prod_{c \in Comp(G)} C^{\infty} (\Sigma, u^*(TM) \otimes \Lambda^{0,1}(\Sigma))$$

%For $v \in V^{deg}(G)$ define 
%$$\delta_v \overline{\mathcal{M}}_{G,m} = \text{ev}_v^{-1} ( L) .$$

Consider collections of perturbations
$( \mathfrak{s}_{(G,m)}^+ )_{(G,m) \in \mathfrak{G}(\beta, \leq \kappa   ) }$ 
of the collection of Kuranishi spaces
 $(\overline{\mathcal{M}}_{G,m} )_{G \in \mathfrak{G}(\beta, \leq \kappa   ) }$ such that 
%For each $(G,m) \in \mathfrak{G}(\beta)$ with $\chi(G) \geq \chi$, there exist a pertubation $\mathfrak{s}_{(G,m)}$ of the Kuranishi map $s$ such that 
%collection of  pertubations  $\{ \mathfrak{s}_{(G,m)} \}_{(G,m) \in  \cup_{\chi' \geq \chi}  \mathfrak{G}(\beta, \chi'), }$ such that 
\begin{itemize}
\item they are transversal to the zero section and small in order for constructions of virtual class to work;
\item they are compatible with the identification of Kuranishi spaces (\ref{corner-faces-multi-interp0});
\item 
%are forgetful compatible with respect to the external edges and $E_0$.
  forgetful compatibility holds;
 \item  they are transversal  when restricted to $\delta_{E}  \overline{\mathcal{M}}_{G,m}$ for each $E \subset (E^{in}(G) \setminus E_l) \sqcup D(G)$ 
 \end{itemize}

Using these perturbations we define the collections of chains $( Z_{(G,m)}^+ )_{(G,m) \in \mathfrak{G}_*(\beta, \leq \kappa   ) }$
\begin{equation} \label{multi-curve-cycle-not-ab}
 Z_{(G,m)}^{not-ab, +} =  (\text{ev}_{G,{m}})_* (\mathfrak{s}_{G,{m}}^{-1}(0)) \in C_*(L^{H(G)} \times X^{D(G) })^{Aut(G,m)}  .
\end{equation}

The not abelian $MC$-cycle $Z^{not-ab}$ is obtained taking the fiber product
\begin{equation} \label{intersect-K-gen}
Z_{G,m}^{not-ab}=  Z_{G,m}^{not-ab,+} \times_{X^{D(G^+) }} K^{D(G^+)}  .
\end{equation}

From $Z^{not-ab}$ we can construct the $MC$-cycle and the nice $MC$-cycle as in \cite{OGW3}.

%--

%We define $\delta_{E}  \overline{\mathcal{M}}_{G,m}$ for each $E \subset (E^{in}(G) \setminus E_l) \sqcup D(G)$  using  ADATTARE
%$$  \delta_{E}  \overline{\mathcal{M}}_{G,m} =  \overline{\mathcal{M}}_{G,m} \times_{X^{E \cap D(G)} } K^{E \cap D(G)} $$

%$$  \delta_{E, D}  \overline{\mathcal{M}}_{G,m} =  \overline{\mathcal{M}}_{G,m} \times_{ L^E \times  X^D }( Diag_E \times K^{D} )$$

\begin{remark}
    In \cite{OGW3} was considered also  the Kuranishi spaces
 \begin{equation}    \label{multi-curve-interp}
\overline{\mathcal{M}}_{G,m}^K  :=   \overline{\mathcal{M}}_{G,m} \times_{ X^{D(G)} }   K^{D(G)} .
\end{equation} 
We define  collections of perturbations 
$( \mathfrak{s}_{(G,m)} )_{G \in \mathfrak{G}(\beta, \leq \kappa   ) }$ 
of the collection of Kuranishi spaces
$(\overline{\mathcal{M}}_{G,m}^K )_{G \in \mathfrak{G}(\beta, \leq \kappa   ) }$
which satisfy suitable transversality conditions,  compatibility on the corner faces and forgetful compatibility. The associated virtual fundamental chain leads to a $MC$-cycle $Z= (Z_{G,m})_{G \in \mathfrak{G}(\beta, \leq \kappa   )}$. 
\end{remark}

\subsection{Coherent Cycles}

For $(w,w')  \in   \mathcal{P}_{\gamma}$, define $(w,w')^k \in Gen(\gamma)$ as
$$(w,w')^k= ( w_1 , w_2, ..., w_k )$$
where $w_i$ are closed one-chains, with $(w_i, w_j) \in  \mathcal{P}_{\gamma}$  for each $i,j$, and $(w_i, w_j)= (w, w')$ as element of $ \mathcal{P}_{\gamma}$.  $(w,w')^k$ is well defined up to small isotopy.

Let $\mathcal{Z}^{\diamond, \diamond}$ be the set of linear combinations of 
$$\text{exp}(\frac{1}{2} (w, w')) = \sum_k \frac{1}{2^k k!}  (w,w')^k $$

\begin{lemma}
$\mathcal{Z}^{\diamond, \diamond}_{\gamma}/isotopies$ is a rank one modulo over $\Q[[g_s]]$.
\end{lemma}

Given two coherent cycles with disjoint support we can define the product
\begin{equation} \label{product-coherent}
\text{exp}(\frac{1}{2} (w_1, w_1'))  \text{exp}(\frac{1}{2}(w_2, w_2'))=  \text{exp}(\frac{1}{2} (w_1+w_2, w_1'+ w_2')) .
\end{equation}

\begin{proposition} \label{coh-construction}
 To each $\beta \in H_2(X,L)$ it is associated an element $Z_{\beta} \in   \mathcal{Z}_{\beta}^{\diamond, \diamond}$ up to isotopy, in the same sense of Theorem \ref{main-theorem}.
\end{proposition}
\begin{proof}
Above we have only considered the notion of coherent nice $MC$-cycle. To adapt the proof of Theorem \ref{main-theorem} to coherent cycles it is necessary to  extend the notion of coherent cycle to  (not-nice) $MC$-cycles. For this we need to adapt the definition of forgetful compatibility.  Recall that the forget compatibility of \cite{OGW3} is defined in terms of chains with coefficients on $\mathfrak{Gen}(\gamma)$. We require that  these chains have coefficients on nice  coherent $MC$-cycles. 

%We remark that from the construction of \cite{OGW3} we obtain a coherent $MC$-cycle.  
We remark that the $MC$-cycle constructed  in \cite{OGW3}  from the moduli space of (multi-) pseudo-holomorphic curves is coherent. Actually it satisfies a stronger property: $w_h=w_{h'}$ for each $h,h'$. 

Finally  the correspondence (up to isotopy) $ \mathcal{Z} \rightsquigarrow \mathcal{Z}^{nice}$ can be proved for coherent $MC$-cycles using the same argument.
\end{proof}

\subsection{MC-cycles including the degrees} \label{nice-MC2-section}

In \cite{OGW3} we also considered nice $MC$-cycles which keep trace of the components of multi-curves. A component will corresponds to an element of $V$ in (\ref{generator2}) below. Here we consider a slightly different version, which forget about the Euler Charactertic of each single component and as before we consider arbitrary total Euler Characteristic at the same time.

The nice $MC$-cycles we define below arise from moduli space of multi-curves associated to graphs without closed components.
(A closed component of a decorated graph $G \in \mathfrak{G}$ considered in \cite{OGW3} is a component $c \in Comp(G)$ with $V_c = D_c = \emptyset$.) That is,  we have discarded  the purely closed contributions from the $MC$-cycle. 
%The closed components are the generators of the Closed Gromov-Witten partition function. 
%Open Gromov-Witten partition function is obtained quotient the Open-Closed partition function by the closed one. This corresponds to consider graphs without  closed  components. 

We fix the following data:
 \begin{itemize}
\item An oriented three manifold $M$,
\item a finite-rank free abelian group $\Gamma$, called \emph{topological charges},
\item an homorphism of abelian groups
$$ \partial : \Gamma \rightarrow H_1(M, \Z) $$
called boundary homomorphism.
\item an homomorphism of abelian groups 
$$ \omega : \Gamma \rightarrow \R $$
called \emph{symplectic area};
\item A one dimensional homology class $\mathfrak{c} \in H_1(M, \Z)$ called Chern-Class.
\end{itemize}

In the geometric context associated to $(X,L)$ we have 
\begin{itemize}
\item $\Gamma= H_2(M,L,\Z)$
\item $\partial: \Gamma \rightarrow H_1(L,\Z)$
is the usual boundary in homology.
\item $\omega(\beta) = \int_{\beta} \omega$
for $\beta \in \Gamma$.
\item $\mathfrak{c} = \mathfrak{c}(K) \in H_1(L, \Z)$ is defined from the four chain $K$ (see \cite{OGW3}).
\end{itemize}

A generator of the $MC$-cycles is defined by an array
\begin{equation} \label{generator2}
(V, H,  ( w_{v,h} )_{v \in V,h \in H} ,  ( w^{ann}_h )_{h \in H}, (  \beta_v )_{v \in V}  )
\end{equation}
where $V$ and $H$ are finite sets, $\beta_v \in \Gamma$,   $ \{ w_{v,h} \}_{h \in H}$ are closed integer one dimensional chain on $M$ on the homology class $\partial \beta_v$ , which are close in the $C^0$-topology.  
 $ \{ w^{ann}_h \}_{h \in H}$  are closed half-integer one dimensional chains on $M$  on the homology class $\mathfrak{c}$,   which are close in the $C^0$-topology.  

We assume that 
\begin{equation} \label{support-property-nice-cycles}
   \beta_v \notin tors(\Gamma) , \quad  \parallel \beta_v \parallel \leq C^{supp}  \omega(\beta_v) \text{   for each  } v \in V(G)
\end{equation}

Denote by $\mathfrak{Gen}(\beta)^{\dagger}$ the set of objects (\ref{generator2}) such that $\sum \beta_v = \beta$, modulo the obvious equivalence relation.   

Observe that (\ref{support-property-nice-cycles}) implies that there exists  $N_{\beta } \in \Z_{> 0}$ such that
\begin{equation} \label{bound-vertices}
    |V| \leq N_{\beta } 
\end{equation}
for each generator (\ref{generator2}).

Set $w_h= \sum w_{v,h} + w^{ann0}_h$.
Define $\mathfrak{Gen}(\beta) $ as the subset of $ \mathfrak{Gen}(\beta)^{\dagger}$ obtained imposing the extra condition  
\begin{equation} \label{generator-transversality2}
w_{h} \cap w_{h'} = \emptyset \text{    if    } h \neq h'.
\end{equation}
%The vector space $\mathcal{Z}_{\beta}^{\diamond}$ of $\mathbf{nice \thickspace MC-cycles}$ in the homology class $\beta$  is the formal vector space generated by $\mathfrak{Gen}(\beta)$.

Let $\mathcal{Z}^{\blacklozenge}_{\beta} $ be the set of formal power series
$$ \sum_i g_s^{k_i} Gen_i$$
with $Gen_i \in \mathfrak{Gen}(\beta) $, $k_i \rightarrow \infty$, and $k_i + |V_i| \geq 0$. The last condition is related to the fact we have discarded the closed components from the $MC$-cycle.

%\begin{remark}
% OPEN CLOSED VERSUS OPEN
%\end{remark}

To define  isotopies of nice-$MC$ cycles $\tilde{Z} $  we consider objects 
\begin{equation} \label{generator2-iso}
(V, H,  ( \tilde{w}_{v,h} )_{v \in V,h \in H} ,  ( \tilde{w}_{h}^{ann} )_{h \in H}   , (  \beta_v )_{v \in V} , [a,b]   )
\end{equation}
where $[a,b] \subset [0,1]$,  $ \tilde{w}_{h,v} : F_v \times [a,b] \rightarrow M$ for some one dimensional compact manifold $F_v$. 
%Set $ \tilde{w}_{h}^{\flat} = \sum_v \tilde{w}_{h,v} $.

Let $\tilde{\mathfrak{Gen}}(\beta)$ be the set of objects (\ref{generator2-iso}) modulo the obvious equivalence relation.

An isotopy of of nice-$MC$ cycles $\tilde{Z} $ is defined by formal power series 
$$\tilde{Z}= \sum_i r_i g_s^{k_i} (V_i, H_i,  ( \tilde{w}_{v,h,i} )_{v \in V_i,h \in H_i} ,  ( \tilde{w}_{h,i}^{ann} )_{h \in H_i}   , (  \beta_v )_{v \in V_i} , [a_i,b_i]   ) ,$$
with $r_i \in \Q$, $k_i \rightarrow \infty$, and $k_i + |V_i| \geq 0$. 

The formal power series are considered modulo gluing of objects (\ref{generator2-iso}).

$\tilde{Z}$ defines a one parameter family of nice $MC$-cycles $(Z^t)_t$ $Z^t \in \mathcal{Z}^{\blacklozenge}_{\beta}$ discontinuous for a finite number of times. The discontinuity at the time $t_0$ is obtained by the formula:
\begin{equation} \label{jump2}
 Z^{t_0^+} - Z^{t_0^-} =  \pm r' g_s^{k'} (V', H',  ( w_{v,h}' )_{v \in V',h \in H'} ,( {w_{h}^{ann0}}' )_{h \in H'}, ( \beta_v' )_{v \in V'}  )
  % \{ \tilde{w}_h(t_0) \}_{h \in H \setminus {h_1, h_2}}   .
\end{equation}
where
\begin{itemize}
    \item If  $\tilde{w}_{v_1, h_1,i}$ crosses $\tilde{w}_{v_2, h_2,i}$  transversely at the time $t_0$ we have two cases:
    \begin{itemize}
        \item $v_1 \neq v_2$: $H'= H_i \setminus \{ h_1,h_2 \} $  , $V' = V_i \setminus \{v_1, v_2 \}  \sqcup v_0$, where $v_0$ is a new vertex with associated data $\beta_{v_0}'= \beta_{v_1} + \beta_{v_2}$,  $ w_{ v_0,h}' = w_{ v_1,h,i}^{t_0^- } +  w_{ v_2,h,i}^{t_0^- } $ for each $h \neq h_1, h_2$. All the other data remain  the same.
\item  $v_1=v_2$:  $H'= H_i \setminus \{ h_1,h_2\} $, $V'=V_i$  .  All the other data remain the same.
    \end{itemize}
    \item If $\tilde{w}_{v, h_1,i}$ crosses $\tilde{w}^{ann}_{ h_2.i}$, $H'= H_i \setminus \{ h_1,h_2\} $, $V'=V_i$  .  All the other data remain the same.
\end{itemize}
$r' = r_i$, $k' = k_i+1$ and the sign is defined by the sign of the crossing.

\subsection{Extension to a full $MC$-cycle} \label{extension-section}

Assume that we have a filtration on the set of decorated graphs $\mathfrak{G}^0 \subset \mathfrak{G}^1 \subset ....$.
%such that each  $ \mathfrak{G}^k$ is a finite set and it is closed by the operation of forgetting, $\delta_e$,... 
Assume that
\begin{itemize}
    \item if $G \in \mathfrak{G}^k$ then $forget_{e}G \in \mathfrak{G}^k$ for each $e \in E^{ex}(G)$;
    \item if $\delta_e G \in \mathfrak{G}^k$ for $e \in E^{in}(G)$, then $G \in \mathfrak{G}^k$.
\end{itemize}

Under this assumptions we can make the truncation of  $MC$-chain complex to $\mathfrak{G}^k$. We denote by  $\mathcal{Z}_k$ the corresponding space of $MC$-cycles.

We have the following two fundamental Lemmas:

 \begin{lemma} \label{extension-isotopies-inductive}
 Let $Z_{ k}^0 \in  \mathcal{Z}_{ k} $,  $Z_{ k}^1 \in  \mathcal{Z}_{ k} $ and    $\tilde{Z} _{ k} \in  \tilde{ \mathcal{Z}}_{ k} $  isotopy between $Z_{ k}^0$ and $Z_{ k}^1$.
 Let $Z_{ k +1}^0    \in  \mathcal{Z}_{ k +1} $ extending $Z_{ k}^0 $.
 
%--- with $\tilde{Z}_{\prec G_0}^{< -T} = Z_{ \prec G_0} \times (-\infty, -T)  \text{  for   }T >>0 .$

There exist
$Z_{ k+1}^1 \in  \mathcal{Z}_{ k+1} $ extending $Z_k $ and 
 $\tilde{Z} _{k+1} \in  \tilde{ \mathcal{Z}}_{k+1} $  isotopy between $Z_{ k+1}^0$ and $Z_{k+1}^1  $ extending $\tilde{Z} _{ k}$.
%\begin{itemize} 
%\item $Z_{\preccurlyeq G_0}^1  $ extends  $Z_{\prec G_0}^1  $;
%\item  $\tilde{Z} _{\preccurlyeq G_0}  $ extends $\tilde{Z} _{\prec G_0}  $. 
% \end{itemize}
 \end{lemma}

%\begin{lemma}

%$\tilde{Z} _{\prec G_0}^0 \in  \tilde{ \mathcal{Z}}_{\prec G_0} $, $\tilde{Z} _{\prec G_0}^1 \in  \tilde{ \mathcal{Z}}_{\prec G_0} $, $\tilde{Z} _{\prec G_0}^{01} \in  \tilde{ \mathcal{Z}}_{\prec G_0} $

%Let 
%$\tilde{\tilde{Z}} _{\prec G_0}^0 \in  \tilde{\tilde{ \mathcal{Z}}}_{\prec G_0} $ be an isotopy of isotopies that agrees with $\tilde{Z}_{\prec G_0}^0$ for $s=0$,   $\tilde{Z} _{\prec G_0}^1$ for $s=1$, $\tilde{Z} _{\prec G_0}^{01} \in  \tilde{ \mathcal{Z}}_{\prec G_0} $ for $t=0$. 

%There exists $\tilde{\tilde{Z}}_{\preccurlyeq G_0}^0 \in  \tilde{\tilde{ \mathcal{Z}}}_{\preccurlyeq G_0} $ extending $\tilde{\tilde{Z}} _{\prec G_0}^0 $ so that its restriction to $s=0,1$ $t=0$, coincides with $\tilde{Z} _{\preccurlyeq G_0}^0$ , $\tilde{Z} _{\preccurlyeq G_0}^1$, $\tilde{Z} _{\preccurlyeq G_0}^{01}$.
% \end{lemma} 

%An isotopy of isotopy $\tilde{\tilde{Z}}$ is collection of chains $\tilde{\tilde{Z}}_{G,m} \in C_*(M^{H(G)} \times \R_t \times \R_s )$.

 \begin{lemma} \label{extension-isotopies-isotopies}

%Let
%\begin{itemize}
%\item $\tilde{Z} _{k+1}^{\bullet, 0} \in  \tilde{ \mathcal{Z}}_{k+1} $ isotopy between  $Z_{k+1}^{0,0}$ and   $Z _{k+1}^{1,0}$
%\item $\tilde{Z} _{k+1}^{\bullet, 1} \in  \tilde{ \mathcal{Z}}_{k+1} $ isotopy between  $Z_{k+1}^{0,1}$ and   $Z _{k+1}^{1,1}$
%\item $\tilde{Z} _{k+1}^{ 0, \bullet} \in  \tilde{ \mathcal{Z}}_{k+1} $ isotopy between  $Z_{k+1}^{0,0}$ and   $Z _{k+1}^{0,1}$
%\item $\tilde{Z} _{k}^{1, \bullet } \in  \tilde{ \mathcal{Z}}_{k} $ isotopy between  $Z_{k}^{1,0}$ and   $Z _{k}^{1,1}$
%\end{itemize}

Let 
$\tilde{\tilde{Z}} _{k} \in  \tilde{\tilde{ \mathcal{Z}}}_{k} $ be an isotopy of isotopies such that
\begin{itemize}
\item $\tilde{\tilde{Z}} _{k} \cap \{s <-S \} =  Z_{ k}^{0 \bullet} \times \R_{s<-S}  \text{  for   }S \gg 0$
\item  $\tilde{\tilde{Z}} _{k} \cap  \{t <-T \} =  Z_{ k}^{ \bullet 0 } \times \R_{s>S}  \text{  for   }T \gg 0$
\item $\tilde{\tilde{Z}} _{k} \cap \{ t>T \} =  Z_{ k}^{ \bullet, 1} \times  \R_{t >T}  \text{  for   }T \gg 0$
\item $\tilde{\tilde{Z}} _{k} \cap \{ t>T\} =  Z_{ k}^{1,  \bullet} \times  \R_{t >T}  \text{  for   }T \gg 0$
\end{itemize}

Let  $Z_{ k+1}^{0, \bullet} , Z_{ k+1}^{ \bullet, 0 }, Z_{ k+1}^{ \bullet, 1} \in \tilde{ \mathcal{Z}}_{k+1} $ extending  $Z_{ k}^{0 \bullet} ,Z_{ k}^{ \bullet 0 }, Z_{ k}^{ \bullet, 1} $.

There exists $Z_{ k+1}^{1, \bullet}   \in  \tilde{ \mathcal{Z}}_{k+1} $ extending $Z_{ k}^{1, \bullet} $ and
$\tilde{\tilde{Z}}_{k+1} \in  \tilde{\tilde{ \mathcal{Z}}}_{k+1} $ extending 
$\tilde{\tilde{Z}} _{k} $ such that 
\begin{itemize}
\item $\tilde{\tilde{Z}} _{k+1} \cap \{s <-S \} =  Z_{ k+1}^{ \bullet, 0} \times  \R_{s<-S }  \text{  for   }S \gg 0$
\item  $\tilde{\tilde{Z}} _{k+1} \cap \{s >S \}=  Z_{ k+1}^{ \bullet, 1 } \times  \R_{s>S } \text{  for   }T \gg 0$
\item $\tilde{\tilde{Z}} _{k+1} \cap \{t <-T \} =  Z_{ k+1}^{0, \bullet} \times \R_{t>T } \text{  for   }T \gg 0$
\item $\tilde{\tilde{Z}} _{k+1} \cap \{ t>T\} =  Z_{ k}^{1,  \bullet} \times  \R_{t >T}  \text{  for   }T \gg 0$
\end{itemize}

%its restriction to $s=0,1$ $t=0$, coincides with  $\tilde{Z} _{\preccurlyeq G_0}^0$ , $\tilde{Z} _{\preccurlyeq G_0}^1$, $\tilde{Z} _{\preccurlyeq G_0}^{01}$.
 \end{lemma}

%\begin{lemma}
%Assume that we have defined $MC$-cycles $Z^k \in \mathcal{Z}^k$ for each $k$ such that  there exists  $\tilde{Z}^k$  isotopy between $Z_{k+1}$ and $Z_k$  for each $k$. We can define inductively a $MC$-cycle $\dot{Z}_k$ and isotopy $\dot{\tilde{Z}}_k$ such that
%\begin{itemize}
%\item $\dot{Z}_k$ extends $\dot{Z}_{k-1}$;
%\item $\dot{\tilde{Z}}_k$ is an isotopy between $Z_k$ and $\dot{Z}_k$. 
%\end{itemize}
%\end{lemma}
%\begin{proof}
%Assume that $\dot{Z}_k$ is defined. Apply Lemma \ref{extension-isotopies-inductive} to the composition of $\tilde{Z}_{k}$ and $\dot{\tilde{Z}}_k$, to get $\dot{Z}_{k+1}$ and an isotopy between  $\dot{\tilde{Z}}_{k+1}$ between $Z_{k+1} $ and $\dot{Z}_{k+1}$.
%\end{proof}

\begin{lemma}  \label{inductive-extension}
Fix $J$. For each $k$, fix perturbation,etc and obtain  a $MC$-cycle $Z_k$ up to $k$. Assume that the perturbations are small enough so that there exists  $\tilde{Z}_k$  isotopy between $Z_{k+1}$ and $Z_k$  for each $k$. 
We can define inductively a $MC$-cycle $\dot{Z}_k$ and isotopy $\dot{\tilde{Z}}_k$ so that
\begin{itemize}
\item $\dot{Z}_k$ extends $\dot{Z}_{k-1}$;
\item $\dot{\tilde{Z}}_k$ is an isotopy between $Z_k$ and $\dot{Z}_k$. 
\end{itemize}
\end{lemma}
\begin{proof}
Assume that $\dot{Z}_k$ is defined. Apply Lemma \ref{extension-isotopies-inductive} to the composition of $\tilde{Z}_{k}$ and $\dot{\tilde{Z}}_k$, to obtain $\dot{Z}_{k+1}$ and an isotopy  $\dot{\tilde{Z}}_{k+1}$ between $Z_{k+1} $ and $\dot{Z}_{k+1}$.
\end{proof}

Let us take two choices $J_0, J_1$. 
As above, for each $k$, fix perturbation,etc and obtain   $Z_k^\mathfrak{0}$ ,  $\tilde{Z}_k^{\mathfrak{0}}$, $Z_k^{\mathfrak{1}}$  , $\tilde{Z}_k^{\mathfrak{1}}$,  for each $k$. Assume also that we have defined isotopies $\tilde{Z}_k$ between $Z_k^\mathfrak{0}$ and $Z_k^\mathfrak{1}$, and isotopy of isotopies $\tilde{\tilde{Z}}_k$ such that
%between  $\tilde{Z}_k$, $\tilde{Z}_k'$ , $\tilde{Z}_k$, $\tilde{Z}_{k+1}$.
\begin{itemize}
\item $\tilde{\tilde{Z}}_k$ agrees with   $\tilde{Z}_k^0$ for $s=0$;
\item $\tilde{\tilde{Z}}_k$ agrees with   $\tilde{Z}_k^1$ for $s=1$;;
\item $\tilde{\tilde{Z}}_k$ agrees with   $\tilde{Z}_k$ for $t=0$;
\item $\tilde{\tilde{Z}}_k$ agrees with   $\tilde{Z}_{k+1}$ for $t=1$.
\end{itemize}
Apply Lemma (\ref{inductive-extension}) to obtain $\dot{Z}_k^0 , \dot{\tilde{Z}}_k^0 , \dot{Z}_k^1,\dot{\tilde{Z}}_k^1$.
\begin{lemma}
We can define inductively a $MC$-cycle isotopy $\dot{\tilde{Z}}_k$ and isotopy of isotopy $\dot{\tilde{\tilde{Z}}}_k$ so that
$\dot{\tilde{Z}}_k$ extends $\dot{\tilde{Z}}_{k-1}$ and
%\item $\dot{\tilde{\tilde{Z}}}_k$ is an isotopy of isotopy between $\tilde{Z}_k^{\bullet}$ , $\dot{\tilde{Z}}_k$, $\dot{\tilde{Z}}_k'$ , $\dot{\tilde{Z}}_k^{\bullet}$. 
\begin{itemize}
\item $\tilde{\tilde{Z}}_k$ agrees with   $\dot{\tilde{Z}}_k^0$ for $s=0$;
\item $\tilde{\tilde{Z}}_k$ agrees with   $\dot{\tilde{Z}}_k^1$ for $s=1$;
\item $\tilde{\tilde{Z}}_k$ agrees with   $\tilde{Z}_k$ for $t=0$;
\item $\tilde{\tilde{Z}}_k$ agrees with   $\dot{\tilde{Z}}_k$ for $t=1$.
\end{itemize}
\end{lemma}
\begin{proof}
The argument is the same of Lemma \ref{inductive-extension}  using  Lemma \ref{extension-isotopies-isotopies}. 
\end{proof}

\section{Boundary Conditions}

We assume that $X$ contains a compact domain $X^{in}$ with boundary a smooth hypersurface $Y$ which is a contact manifold, such that $X \setminus X^{in}$ can be identified with $\R^+ \times Y$. We write
%such all positive area curves are in $X_{comp}$ and 
\begin{equation} \label{decomposition}
X = X^{in} \sqcup_{Y}   ( \R^+  \times Y),
\end{equation}
where $[1, + \infty) \times Y$ is the symplectization of $Y$.

---- such that the translation  in the $\R^+$ on $ \R^+ \times Y$  preserves  all (???) the geometric data (such as the complex structure.)

We assume that the lagrangians submanifolds $L$ which are  cylindrical Legendrian $\R^+  \times \partial X_{c}$. Compatible with (\ref{decomposition}), we assume that,    
$$ L = L^{in} \sqcup_{\Sigma}   (\R^+  \times \Sigma)$$ 
where $L^{in} = X^{in}  \cap L$ and  $\Sigma = \partial L^{in}= Y \cap L$ is a Legendrian submanifold of $\partial X_c$.
%Outside a compact set $X$ is equal to $X_{\partial \times \R^+}$.

We also assume that the four chain is written as
$$ K = K^{in} \sqcup_{T}   (\R^+  \times T)$$ 
where $K^{in} = X^{in}  \cap K$ and  $T = Y \cap K$.
% $\partial K_c = L_c + K_{\partial}$. 

We make the following assumption
\begin{itemize}
  \item For any real number $E >0$,   we can chose $X^{in}$ such that all the not-constant holomorphic curves with boundary on $(X,L)$ of area less than $E$ are mapped in $X^{in}$.  
 \end{itemize}
 We shall consider objects in the Novikov-Ring, and consider the objects up to $E$. 
 \begin{remark}
A standard way to achieve the assumption above is to impose convexity properties on the pair $(X,L)$.
\end{remark}

%This problem is present also in the case of closed Gromov-Witten invariants, that is $(L = \emptyset)$, where $K$ is just a closed four chain. In compact geometry the moduli of area zero spheres with three marked points mapped on $K$  leads to the intersection pairing $K \cap K \cap K$. In the not compact case the intersection depends on the choice of a boundary condition at infinity, that can be rephrased as the choice of a small perturbation of the diagonal of $X_{\infty}^3$ transversal to $K_{\infty}^3$.

As usual  in non-compact geometries, in order to obtain  well defined objects it is necessary to fix suitable boundary conditions at infinity . The assumptions above assure that the moduli space of pseudo-holomorphic curves of a fixed \emph{positive} area live in a compact region of $X$. However the moduli space of pseudo-holomorphic curves of area zero is always \emph{not} compact if $X$ is not compact.

The moduli space of pseudo-holomorphic curves  $\overline{\mathcal{M}}^{main}_{(g,h),(n, \overrightarrow{m})}(\beta=0)$ in the class $\beta=0$ set theoretically can be identified with $X \times  \overline{\mathcal{M}}^{main}_{(g,h),(n, \overrightarrow{m})}$ if $h=0$, or $L \times  \overline{\mathcal{M}}^{main}_{(g,h),(n, \overrightarrow{m})}$ if $h \neq 0$, where $\overline{\mathcal{M}}^{main}_{(g,h),(n, \overrightarrow{m})}$ is the Deligne-Munford compactification of the moduli of the Riemmanian surfaces. The   obstruction bundle (\ref{obstruction-bundle}) becomes a sub-vector space
\begin{equation} \label{obstruction-bundle0}
E_{\mathbf{p}} \subset   T_{p}X \otimes \Lambda^{0,1} (\Sigma),
\end{equation}
where $p \in X$ is the constant value of the pseudo-holomoprhic curve $\mathbf{p}$. 
%are defined locally on $X$, thus it makes sense to consider the moduli of area zero curves with target $Y \times \R^+$ . 

It makes sense to consider area zero curves whose target is an open subset of $X$. In particular we shall be interested to  area zero curves with target $ \R^+ \times Y$.

 In closed Gromov-Witten theory (which is the particular case of our problem where $L = \emptyset$ and $K$ is a closed four chain) a choice of a boundary condition consists on the choice of a perturbation of the moduli space of area zero curves on  $ \R^+ \times Y$ that is invariant under the translation on the $\R$-direction. Note that this makes sense since there is an obvious action of $ \R^+ $ on RHS of (\ref{obstruction-bundle0}), which is trivial on the second factor. 
Note that in the closed case the boundary condition affects only the degree zero invariant.  
%\begin{remark} \label{boundary-closed}
% For example ,  for area zero spheres with three marked points mapped on $K$, it defines the intersection pairing $K \cap K \cap K$, which depends on the choice on the choice of a small perturbation of the diagonal of $Y^3$ transversal to $K_{\infty}^3$. 
%\end{remark}

In the context of open-closed Gromov-Witten theory ( $L \neq \emptyset$) we need to work with the moduli of multi-curves $  \overline{\mathcal{M}}_{G,m}$. Since  there can be     area zero   multi-curve components in each degree,  the choice of a boundary condition affects the moduli of multi-curves  of degree not zero also, in contrast with what happen in the closed case .
% all the components of the moduli of multi-curves of positive area live in a compact region of $X$. However the moduli of the area zero components are not compact. 
%Therefore we need  to fix boundary conditions on the perturbation on the moduli of area zero pseudo-holomorphic curves. 

The generalization of the notion of  boundary condition used in the context of closed Gromov-Witten theory  is not straightforward since the  set of area zero components of the multi-curve  mapped on $ \R^+ \times Y$ can be different in different regions of the moduli space $\overline{\mathcal{M}}_{G}$. Hence it is necessary to set the invariance condition in the $\R$ direction consistently in a suitable sense.
%in the different regions and compatible with the boundary of the moduli. 
%We first explain how to solve the problem in the case $|m|=0$. 
 
%I we would like to do something similar, however the situation is more complicate since we need to deal with multi-curves. In this case there are different regions of the moduli space according to the set of components of the multi-curve mapped on $Y \times \R^+$.  It is necessary to interpolate among the different regions in a consistent way. 

\subsection{Boundary Conditions}

We first explain the definition of the boundary conditions for decorated graphs with  $|m|=0$.

%Before to consider a general manifold with boundary we 
% consider  the special case
% $$ Y \times \R .$$ 
   Let $G \in \mathfrak{G}$ with $\beta(G)=0$.
   Let  $\mathcal{P}= \{ P_i \}$ be an ordered partition of $Comp(G)$. 
   
  For each $i$, let $H_i = \sqcup_{c \in P_i}   H_c $.  
  % Let $E(\mathcal{P})$ be the set of pairs $(h,h') \in H \times H$ with $c,c' \subset P_i$ for some $i$, where $c,c'$ are the components of $G$ where $h,h'$ are attached respectively. 
 Set
  $$E(\mathcal{P}) = \sqcup_i \{ (h,h') \in H_i \times H_i | h \neq h' \}/ \Z_2 . $$

%  Denote by $E(\mathcal{P})$ the set of pairs $(h,h') \in H \times H$ with $c,c' \subset P_i$ for some $i$, where $c,c'$ are the components of $G$ where $h,h'$ are attached respectively. 
  
 Using the partition $\mathcal{P}$  define the decorated graph 
 $$cut_{\mathcal{P}}(G)= (Comp(G), E(  cut_{\mathcal{P}}(G) )  , H(G) ,V(G) , D(G)), $$
  obtained cutting the edges of $G$ which do not belongs to $E(\mathcal{P})$:
 $$ E^{in}(  cut_{\mathcal{P}}(G) ) = E^{in}(G) \cap E(\mathcal{P}),$$
 % defines  graphs   $G_i$  with $Comp(G_i) = P_i$, $H_v$ the same as $G$ for each $v$, $E^{in} (G_i) = E(P_i) \cap E^{in}(G)$.  
%To each $P_i$ it is associated a sub-graph $G^i$ of $cut_{\mathcal{P}}(G)$ and
We have a unique decomposition  of $G$ is sub-graphs
 \begin{equation}  \label{cut-graph}
cut_{\mathcal{P}}(G) = \sqcup_i G^i .
\end{equation} 
such that $Comp(G_i)= P_i$.

 Let 
 $$ \mathfrak{t} :  \R  \times Y \rightarrow \R$$
 be the projection on the $\R$ factor, which we refer as the time function.

For $\rho,\eta>0$, let  $ Y[\mathcal{P}, \rho, \eta]$ be the set of $  ( z_h )_h \in  ( \R  \times Y)^H$ such that 
\begin{itemize}
    \item  $\max_{h \in P_i} \mathfrak{t}(z_h)+ \eta < \min_{h \in P_{i+1}} \mathfrak{t}(z_h)$ for  $1 \leq i \leq f-1$ ;
    \item $  \max_{h \in P_i} \mathfrak{t}(z_h)  < \min_{h \in P_{i}} \mathfrak{t}(z_h) + \rho  $
for  $1 \leq i \leq f$ .
\end{itemize}

%For $G$, $\mathcal{P}= \{ P_i \}_i$ and $\rho, \zeta$ positive real numbers, set
 % $$  Y[\mathcal{P}, \rho, \zeta] = \{ \{ z_h \}_h \in (\R \times Y )^H   | \forall i ,   \max_{h \in P_i} \mathfrak{t}(z_h)+ \eta < \min_{h \in P_{i+1}} \mathfrak{t}(z_h),    \max_{h \in P_i} \mathfrak{t}(z_h)  < \min_{h \in P_{i}} \mathfrak{t}(z_h) + \rho  \}   $$

A $ \mathbf{boundary \thickspace data}$  $\mathcal{B }$ is specified by a sequence of  collections      of Kuranishi structures  and perturbations $\{  \overline{\mathcal{M}}_{G}^{\spadesuit,n}, \mathfrak{s}_{G}^{\spadesuit,n}, \}_{G}$
 with the following properties:
%For the  Kuranishi structures $\{  \overline{\mathcal{M}}_{G,m} \}_{G,m}(Y)$  we require that
\begin{itemize}
\item  conditions  of Subsection \ref{review-section} hold ;
\item  $\{  \overline{\mathcal{M}}_{G}^{\spadesuit,n}, \mathfrak{s}_{G}^{\spadesuit,n} \}_{G}$ are invariant under translation in the $\R$-direction;
\item For each $G,\mathcal{P},\rho >0 $ there exists $\eta >0$ such that on the region $ Y[\mathcal{P}, \rho/n, \eta/n]$ we have an identification of Kuranishi spaces
\begin{equation} \label{decomposition0-kuranishi}
\overline{\mathcal{M}}_{G}^{\spadesuit,n} = \prod_i \overline{\mathcal{M}}_{G_i}^{\spadesuit,n}
\end{equation}
and accordingly
\begin{equation} \label{decomposition0-perturbation}
\mathfrak{s}^{ \spadesuit,n}_{G} = \prod_i \mathfrak{s}^{\spadesuit , n}_{G_i},
\end{equation}
% we have an homomorphism XXX and an identification of Kuranishi spaces   $\overline{\mathcal{M}}_{G,m}  \cong \overline{\mathcal{M}}_{G,m, \mathcal{P}} $, when  $\overline{\mathcal{M}}_{G,m, \mathcal{P}} $ is endowed with the Kuranishi structure induced by (\ref{factorization-moduli}) . This identification involves a choice of TRIANGLE.
where the graphs $(G_i)_i$ are defined  in (\ref{cut-graph}).
\end{itemize}

%For the perturbations  $\mathfrak{s}_{G,m}$ we require that
%\begin{itemize}
%\item conditions ... hold;
%\item they are invariant under the translation in the $\R$-direction;
%\item For each $\rho $ there exists $\eta$ such that
%$\mathfrak{s}^n_{(G,m)}$ is compatible with  the identification of Kuranishi spaces (\ref{factorization-moduli}) when restricted to $ M[\mathcal{P}, \epsilon/n, \delta/n]$.
%\end{itemize}

\begin{lemma} \label{existence-boundary-data}
There exist boundary data. Two boundary data are isotopic.  
\end{lemma}
\begin{proof}
We proceed  by induction.

%Assume that as in the inductive argument of the proof of Lemma (\ref{basic-cycle}) we have defined $Z_{\prec G_0}$ and $Z_{G_0,m}$ for $|m| <l$ such that properties ... hold and moreover (\ref{coherent-limit}) holds.

Assume that we have defined $\{  \overline{\mathcal{M}}_{G}^{\spadesuit,n}, \mathfrak{s}_{G}^{\spadesuit,n} \}$ for $G \prec G_0$.

For each $\mathcal{P}$ we define $\rho_{\mathcal{P}}$,   $\zeta_{\mathcal{P}}$ ,  Kuranishi structure $\overline{\mathcal{M}}_{G, \mathcal{P}}^{\spadesuit,n}$, perturbations  $\mathfrak{s}^{ \spadesuit,n}_{G,\mathcal{P}} $  on the region $ Y(\mathcal{P}, \rho_{\mathcal{P}},   \zeta_{\mathcal{P}} )$ such that
\begin{enumerate}

\item  $\overline{\mathcal{M}}_{G, \mathcal{P}}^{\spadesuit,n}, \mathfrak{s}^{ \spadesuit,n}_{G,\mathcal{P}} $  satisfy the conditions of Subsection \ref{review-section}. 
% (forgetful compatible and  weakly-submersive ,  compatible with the corner faces);

%\item $$ Z_{G,m, \mathcal{P}_1 }^{\epsilon} =  Z_{G,m, \mathcal{P}_2 }^{\epsilon}$$ on $ T_{ R_{\mathcal{P}_1} \epsilon}(Diag_{\mathcal{P}_1}^{C_{\mathcal{P}_1} \epsilon}) \cap T_{ R_{\mathcal{P}_2} \epsilon}(Diag_{\mathcal{P}_2}^{C_{\mathcal{P}_2} \epsilon})   $;

\item $\overline{\mathcal{M}}_{G, \mathcal{P}}^{\spadesuit,n}, \mathfrak{s}^{ \spadesuit,n}_{G,\mathcal{P}} $ 
are translation invariant in the $\R$ direction;

\item
\begin{equation} \label{covering-boundary-condition}
\partial  Y(\mathcal{P}, \rho_{\mathcal{P}},   \zeta_{\mathcal{P}} ) \subset \bigsqcup_{\mathcal{P}' \prec \mathcal{P}}  Y(\mathcal{P}', \rho_{\mathcal{P}'},   \zeta_{\mathcal{P}'} ) ,
\end{equation}
\begin{equation} \label{not-overlapping-condition}
 Y(\mathcal{P}_1, \rho_{\mathcal{P}_1},   \zeta_{\mathcal{P}_1} ) \cap Y(\mathcal{P}_2, \rho_{\mathcal{P}_2},   \zeta_{\mathcal{P}_2} ) \neq \emptyset \Longrightarrow
  \mathcal{P}_1 \preccurlyeq \mathcal{P}_2 \text{   or   }  \mathcal{P}_2 \preccurlyeq \mathcal{P}_1  .
  \end{equation}

\item $\overline{\mathcal{M}}_{G, \mathcal{P}}^{\spadesuit,n}, \mathfrak{s}^{ \spadesuit,n}_{G,\mathcal{P}} $  extends the Kuranishi structures and perturbations defined on the LHS of (\ref{covering-boundary-condition}) by $\overline{\mathcal{M}}_{G, \mathcal{P}'}^{\spadesuit,n}, \mathfrak{s}^{ \spadesuit,n}_{G,\mathcal{P}'} $ for  $\mathcal{P}' \prec \mathcal{P}$.
% $ \partial M(\mathcal{P}, \rho_{\mathcal{P}},   \zeta_{\mathcal{P}} )$  by $\mathcal{P}' \prec \mathcal{P}$.

\end{enumerate}

Now we proceed by induction on $\mathcal{P}$.
For each $\mathcal{P}' \prec \mathcal{P}$, assume that we have defined  $\rho_{\mathcal{P}'},   \zeta_{\mathcal{P}'} $ and $\overline{\mathcal{M}}_{G, \mathcal{P}'}^{\spadesuit}, \mathfrak{s}^{ \spadesuit,n}_{G,\mathcal{P}'} $ .
% the Kuranishi structure on $M(\mathcal{P}', \rho_{\mathcal{P}'},   \zeta_{\mathcal{P}'} )$.
%$Z_{G,m, \mathcal{P}'}^{\epsilon} $, $C_{\mathcal{P}'}$,   $R_{\mathcal{P}'}$ with the property above.

Take $\rho_{\mathcal{P}}$ big enough such that (\ref{covering-boundary-condition}) holds.

Take $\eta_{\mathcal{P}}$ such that  (\ref{decomposition0-kuranishi})abd  (\ref{decomposition0-perturbation}) hold for $G \prec G_0$, 
and moreover property  (\ref{not-overlapping-condition}) holds.

%If $|\mathcal{P}| >1$, on $\partial M[\mathcal{P}, \rho/n, \eta/n] $ there exists a unique homeomorphism ... compatible with  $ \overline{\mathcal{M}}_{G,m, \mathcal{P}'} \cong  \overline{\mathcal{M}}_{G,m} $ $\mathcal{P}'$ for each $\mathcal{P}' \prec \mathcal{P}$. Pick an extension of it on $ M[\mathcal{P}, \rho/n, \eta/n] $ which is compatible with the one defined in the preview steps on the boundary strata $\overline{\mathcal{M}}_{G,m} (G',m',E')$ for each $(G',E',m') \in \mathfrak{G}(G,m)$ . Use it to  define the Kuranishi structure on $M[\mathcal{P}, \rho/n, \eta/n]$ compatible with (\ref{factorization-moduli}).

If $|\mathcal{P}| >1$, condition  (\ref{decomposition0-kuranishi}) defines the Kuranishi structure on $ Y[\mathcal{P}, \rho/n, \eta/n] $.

If $\mathcal{P}$ is the trivial partition, the Kuranishi structures defined in the preview steps and condition (\ref{not-overlapping-condition}) fixes the Kuranishi structure on the RHS of   (\ref{covering-boundary-condition}).
%$\partial Y[\mathcal{P}, \rho/n, \eta/n] $. 
Extend this Kuranishi structure to a Kuranishi structure on  $M[\mathcal{P}, \rho/n, \eta/n]$ which is compatible with one defined in the preview steps on the corner faces $\overline{\mathcal{M}}_{G} (G',E')_{(G',E') \in \mathfrak{G}(G)}$.
In the same way we can construct the perturbation $\mathfrak{s}^{ \spadesuit , n}_{G}$.

\end{proof}

\subsubsection{Boundary Conditions} \label{boundary-conditions-subsection}

An open partition $\mathcal{P}^{open}$  consists of an array $(P_{in}, P_1,...,P_n )$  of disjoint sub-sets of $Comp(G)$ with
\begin{itemize}
    \item $P_i \neq \emptyset $ for each $i$,
    \item  $ P^{in} \sqcup \sqcup_i P_i = Comp(G)  .$
\end{itemize}
There is an obvious partial ordering on the set of open partitions.

%Extend the function $\mathfrak{t}$ used above to $X$ setting $\mathfrak{t}(z)=0$ for each $z \in M^{in}$. Use this to define the functions $max_S$ and $min_S$ on $M^H$, and thus the regions  $X[\mathcal{P}, \rho, \zeta] \subset X^H$ using the same formula (\ref{interpolating-regions}).

Extend  to $M$ the function $\mathfrak{t}$ used in the preview subsection   setting $\mathfrak{t}(z)=0$ for each $z \in M^{in}$. 
%With this change the definition (\ref{interpolating-regions}) applies to (\ref{end-manifold}) to obtain the regions  $M[\mathcal{P}, \rho, \eta] \subset M^H$.
For $\rho,\eta>0$, let  $ M[\mathcal{P}, \rho, \eta]$ be the set of $  \{ z_h \} \in  M^H$ such that 
\begin{itemize}
    \item   $\max_{h \in H^{in}} \mathfrak{t}(z_h)+ \eta < \min_{h \in P_1} \mathfrak{t}(z_h)$, 
$\max_{h \in H^{in}} \mathfrak{t}(z_h) < \rho$;
    \item  $\max_{h \in P_i} \mathfrak{t}(z_h)+ \eta < \min_{h \in P_{i+1}} \mathfrak{t}(z_h)$ for  $1 \leq i \leq f-1$ ;
    \item $  \max_{h \in P_i} \mathfrak{t}(z_h)  < \min_{h \in P_{i}} \mathfrak{t}(z_h) + \rho  $
for  $1 \leq i \leq f$ .
\end{itemize}

We say that a collection of Kuranishi structures and perturbations $\{  \overline{\mathcal{M}}_{G}, \mathfrak{s}_{G} \}_{G}$ are compatible with the $\mathbf{boundary \thickspace condition}$ specified by the boundary data $\mathcal{B}$ if
\begin{itemize}
\item  the conditions of Subsection \ref{review-section} hold  ;
\item for each $G,\mathcal{P},\rho $ there exists $\eta$ such that on the region $ M[\mathcal{P}, \rho/n, \eta/n]$ we have an identification of Kuranishi spaces
$$\overline{\mathcal{M}}_{G} = \overline{\mathcal{M}}_{G^{in}}\times  \prod_i \overline{\mathcal{M}}_{G_i}^{\spadesuit}$$
and accordingly
$$\mathfrak{s}^n_{G} = \mathfrak{s}^n_{G^{in}} \times  \prod_i \mathfrak{s}^{\spadesuit , n}_{G^i}.$$
% we have an homomorphism XXX and an identification of Kuranishi spaces   $\overline{\mathcal{M}}_{G,m}  \cong \overline{\mathcal{M}}_{G,m, \mathcal{P}} $, when  $\overline{\mathcal{M}}_{G,m, \mathcal{P}} $ is endowed with the Kuranishi structure induced by (\ref{factorization-moduli}) . This identification involves a choice of TRIANGLE.
\end{itemize}

The following is proved adapting the inductive argument used in the proof of Lemma \ref{existence-boundary-data}
\begin{lemma}
For each boundary data $\mathcal{B}$, there exist Kuranishi structures and perturbations compatible with $\mathcal{B}$.
\end{lemma}

\subsubsection{Extension to $|m|>0$}

We now adapt the definition of boundary condition provided above  to graphs $(G,m) \in \mathfrak{G}_*$ with $|m|>0$. We need to extend (\ref{decomposition0-kuranishi}), (\ref{decomposition0-perturbation}) in a way compatible with the boundary of the moduli space. 
For this we need to introduce some notation. 

Consider an order set $\mathbf{k}= \{ k_0, k_1, ..., k_f \}$ of integer numbers with $k_i \leq k_{i+1}$, $0 \leq k_i \leq l$. We assume 
\begin{equation} \label{assumption-k-0f}
k_0=0,k_f=l .
\end{equation}

Denote by $\Delta^{l}$ the standard simplex of dimension $l$. 
Set
$$\Delta^{l,\mathbf{k}} = \prod_{i=1}^{f}   \Delta^{k_i - k_{i-1}}.$$

We label the boundary faces of  $ \Delta_{l,\mathbf{k}} $ in the following way:
\begin{itemize}
    \item for each $i \notin \mathbf{k}$, $\partial_i \Delta_{l,\mathbf{k}} $ ;
    \item if   $k_j > k_{j-1}$, $\partial_{k_j}^- \Delta^{l,\mathbf{k}} $;
    \item  if $k_j < k_{j+1}$, $\partial_{k_j}^+\Delta^{l,\mathbf{k}} $.
\end{itemize}

%it is associated a boundary face  $\partial_i \Delta_{l,\mathbf{k}} $ of $\Delta^{l,\mathbf{k}} $ for each $0 \leq i \leq l$ such that $i \neq k_j$ for each $j$ . There is a boundary face    $\partial_{k_j}^- \Delta^{l,\mathbf{k}} $  if   $k_j > k_{j-1}$ and a boundary face $\partial_{k_j}^+\Delta^{l,\mathbf{k}} $ if $k_j < k_{j+1}$. 
%%of $\Delta_{k_0,k_1,..., k_f} $.

We have canonical identifications
\begin{equation} \label{identification-faces-simplex}
 \partial_{k_j}^-  \Delta^{l,\mathbf{k}} = \partial_{k_j}^+\Delta^{l,\mathbf{k}}   
 \end{equation}
whenever $k_j > k_{j-1}$, and
% and
%$$d_{i,l} (\partial_i \Delta_{k_0,k_1,..., k_f} ^l) \cong  \Delta_{k_0,k_1,..., k_f}^{l-1}  $$
%where $d_{i,l} : XXX$.
\begin{equation} \label{identification-faces-simplex1}
    \partial_i \Delta^{l,\mathbf{k}}  = \Delta^{l-1, \partial_i \mathbf{k}} 
\end{equation}
where we have denoted $\partial_i \mathbf{k} = \{ k_0, ...,k_{j-1}, k_j -1, ...,k_f-1  \}$ if $ k_{j-1} < i < k_j$.
%$k_j' = k_j$ if $k_j <i$,and $k_j' = k_j-1$ if $k_j >i$.

It is elementary to check that 
\begin{equation} \label{decomposition-simplex}
 ( \bigsqcup_{\mathbf{k}}   \Delta^{\mathbf{k}}  )/\sim \cong \Delta^l .
\end{equation}
where the LHS is defined making the identifications (\ref{identification-faces-simplex}). The identification  (\ref{decomposition-simplex}) is far to be unique.

%Pick an identification  (\ref{decomposition-simplex}) for each $l$. We say that they are compatible if the follow

Note that a choice of  (\ref{identification-faces-simplex}) for $l$ and $l-1$ together with  (\ref{identification-faces-simplex1}) yield to an identification
$ \partial_i \Delta^l = \Delta^{l-1} .$
We shall consider identification  (\ref{decomposition-simplex}) for different $l$ such that the last equation agrees with the standard embedding of $ \Delta^{l-1}$ with the boundary face $\partial_i \Delta^l$.

%There exists a decomposition of the symplex  $\Delta^l$
%\begin{equation} \label{decomposition-simplex}
%\Delta^l = ( \bigsqcup_{\mathbf{k}}   \Delta^{\mathbf{k}}  )/\sim
%\end{equation}
%where the RHS is defined using  the identifications (\ref{identification-faces-simplex}). The decomposition \ref{decomposition-simplex} is far to be unique.

Consider a decorated graph  $(G,m)$ and a partition $\mathcal{P}$ as above. Write $m = \{E_0, E_i, ...,E_l \}.$ 
%For $m = \{E_0, E_1,...,E_l  \}$,    assume that
% \begin{equation} \label{condition-edges}
% E_l  \cap E(\mathcal{P}) = \emptyset
%  \end{equation}
% If 
%$E_l \subset E(\mathcal{P})$,  set
%$$ m_{[a,b]}^i= \{ E_a \cap E(G_i),  E_{a+1} \cap E(G_i) , ...,  E_b \cap E(G_i) \} $$

Set 
$$  cut_{\mathcal{P}}(G,m) = (cut_{\mathcal{P}}(G), cut_{\mathcal{P}}(m) )  $$
 where 
 $$cut_{\mathcal{P}}(m) = \{ E^{ex}(cut_{\mathcal{P}}(G)),  E_{1} \cap E(\mathcal{P}) , ...,  E_l \cap E(\mathcal{P}) \}   $$
 if $m = \{E_0, E_1,...,E_l  \}$ with $E^0= E^{ex}(G)$.
Define $m^i$ from  the relation 
 \begin{equation}  \label{cut-graph2}
cut_{\mathcal{P}}(G,m) = \sqcup_i (G^i , m^i ).
\end{equation} 

Set
\begin{equation} \label{decomposition-moduli}
\overline{\mathcal{M}}_{G,m, \mathcal{P}, \mathbf{k}}  =       \prod_{0 \leq i < f} \overline{\mathcal{M}}_{(G^i, m_{[k_i,k_{i+1}]}^i )}   ,
\end{equation}
where we have used the notation 
 %$E_l \subset E(\mathcal{P})$,  set
$$ m_{[a,b]}= \{ E_a ,  E_{a+1}  , ...,  E_b  \} ,$$
if $m = \{E_0, E_1,...,E_l  \}$.

If the spaces associated to $(G^i,m^i)$ are endowed with a Kuranishi structure, relation (\ref{decomposition-moduli}) defines a Kuranishi structure on $\overline{\mathcal{M}}_{G,m, \mathcal{P}, \mathbf{k}}$.

A choice of the identification  (\ref{decomposition-simplex})  yields to an identification 
\begin{equation} \label{identification-factorization-moduli}
 \overline{\mathcal{M}}_{G,m} = ( \sqcup_{ \mathbf{k}} \overline{\mathcal{M}}_{G,m, \mathcal{P}, \mathbf{k}} )/\sim  .
 \end{equation}

We want to refine condition (\ref{decomposition0-kuranishi})  asking that  (\ref{identification-factorization-moduli}) becomes an identification of Kuranishi spaces compatible with the corner faces. In particular this involves a choice of an identification ( \ref{decomposition-simplex}) depending on the point on moduli space $\overline{\mathcal{M}}_{G}$.  This identification has to be compatible at the corner faces with the identification made in the the previews inductive steps ,
%Note that it is not possible to choice a decomposition ( \ref{decomposition-simplex}) globally on the moduli space respecting the compatibility with the corner faces. 

%an identification of Kuranishi spaces   $\overline{\mathcal{M}}_{G,m}  \cong \overline{\mathcal{M}}_{G,m, \mathcal{P}} $, when  $\overline{\mathcal{M}}_{G,m, \mathcal{P}} $ is endowed with the Kuranishi structure induced by (\ref{factorization-moduli}) . 
%This identification involves a choice of TRIANGLE.

Now we can adapt the inductive argument  of Lemma \ref{existence-boundary-data}  in this setting. The $|\mathcal{P}|=1$ is the same of Lemma \ref{existence-boundary-data}.  
Assume $|\mathcal{P}| >1$.  The choice of an identification ( \ref{decomposition-simplex}) for $\mathcal{P}' \prec \mathcal{P}$ induces an identification ( \ref{decomposition-simplex}) on the region defined by the RHS of (\ref{covering-boundary-condition}).  Pick an extension of it on 
$\overline{\mathcal{M}}_{G, \mathcal{P}}$
%$ M[\mathcal{P}, \rho/n, \eta/n] $ 
which is compatible with the identification inducted on the boundary strata $\overline{\mathcal{M}}_{G,m} (G',m', E')$ for each $(G',m', E') \in \mathfrak{G}(G,m)$  by the one defined in the preview steps . Use it to  define the Kuranishi structure of $\overline{\mathcal{M}}_{G,m}$ on the region $ Y[\mathcal{P}, \rho/n, \eta/n] $  compatible with  (\ref{identification-factorization-moduli}).

\subsubsection{ Change of Boundary Data }

%We shall consider perturbations which are invariant under the action of $\R^+$.  A choice of a boundary condition $\mathcal{B}$ includes the choice of $K_{\infty}$ and a choice of the perturbation of the associated moduli space of area zero curves. 

After we fix a boundary data $\mathcal{B}$, the construction of section (\ref{review-section}) applied to a perturbation of the moduli space considered on subsection \ref{boundary-conditions-subsection}  yield to a $MC$-cycle 
$$ Z^{\mathcal{B} } \in \mathcal{Z}  $$
which is well defined up to isotopy.

Now we want to consider how $Z^{\mathcal{B} }$ changes when we change $\mathcal{B}$.
Fix two boundary data  $\mathcal{B}_0$  and  $\mathcal{B}_1$. 
In a similar way to what we have done in subsection \ref{boundary-conditions-subsection}, we consider the problem on $\R \times Y$ and define  Kuranishi structures and perturbations compatible with the boundary condition $\mathcal{B}_0$ for $t \ll 0$ and $\mathcal{B}_1$ for $t \gg 0$. 
This yield to a $MC$-cycle
$$ Z^{\mathcal{B}_0 , \mathcal{B}_1 } \in \mathcal{Z}   .$$

A Kuranishi structure and pertubation on the moduli space on $X$ compatible with the boundary data $\mathcal{B}_0$ and a Kuranishi structure on the moduli $\R \times Y$ compatible with with the boundary data $\mathcal{B}_0$  and  $\mathcal{B}_1$ as above can be glued to obtain a Kuranishi structure and perturbation on the moduli space on $X$ compatible with the boundary data $\mathcal{B}_1$. The corresponding $MC$-cycles satisfies  a relation which we write schematically as
$$ Z^{\mathcal{B}_1} = Z^{\mathcal{B}_0 } \boxtimes  Z^{\mathcal{B}_0 , \mathcal{B}_1 } . $$

The same considerations yield to the relation
$$ Z^{\mathcal{B}_0 , \mathcal{B}_2 } = Z^{\mathcal{B}_0 , \mathcal{B}_1 } \boxtimes Z^{\mathcal{B}_1 , \mathcal{B}_2 } $$
when we consider three boundary datas $\mathcal{B}_0$  ,  $\mathcal{B}_1$ and  $\mathcal{B}_2$.

%----

%For boundary conditions $\mathcal{B}_1$ and $\mathcal{B}_2$  we can consider on $X_{\partial} \times [a_1,a_2]$  the pertubation defined by  $\mathcal{B}_1$ on $X_{\partial} \times [a_1,b_1]$ and $\mathcal{B}_2$ on $X_{\partial} \times [b_2,a_2]$. Using the general machine this define an element
%$$ Z_{\mathcal{B}_1 , \mathcal{B}_2 } \in MCH(X_{\partial} \times [a_1,a_2])  .$$
%This satisfies the relation
%$$ Z_{\mathcal{B}_1 , \mathcal{B}_3 } = Z_{\mathcal{B}_1 , \mathcal{B}_2 } + Z_{\mathcal{B}_2 , \mathcal{B}_3 } $$
%for each choice of boundary conditions $\mathcal{B}_1$ , $\mathcal{B}_2$ and $\mathcal{B}_3$.

%From the definition it follows directly the fundamental relation
%$$ Z_{\mathcal{B}_2} = Z_{\mathcal{B}_1 } +  Z_{\mathcal{B}_1 , \mathcal{B}_2 } . $$

\subsection{Frame refined Four Chain}

%In this subsection we refine the usual relative homology groups in order to fix partially the boundary conditions at infinity. 

%We have said that in order to define Open Gromov-Witten we need to fix boundary condition at infinity. 

%The problem of fixing the boundary conditions at infinity is present also in the case of closed Gromov-Witten invariants, where the invariant depends on the choice of a boundary condition of the moduli of area zero spheres with three marked points mapped on $K$, that defines the intersection pairing $K \cap K \cap K$. In the closed case $(L = \emptyset)$ and $K$ is a closed four chain. The boundary condition can be rephrased as the choice of a small perturbation of the diagonal of $X_{\infty}^3$ transversal to $K_{\infty}^3$.

In order to fix partially the boundary conditions at infinity 
In this section   we define a refined version of the four chain.
%extension of the relative homology of $K_{\infty} \in H_3(X_{\infty}',\partial X_{\infty}')$, 
The residue ambiguity will be analogous to the one arising in Closed Gromov-Witten invariants.  
%( i.e., $L = \emptyset$)  of remark \ref{boundary-closed}. 

%We shall fix topological data at infinity such that the residue ambiguity is the same of the closed case.
%This will involve the relative homology of $K_{\infty} \in H_3(X_{\infty}',\partial X_{\infty}')$ improved by the extra data provided by a frame of $T^*L$. We now define the topological object we need. Let us first recall some basic fact about frames on $3$-manifolds.

%We want to fix partially the boundary condictions in terms of topological data. In order to do this we need to introduce frames at infinity in the following sense. 

\subsubsection{Frames and Spin Structures}

We start recalling some standard and elementary fact about frames on a compact oriented $3$-manifold $M$. Fix a riemmanian metric on $M$ (up to isotopy , all the objects we will consider do not dependent on the metric).

Denote by $Spin(M)$ the set of spin structures on $M$. Recall that $Spin(M)$ is a torsor on the group $H^1(M, \Z_2)$.
Since $Spin(3) \approx S^3$ as differential manifold, from the degree of the map we obtain an isomorphism of groups
\begin{equation} \label{degree-spin}
Maps(M, Spin(3))/isotopies \cong \Z .
\end{equation}
Moreover, since $Spin(3)$ is a $\Z_2$-covering of $SO(3)$, there is  an exact sequence  of groups
\begin{equation} \label{exact-sequence-spin}
 0 \rightarrow Maps(M, Spin(3))/isot \rightarrow Maps(M,SO(3))/isot \rightarrow  H^1(M, \Z_2) \rightarrow 0  .
 \end{equation}

Denote by $\mathfrak{Fr}(M)$ the set of orthogonal frames of $M$, i.e., the set of orthogonal trivializations of the tangent bundle of $M$. $\mathfrak{Fr}(M)$ is a torsor on $Maps(M, SO(3))$.
Denote by $\mathfrak{Fr}(M)/isot$ the set of orthogonal frames modulo isotopies, which is
  a torsor on $Maps(M, SO(3))/isot$. 

Since a frame in particular determines a spin structure, we have a natural map of sets
$$\mathfrak{Fr}(M)/isot \rightarrow Spin(M) .$$
%$\mathfrak{Fr}(M)$ is a torsor on $Maps(M, SO(3))$, and $\mathfrak{Fr}(M)/isot$ is a torsor on $Maps(M, SO(3))/isot$. 
From (\ref{exact-sequence-spin})  we obtain the identification of sets
\begin{equation} \label{frames-spin}
(\mathfrak{Fr}(M)/isot) / \Z \cong Spin(M).
\end{equation}
where $\Z$ acts on $\mathfrak{Fr}(M)/isot$ throughout the isomorphism (\ref{degree-spin})

%Since $\mathfrak{Fr}(M)$ is a torsor on $Maps(M, SO(3))$, and $\mathfrak{Fr}(M)/isot$ is a torsor on $Maps(M, SO(3))/isot$. from (...)  it follows the sequence
%$$0 \rightarrow \Z \rightarrow \mathfrak{Fr}(M)/isot \rightarrow Spin(M) \rightarrow 0.$$

For our purpose it will  be useful to consider the restriction to $Maps(M, SO(2))$ of  the action of $Maps(M, SO(3))$ on $\mathfrak{Fr}(M)$ , where 
$SO(2)$ is embedded on $SO(3)$ from the inclusion $\R^2 \cong \{ 0 \} \times \R^2 \subset \R^3$.  Note that
$$ Maps(M, SO(2))/isot \cong H^0(M, U(1))/ H^0(M, \R )  \cong H^1(M, \Z )  $$
From this, we obtain an action of the group $H^1(M, \Z ) $ on  $\mathfrak{Fr}(M)/isot$:
\begin{equation} \label{action-H^1(M, Z )}
H^1(M, \Z ) \times \mathfrak{Fr}(M)/isot \rightarrow \mathfrak{Fr}(M)/isot .
\end{equation}

%---

 %Since $Spin(3)$ is the universal covering of $SO(3)$, we have the exact sequence
%$ 0 \rightarrow Maps(M, Spin(3))/isot \rightarrow Maps(M,SO(3))/isot \rightarrow  Hom(\pi_1(M), \pi_1(SO(3))) \rightarrow 0  .$
%Since $Maps(M, Spin(3))/isot \cong \Z$ with isomorphism provided by the degree of the map, and $Hom(\pi_1(M), \pi_1(SO(3))) \cong H^1(M, \Z_2)$, it follows the exact sequence
%$$ 0 \rightarrow \Z \rightarrow Maps(M,SO(3))/isot \rightarrow  H^1(M, \Z_2) \rightarrow 0  .$$

\subsubsection{The case $\Sigma \times \R$}
Consider the particular case   $M=\Sigma \times \R_t$ (where $t$ denotes the coordinate on the factor $\R$, which we refer as time coordinate).
% In this case there is an obvious action of the group $\R $ on $\mathfrak{Fr}(\Sigma \times \R)$ by translation in the $t$-direction. 
Let $\mathfrak{Fr}(\Sigma \times \R)^{\R}$ be the set of elements of $\mathfrak{Fr}(\Sigma \times \R)$ invariant under translation in the $t$-direction. 
$\mathfrak{Fr}(\Sigma \times \R)^{\R}$ can be identified with $\mathfrak{Fr}(T\Sigma \oplus \R)$, where $\R$ denote the trivial real bundle of dimension one over $\Sigma$.

We say that two frames $\mathbf{fr}_1, \mathbf{fr}_2 \in \mathfrak{Fr}(\Sigma \times \R)^{\R} $ are isotopic if there exists $\tilde{\mathbf{fr}} \in \mathfrak{Fr}(\Sigma \times \R)$ such  that 
%agrees with $f_1$ for $t <<0$ and with  $f_2$ for $t>> 0$:
$$   \tilde{\mathbf{fr}} \equiv \mathbf{fr}_1  \text{    for      } t  \ll 0   ,  \quad  \tilde{\mathbf{fr}} \equiv \mathbf{fr}_2  \text{    for      } t  \gg 0 .$$
%For the same reason of  (\ref{frames-spin}),  
Two frames are isotopic if and only if they induce the same spin structure on $\Sigma$:
$$   (  \mathfrak{Fr}(\Sigma \times \R)^{\R}) /isotopies  \cong Spin(\Sigma) .$$

%\subsubsection{Action of $H_1(\Sigma, \Z)$ on $ \mathfrak{Fr}(\Sigma \times \R)^{\R}$} 
 
The group action (\ref{action-H^1(M, Z )})  induces on $\Sigma \times \R$ the action of $ H^1(\Sigma, \Z)$ on $\mathfrak{Fr}(\Sigma \times \R)^{\R}/iso$.
%In particular we have an action of $Maps(\Sigma, SO(2))$ on $\mathfrak{Fr}(\Sigma \times \R)^{\R}$.
%Since  $Maps(\Sigma, SO(2))/isotopies = H^1(\Sigma, \Z)$,  this defines an action of $ H^1(\Sigma, \Z)$ on $\mathfrak{Fr}(\Sigma \times \R)^{\R}/iso$.

%To $V_0, V_1 \in \Gamma(\Sigma, T\Sigma \oplus \R )$ compatible with $\mathbf{fr}_0,\mathbf{fr}_1 \in  \mathfrak{Fr}(T\Sigma \oplus \R)$ we can associate an element  $[graph(\tilde{V})] \in C_3(S(T\Sigma \oplus \R))/ \partial C_4(S(T\Sigma \oplus \R))$ with $\partial [graph(\tilde{V})] = graph(V_1) - graph(V_0)$.

Now consider pairs 
\begin{equation} \label{pair-frame-vector}
    (\mathbf{fr}, V)
\end{equation}
 where $ \mathbf{fr} = (fr_1,fr_2,fr_3) \in \mathfrak{Fr}(T\Sigma \oplus \R)$ and 
$V \in \Gamma(\Sigma, T\Sigma \oplus \R )$ is compatible with $ \mathbf{fr}$ in the following sense:
 there exit $a_1,a_2,a_3 \in \R$ such that
\begin{equation} \label{compatible-vector-frame}
     V|_{\Sigma} = a_1 fr_1 + a_2 fr_2 + a_3 fr_3 .
\end{equation}
 
%We say that 
%$V \in \Gamma(\Sigma, T\Sigma \oplus \R )$ is compatible with $\mathbf{fr} = (fr_1,fr_2,fr_3) \in \mathfrak{Fr}(T\Sigma \oplus \R) $
%if there exit $a_1,a_2,a_3 \in \R$ such that
%$$  V|_{\Sigma} = a_1 fr_1 + a_2 fr_2 + a_3 fr_3 .$$

%%Denote by $C_3(S(TM), \Z| \mathbf{fr})$ the vector space of chains which are compatible with $\mathbf{fr}$.

%Now consider pairs $(\mathbf{fr}, V)$ where $ \mathbf{fr} \in \mathfrak{Fr}(T\Sigma \oplus \R)$ and 
%$V \in \Gamma(\Sigma, T\Sigma \oplus \R )$ is compatible with $ \mathbf{fr}$. 
The difference of two pairs $(\mathbf{fr}_0, V_0), (\mathbf{fr}_1, V_1)$ is defined as
\begin{equation} \label{difference-frames}
    (\mathbf{fr}_0, V_0) - (\mathbf{fr}_1, V_1) = [graph(\tilde{V})] \in C_3(S(T\Sigma \oplus \R)) / \partial C_4(S(T\Sigma \oplus \R)) ,
\end{equation}
%where we quotient by $\partial C_4(S(T\Sigma \oplus \R))$.
where $\tilde{V}=(V_t)_{t \in [0,1]}$, $ V_t \in \Gamma(S(TM)|_{\Sigma})$ and
 there exists an isotopy of frames  $(\mathbf{fr}_t)_{t \in [0,1]} $ with   $V_t$  compatible with  $\mathbf{fr}_t$ for each $t$. Here $graph (\tilde{V})$ is  considered as a map $ \Sigma \times [0,1] \rightarrow S(TM)$.

There is an obvious action of the group $ Maps(\Sigma, SO(2))$ on the set of pairs (\ref{pair-frame-vector}) defined by
\begin{equation} \label{action-pairs}
    A  \cdot (\mathbf{fr}, V)  = (A_V \cdot \mathbf{fr}, V)
\end{equation}
where $A_V$ is the rotation in the direction $V$ with angle defined by $A \in Maps(\Sigma, SO(2))$ .
 
It is easy to check  the following formula
%If $A \in Maps(\Sigma, SO(2))/isot \cong H^1(M, \Z)$ fixes $V_0$, 
%$$ ( A \cdot fr , fr) \mapsto PD([A])  .$$
\begin{equation} \label{action-pairs-difference}
    A \cdot (\mathbf{fr}_0, V_0) -  (\mathbf{fr}_1, V_1) = ( (\mathbf{fr}_0, V_0) - (\mathbf{fr}_1, V_1) ) +  PD([A]) ,
\end{equation}
 where $[A] \in Maps(\Sigma, SO(2))/isot \cong H^1( \Sigma, \Z)$, and 
$PD([A])  \in H_1( \Sigma, \Z) = H_3^-(S(T\Sigma \oplus \R))$.

%---

%There is a map 
%$$ \mathfrak{Fr}(\Sigma \times \R)_0^{\R} \times  \mathfrak{Fr}_0(\Sigma \times \R)^{\R} \rightarrow H_3(S(TM)),$$
%$$(fr, fr') \mapsto (\tilde{fr} \cdot (1,0,0))_*([\Sigma]).$$

%If $A \in Maps(\Sigma, SO(2))/isot \cong H^1(M, \Z)$, 
%$$ ( A \cdot fr , fr) \mapsto PD([A])  .$$

%To a map $A: \Sigma \rightarrow SO(2)$ it is associated an action on $Frames(\Sigma \times \R)^{\R}$. Namely   the frame $A \cdot (e_1,e_2)$ is given by rotating by $A$ the vector $e_2$ in the plane orthogonal to $e_1$.
%Since $H_1(\Sigma, SO(2))= H_1(\Sigma , \Z)$, we get an action of $H_1(\Sigma, \Z )$ on $Frames(\Sigma \times \R)^{\R}$. 

%Identifying $SO(2)$ with the subgroup of $SO(3)$ that fixes the vector $(1,0,0)$, we have a map $H_0(\Sigma, SO(2)) \rightarrow H_0(\Sigma, SO(3))$.
%This defines an action of the group $H_1(\Sigma , \Z) = H_0(\Sigma, SO(2))$ on $Frames(\Sigma \times \R)^{\R}$.

%For $ (e_1,e_2) \in Frames(\Sigma \times \R)^{\R}$ and  $A: \Sigma \rightarrow SO(2)$, we see  we identify $SO(2)$ with the subgroup of $SO(3)$ that fixes $e_1$. 
% the frame $A \cdot (e_1,e_2)$ is given by rotating by $A$ the vector $e_2$ in the plane orthogonal to $e_1$.
 
\subsubsection{Euler Structures: the compact case}

An Euler Structures on a compact oriented three manifold $M$ is a nowhere vanishing vector field on $M$.  Two Euler Structures on $M$ are called homologous if they are isotopic outside a ball as nowhere vanishing vector fields. We denote by $\mathfrak{Eul}(M)$  the set of homology classes of Euler Structures.  $\mathfrak{Eul}(M)$ is a torsor on $H_1(M, \Z)$.
The Chern Class of an Euler class $[V]$ is defined as $\mathfrak{c}([V])= V - opp(V) \in H_1(M,\Z)$, where $opp(V)$ denotes the opposite of the vector field $V$.

 We now introduce  an  equivalent  definition of $\mathfrak{Eul}(M)$. 
%There is an obvious map $\mathfrak{Eul}(M) = \Gamma(S^2(TM) \rightarrow M)$.
Denote by $S(TM)$ the spherical tangent bundle of $M$ and let
$pr:S(TM) \rightarrow M $
be the projection. 
Set 
$$\mathfrak{Eul}(M)^{\blacklozenge} := (pr_*)^{-1}([M])  \subset H_3(S(TM), \Z).$$
It is immediate that $\mathfrak{Eul}(M)^{\blacklozenge}$ is a torsor on 
$$ H_3^-(S(TM), \Z)  := \text{Ker} \{ pr_*: H_3(S(TM), \Z) \rightarrow H_3(M , \Z) \}.$$
This group can be also identified with the set of  $r$-anti-invariant elements of $ H_3(S(TM), \Z)$, where  $r : TM \rightarrow TM$ is the reverse map .
% $ H_3^-(S(TM), \Z) $ can be identified with  the kernel of the map $\pi_*: H_3(S(TM), \Z) \rightarrow H_3(M , \Z)$ induced by the projection map $\pi:S(TM) \rightarrow M$:
 % It is trivial to check that  $H_3^-(S(TM), \Z) \cong H_1(M,\Z)$.

%We have the exact sequence:
%$$ 0 \rightarrow H_3^-(S(TM), \Z) \rightarrow  H_3(S(TM), \Z) \rightarrow H_3(M , \Z) \rightarrow 0.$$

%where $S^2(TM)$ is the spherical tangent bundle of $M$, and $ H_3^-(S^2(TM), \Z) $ are the of element anti-invariant  $ H_3(S^2(TM), \Z)$ by the reverse isomorphism of $TM$. It is trivial to prove that  $H_3^-(S^2(TM), \Z) \cong H_1(M,\Z)$.

% $ H_3^-(S(TM), \Z) $ is the group of elements $r$-anti-invariant of $ H_3(S(TM), \Z)$, where  $r : TM \rightarrow TM$ is the reverse homomorphism and $H_3^-(S(TM), \Z) \cong H_1(M,\Z)$. Hence  

Since $H_3^-(S(TM), \Z) \cong H_1(M,\Z)$,
$\mathfrak{Eul}(M)$ and  $\mathfrak{Eul}(M)^{\blacklozenge}$ are torsors on the same group.  
Actually there is a natural bijection of sets  
\begin{equation} \label{euler-iso}
 \mathfrak{Eul}(M) \xrightarrow{\sim} \mathfrak{Eul}(M)^{\blacklozenge}
  \end{equation}
compatible with the action of the group. This is defined by
  $$  V \mapsto [graph(V)]  \in H_3(S(TM), \Z) $$
  where $V $ is considered as a section of the fibration $S(TM) \rightarrow M$.
  %To define  the map (\ref{euler-iso}),   identify the set of nowhere vanishing field on $M$ with the sections of the fibration $S(TM) \rightarrow M$, and map an element $V \in \mathfrak{Eul}(M)$ to the homology class $V_*([M]) \in H_3(S(TM), \Z)$.  
It is immediate to check that  since  $ [graph(V_1)] =  [graph(V_2)]$  if  $[V_1]=[V_2]$. 
%are homologous nowhere vanishing  vector fields of $M$.

In $\mathfrak{Eul}(M)^{\blacklozenge}$ the difference between two Euler structures $[U_1], [U_2] \in \mathfrak{Eul}(M)^{\blacklozenge}$ is obtained by
the difference of the representative
\begin{equation} \label{difference-euler}
[U_2] - [U_1] = [ U_2 - U_1] \in  H_3^-(S(TM), \Z).
\end{equation}
Hence the Chern Class $\mathfrak{c}(U)$ of an Euler Structure $[U] \in \mathfrak{Eul}(M)^{\blacklozenge}$ is given by  
\begin{equation} \label{chern-euler}
\mathfrak{c}([U])= [U] - r_*[U] \in H_3^-(S(TM), \Z) .
%\cong H_1(M,\Z) .
\end{equation}
%where $U$ is considered as an element of  $\mathfrak{Eul}(M)'$.

%----

%This can also obtained from the intersection of the three-chains
%\begin{equation} \label{difference-euler0}
%[U_2] - [U_1] = [\pi_* ( U_2 \cap U_1)] \in H_1(M, \Z),
%\end{equation}
%where we have taken transversal representative for $U_1$ and $U_2$. 

%It is useful to write $\mathfrak{c}(U)$ directly as element of $H_1(M,\Z)$ using
%\begin{equation} \label{chern-class}
%\mathfrak{c}(U)= U^{\bullet} \cap M^{pert} \in H_1(D^2(M),\Z) \cong H_1(M, \Z) 
%\end{equation}
%where $U^{\bullet}$ is the chain on the disk bundle $D^3(M)$ which is invariant by dilatation and $\partial U^{\bullet} = U - r(U)$, $M^{pert}$ is any three chain on $D^3(M)$ close on the $C^0$-topology to the zero section $M \subset D^3(M)$ and transversal to $U^{\bullet}$.
%$V'$ is the graph of any vector field of $M$ transversal to $\overline{U}$.

\subsubsection{Framed Euler Structures}

Now consider the case that $M$ is not compact. We assume that $M$ has a collar at infinity given by $\Sigma \times \R^+$ for some (not necessary connected) surface $\Sigma$, i.e. there exists a compact manifold $M^{in}$ with $\partial M^{in} = \Sigma$ such that 
$$M= M^{in} \sqcup_{\Sigma} ( \R^+ \times  \Sigma ).$$
%$$H^{Fr}(M, \partial X'_{\infty}) $$
%$$M = M^{comp} \sqcup \partial M^{comp} \times \R^+.$$

%---

%Now we allow $M$ that can have a not trivial  boundary $\partial M= \Sigma$.

%Given a frame $\mathbf{fr} \in \mathfrak{Fr}(TM|_{\Sigma})  $ we say that a not vanishing vector field $V$ is compatible with $\mathbf{fr}$ if  there exit $a_1,a_2,a_3 \in \R$ such that
%$$  V|_{\Sigma} = a_1 fr_1 + a_2 fr_2 + a_3 fr_3 $$
%where $\mathbf{fr} = (fr_1,fr_2,fr_3)$.

%A framed Euler Structure is given by a pair $(V , \mathbf{fr})$, where $V$ is not vanishing vector field on $M$ compatible with $\mathbf{fr}$.

%---

A framed Euler Structure consists in a pair $(V , \mathbf{fr})$, where $V$ is not vanishing vector field on $M$ whose restriction to $\Sigma$ is compatible with $\mathbf{fr}$ (see (\ref{compatible-vector-frame})).

We say that two framed Euler Structures $(V_0 , \mathbf{fr}_0), (V_1 , \mathbf{fr}_1)$ are homologous if $V_0, V_1$ are isotopic outside a ball with isotopy compatible with an isotopy of frames $(\mathbf{fr}_t)_{t \in [0,1]}$.

We denote by  $ \mathfrak{FrEul}(M) $ the set of homology classes of Euler Structures and with $ \mathfrak{FrEul}(M, \mathbf{fr}) $ the set of homologous classes of Euler Structures compatible with 
$ \mathbf{fr}$. 
It is easy to show that $ \mathfrak{FrEul}(M, \mathbf{fr}) $ is a torsor on $H_1(M, \Z)$.

%A framed Euler Structure is given by a pair $(V , \mathbf{fr})$, where $V$ is not vanishing vector field on $M$ and $\mathbf{fr} \in \mathfrak{Fr}(TM|_{\Sigma})  $ such that there exit $a_1,a_2,a_3 \in \R$ with
%$$  V|_{\Sigma} = a_1 fr_1 + a_2 fr_2 + a_3 fr_3 $$
%where $\mathbf{fr} = (fr_1,fr_2,fr_3)$.

%We can extended the notion of equivalence of Euler Structure in this context. 

%We consider equivalent two Euler structure $(V_1 , \mathbf{fr})$ and $(V_2 , \mathbf{fr})$ if $V_1$

%Given a frame $\mathbf{fr} \in \mathfrak{Fr}(\Sigma \times \R)^{\R}  $ we say  that a vector field $V$on $M$ is compatible with $\mathbf{fr}$ is the restriction to $ \R^+ \times  \Sigma$ is constant with respect to $\mathbf{fr}$. 

%Given a frame $\mathbf{fr} \in \mathfrak{Fr}(TM|_{\Sigma})  $ we say that an element $S \in C_3(S(TM),  \Z)$ is compatible with $\mathbf{fr}$ if, up to triangulation, 
%$$S \cap S(T( \R^+ \times  \Sigma )) = graph (V) $$
%for some vector field $V$  on $ \R^+ \times  \Sigma$ which is constant with respect to $\mathbf{fr}$. 

Now we define $\mathfrak{FrEul}(M)^{\blacklozenge}$.
We say that an element $S \in C_3(S(TM),  \Z)$ is compatible with $\mathbf{fr} \in \mathfrak{Fr}(TM|_{\Sigma})  $ if, up to triangulation, 
\begin{equation} \label{boundary-vector}
\partial S  = graph (V) ,
\end{equation}
for some  $V \in \Gamma(S(TM)|_{\Sigma} )$  compatible with  $\mathbf{fr}$. Here $graph (V)$ is  considered as a map $\Sigma \rightarrow S(TM)$.

Denote by $C_3(S(TM), \Z| \mathbf{fr})$ the vector space of chains compatible with $\mathbf{fr}$.
Consider pairs $(S, \mathbf{fr})$ with  $ \mathbf{fr} \in  \mathfrak{Fr}(TM|_{\Sigma}) $ and $S \in  C_3(S(TM), \Z | \mathbf{fr})   $. We define the difference of two pairs $(S_1, \mathbf{fr}_1),  (S_0, \mathbf{fr}_0)$ as
$$  (S_1, \mathbf{fr}_1) -  (S_0, \mathbf{fr}_0) = S_1 - S_0 + ( (V_1, \mathbf{fr}_1) - (V_0, \mathbf{fr}_0)  ) \in  H_3(M^{in})$$ RIMUOVERE IN?
%\in C_3(M^{in})/ \partial C_4(M^{in}) ,$$
where $V_0,V_1$ are as in (\ref{boundary-vector}) and we have used (\ref{difference-frames}).

Define  $\mathfrak{FrEul}(M)^{\blacklozenge}$ by
\begin{equation} \label{framed-euler-definition}
 \mathfrak{FrEul}(M)^{\blacklozenge} := \{ (S, \mathbf{fr})|   \mathbf{fr} \in  \mathfrak{Fr}(TM|_{\Sigma}), S \in  C_3(S(TM), , \Z | \mathbf{fr}) ,   \pi_*([S])=[M]  \}/\sim  
 \end{equation}
where now
$pr_*: H_3(S(TM), \partial S(TM),\Z) \rightarrow H_3(M, \partial M, \Z)$.
% is the map induced by the projection map  $\pi: S(TM) \rightarrow M$, 
In (\ref{framed-euler-definition}) we consider the equivalence relation 
$$  (S_0, \mathbf{fr}_0) \sim (S_1, \mathbf{fr}_1) \text{   if     }(S_1, \mathbf{fr}_1) -  (S_0, \mathbf{fr}_0) =0.$$

The  action (\ref{action-pairs}) induces an action of $H^1(\Sigma, \Z )$ on $\mathfrak{FrEul}(M)^{\blacklozenge}$. From formula (\ref{action-pairs-difference}) we obtain 
$$ (2 \gamma) \cdot [ (S, \mathbf{fr}) ]= [ (S + i_*(\gamma), \mathbf{fr})],$$
where $i: \Sigma \rightarrow M$ is the inclusion of the boundary, and on the RHS we have identified $H_3^-(S(TM), \Z) =H_1(M, \Z)$. Here we have also used the canonical isomorphism $H_1(\Sigma, \Z ) = H^1(\Sigma, \Z)$.

The spin structure associated to the frame defines a surjective map 
\begin{equation} \label{Euler-spin}
\mathfrak{FrEul}(M)^{\blacklozenge} \rightarrow Spin(\Sigma) 
\end{equation}
whose fibers are  torsors on $H_1(\Sigma, \Z )$.
%$H_3^-(S(TM), \Z)$.

%that induces a canonical identification of sets
%$$ \mathfrak{Eul}^{Fr}(M)/ H_3^-(S(TM), \Z) \cong Spin(\Sigma) .$$

%The Framed refined Chern Class is defined as
%$$\mathfrak{c}(V,fr)= [ \R V \cap V'  ]  \in H_1(TM,\Z) =H_1(M, \Z) ,$$
%where $V'$ is a section of $TM$ (not necessary nowhere vanishing)  transversal to $\R V \subset TM$  and satisfying the boundary condition 
%$$V'|_{\partial M_c}= fr \cdot (0,1,0).$$

%The difference between two Framed Euler Structures associate to the same frame $\mathbf{fr}$ is obtained adding to the RHS of  (\ref{difference-euler}) the correction term of the boundary: 
%$$ [U_1] -[U_2] = [U_2- U_1 + graph(\tilde{V})] $$
%where $\tilde{V}=(V_t)_{t \in [0,1]}$ is, as in \ref{framed-euler-definition}, compatible with $\mathbf{fr}$ and $ \partial U_1  = graph (V_1),  \partial U_2  = graph (U_2)$. Thus, we have used the frame $\mathbf{fr}$ to pick a particular  class of chains on $(S(TM)|_{\Sigma} )$ whose boundary is $\partial U_1 -  \partial U_2$.

%The  Chern Class of a Framed refined Euler structure is obtained again from (\ref{chern-class}).

The Chern class is defined as in (\ref{chern-euler}). 
In contrast with the compact case, the Chern class  is not necessarily  even . Its parity is determinate by the relative Stiefel–Whitney class
$$ w_2(M , \sigma_{\mathbf{fr}}) = \mathfrak{c}([U,  \mathbf{fr}])  \text{   mod   } 2 H_1(M, \Z) ,$$
where $ \sigma_{\mathbf{fr}}$ is the spin structure on $\Sigma $ defined by $\mathbf{fr}$.

The Chern class changes under the  action of $H_1(M, \Z)$ according the formula
\begin{equation} \label{chern-shift}
 \mathfrak{c}(\gamma \cdot [ U,  \mathbf{fr}]) =  \mathfrak{c}([U,  \mathbf{fr}]) + \gamma .
 \end{equation}

\subsubsection{Framed Three Chains}

Let $Y$ be an oriented  five manifold. Let $\Sigma \subset Y$ be a two dimensional submanifold of $Y$. Assume that the normal bundle $N\Sigma$ of $\Sigma$ is trivial. Denote by
$\hat{Y}$ the real blow-up of $Y$ along $\Sigma$.   $\hat{Y}$ is a five dimensional manifold whose boundary is the spherical bundle $S(N\Sigma)$.

%We define $H_3(\check{Y}, \partial \check{Y})$.

%Set
%$$ FrH_3(\check{Y}):= \{ ( T,fr)  \in  C_3(\check{Y}) \times \mathfrak{Fr}(N \Sigma) | ( fr \cdot (1,0,0) )_*([\Sigma]) = \partial T \}/\sim ,$$
%%where $f \in  Frames(L)$ and $T$ is a submanifold of dimension $3$ of $X_{\infty}$ such that $ IM_1(f) = \partial T $. 
%where 
%$$(fr_0, T_0) \sim (fr_1, T_1) \text{   if     } \exists\tilde{fr} \in Isot(fr_1,fr_2) \text{    s.t.   }  T_1-T_2 - ( \tilde{fr} \cdot (1,0,0) )_*([\Sigma]) \in \partial  C_4(\check{Y}). $$
%%$$ (f_0, T_0) \cong (f_1, T_1) \text{   if  $ \exists  (\tilde{f},Q)$ s.t.} IM_1(\tilde{f})= T_1-T_2 + \partial Q$$
%%where  $ \tilde{f}  $ is an  isotopy between  $f_1$ and $f_2$ and $Q \in C_4(X_{\infty})$.

%---

For $ \mathbf{fr} \in   \mathfrak{Fr}(N \Sigma) $
the definition of $  C_3(\hat{Y} , \Z | \mathbf{fr})   $  is analogous to the definition of  $  C_3(S(TM), \Z | \mathbf{fr})   $.

A Framed Three Chain consists in a pair $(T, \mathbf{fr})$ with  $ \mathbf{fr} \in   \mathfrak{Fr}(N \Sigma) $ and $T \in  C_3(\hat{Y} , \Z | \mathbf{fr})   $. We define the difference of two   Framed Three Chains using the formula
$$  (T_1, \mathbf{fr}_1) -  (T_0, \mathbf{fr}_0) = (T_1 - T_0) + ( (V_1, \mathbf{fr}_1) - (V_0, \mathbf{fr}_0)  ) \in C_3(\hat{Y})/ \partial C_4(\hat{Y}) ,$$
where $V_0,V_1$ are as in (\ref{boundary-vector}) and we have used (\ref{difference-frames}).

Two Framed Three Chains $(T_0, \mathbf{fr}_0) , (T_1, \mathbf{fr}_1)  $ are said homologous if $(T_1, \mathbf{fr}_1) -  (T_0, \mathbf{fr}_0) =0 $. 
Denote by $\mathfrak{FrThree}(Y)$ the set of homology classes of Framed Three Chains.

There is an obvious free action of $H_3(\check{Y}) $  on $\mathfrak{FrThree}(Y)$ by addiction.
The spin structure associated to the frame defines a map of sets
\begin{equation} \label{three-spin}
\mathfrak{FrThree}(Y) \rightarrow  Spin(N \Sigma)
 \end{equation}
whose fibers are torsors on $H_3(\check{Y}) $ .

The action of  $H^1(\Sigma, \Z)$ on  $\mathfrak{FrThree}(Y)$ is inducted by the action of $Maps(\Sigma,SO(2))$ analogously to   (\ref{action-H^1(M, Z )}).

\subsubsection{Frame refined Four chain}

A Framed Four Chain consist in a pair
\begin{equation} \label{pairs-framed-four-chains}
    (K^{in}, \mathbf{fr}) \in C_4(\hat{X}^{in}, \Z) \times \mathfrak{Fr}(\Sigma \times \R)^{\R}
\end{equation}
with the following properties.  
Let $S \in C_3({S(N_*L)})$, $T \in C_3(\check{Y},\Z)$ defined by the relation
$$ \partial K^{in} = S +T  . $$
We require that
%exists $E \in C_3(S^2(TL_{comp}),\Z)$ with $[E] \in  H_3'(S^2(TL_{comp}),\Z) $
$$ pr_* ([S]) = [L^{in}] ,$$
%$$ \partial T = ( fr \cdot (1,0,0) )_*([\Sigma]) , $$
%$$ \partial S = - ( fr \cdot (1,0,0) )_*([\Sigma]) .$$
$$ \partial T =  graph (V)  , $$
$$ \partial S = - graph (V) ,$$ 
for some  $V \in \Gamma(S(TL)|_{\Sigma} )$ compatible with  $\mathbf{fr}$ .

%We define $(K,fr) \sim (K', fr')$ if there exists $\tilde{fr} \in Isot(fr_1,fr_2) $, $P \in C_4( \hat{X}_{\infty})$, $Q \in C_4(S(TL_{c}),\Z) $, $R \in C_5(\hat{X}_{c})$ such that 
%$$ T - T' -  ( \tilde{fr} \cdot (1,0,0) )_*([\Sigma])  = \partial P , $$
%$$ S- S' +   ( \tilde{fr} \cdot (1,0,0) )_*([\Sigma])  = \partial Q ,$$
%$$ K_{c} - K_{c}' -  P - Q  = \partial R.$$

We say $(K^{in},\mathbf{fr}) \sim ({K^{in}}', \mathbf{fr}')$ if there exists $P \in C_4( \hat{Y})$, $Q \in C_4(S(TL^{in}),\Z) $, $R \in C_5(\hat{X}^{in})$ such that 
$$ T - T' -  ( (\mathbf{fr}, V) - (\mathbf{fr}', V') ) = \partial P , $$
$$ S- S' +   ((\mathbf{fr}, V) - (\mathbf{fr}', V'))  = \partial Q ,$$
$$ K^{in} - {K^{in}}' -  P - Q  = \partial R.$$

Let  $\mathfrak{FrFour}(\hat{X})$ be the set of the equivalence classes of Framed Four Chains.
%pairs (\ref{pairs-framed-four-chains}).

%$$\partial T_{\infty} = (fr \ast e_1)(\Sigma)$$
%$$ \partial T_{comp} - T_{\infty} \in C_3(\partial X_{comp})  $$
%Isotopies:
%$$  T_{\infty} - T_{\infty}' -  (\tilde{fr} \ast e_1)(\Sigma)  \in \partial C_4(X_{\infty})$$
%$$ T_{comp} - T_{comp}' -   \partial (\tilde{fr} \ast e_1)(\Sigma) \in \partial C_5(X_{comp})$$

A framed Four Chain in particular determine a framed Euler Structure of $L$
$$\mathfrak{FrFour}(X,L)  \rightarrow \mathfrak{FrEul}(L) $$
$$   (K^{in} , \mathbf{fr}) \mapsto (S,\mathbf{fr}) .$$
It also determine a framed three chain of $Y$: 
$$ \mathfrak{FrFour}(\hat{X})  \rightarrow  \mathfrak{FrThree}(\hat{Y}) $$
$$   (K^{in} , \mathbf{fr}) \mapsto (T,\mathbf{fr}) .$$

Consider the fiber product of the two maps above  using  (\ref{Euler-spin}) and (\ref{three-spin})
\begin{equation} \label{boundary-four-chain}
 \mathfrak{FrFour}(X)  \rightarrow \mathfrak{FrEul}(L) \times_{Spin(\Sigma)}    \mathfrak{FrThree}(Y) .
  \end{equation}
The not empty fibers of  (\ref{boundary-four-chain}) are torsors on $\text{Ker} \{ H_4(X) \rightarrow H_1(L) \}$.
The image of  (\ref{boundary-four-chain})  is identified with the pre-image of zero of the map
$$\mathfrak{FrEul}(L) \times_{Spin(\Sigma)}   \mathfrak{FrThree}(Y)  \rightarrow H_3( \hat{X}).$$
$$ [(S_1, fr_1)]  \times [(T_2, fr_2)] \mapsto S_1- T_2 -( (\mathbf{fr}_1, V_1) - (\mathbf{fr}_2, V_2)) .$$

\subsection{Area zero annulus and Euler Structures}

We now compare the Euler structure associated to the Four chain with the $MC$-cycles associated  to the perturbation of the area zero annulus. 
In order to make this comparison the following description of $\mathfrak{Eul}(M)$ is useful.
Consider the set of  pairs $(Z_0 , Z_1)$ where
\begin{itemize}
\item
 $Z_0$ is a three chain on $M \times M$ close in the $C^0$ topology to $Diag_M$ and transversal to $Diag_M$;
\item $Z_1$ is a one-chain on $M$, with
$$\partial Z_1 + Z_0 \cap Diag_M=0,$$
where $Z_0 \cap Diag_M$ is considered as a chain on $Diag_M \approx M$.  
\end{itemize}
For the pairs $(Z_0,Z_1)$ there is an obvious notion of isotopy .
Let $\mathfrak{Eul}(M)^{\blacktriangle}$  be the set of pairs  $(Z_0 , Z_1)$ modulo isotopy.  The group $H_1(M, \Z)$ acts  on $\mathfrak{Eul}(M)^{\blacktriangle}$ by addiction on $Z_1$.
With this action, $\mathfrak{Eul}(M)^{\blacktriangle}$ is a torsor on $H_1(M, \Z)$

Using a small tubular neighbourhood of $Diag_M$, up to isotopy, we can consider $Z_0$ as a three chain on $TM$ close to the zero section. 

If $V$ is a vector field on $M$ with  small norm,  the map
$$V \mapsto (graph(V), 0)$$
defines an element of $\mathfrak{Eul}(M)^{\blacktriangle}$.
%where $graph(V)^{\spadesuit}$ is the chain on $M \times M$ obtained from $graph{V}$ using a small tubular neighbourhood of $Diag_M$.
Since the map is compatible with the action of $H_1(M,\Z)$, it yields to an isomorphism  
 $$\mathfrak{Eul}(M) \cong \mathfrak{Eul}(M)^{\blacktriangle}.$$
%  since both the sides are torsors on $H_1(M,\Z)$.

% ---

% We can define a map 
% $$ \mathfrak{Eul}(M)^{\blacktriangle} \mapsto  \mathfrak{Eul}(M)^{\blacklozenge}$$
% $$  (Z_0 , Z_1) \mapsto trans(Z_0) + (\pi)^*(Z_1) $$
% which is an isomorphism.

% ---

Now let $Ann0$  be the set of area zero decorated graphs $G$ with one only component and $\Sigma_G$ be an annulus with one boundary marked points. The elements of $Ann0$  are characterized as follows:
\begin{enumerate} \label{unstable-cases}
\item $G_0$: $|V_c|=1, |D_c|=0$,  $|H_c| =3 $,  $|E^{in}_c|=1$;  \label{disk-edge}
\item  $G_1$: $|V_c|=2, |D_c|=0$ , $|H_c| = 1 $;   \label{annelus} 
\item    $G_2$: $|V_c|=1, |D_c|=1$ ,$|H_c| = 1 $.   \label{annelus-deg}
\end{enumerate}
%Denote by $G_0, G_1, G_2$ these graphs respecti
Denote by $\mathcal{C}_{Ann0}$ the restriction of the $MC$-chain complex to the graphs $Ann0$ and by $\mathcal{Z}_{Ann0}$ the vector space of the corresponding $MC$-cycles.
$\mathcal{Z}_{Ann0}/isotopy $ is a torsor on $H_1(M, \Z)$, where $H_1(M, \Z)$ acts by addiction on $Z_{G_1}$.

\begin{lemma} 
\begin{equation}    \label{isotopy-classes-Ann}
    \mathcal{Z}_{Ann0}/isotopy \cong \mathfrak{Eul}(M)^{\blacktriangle}.
\end{equation}
\end{lemma}
\begin{proof}
 We already observed that both the sides are torsors on $H_1(M, \Z)$. To prove the lemma we need to define a map between the two sets  compatible with $H_1(M, \Z)$.
 
From an element  $(Z_{G_0}, Z_{G_1}, Z_{G_2}) \in \mathcal{Z}_{Ann0}$ we obtain an element of $\mathfrak{Eul}(M)^{\blacktriangle}$ as follows.
Let $e $ be the only internal edge of $G_0$, and let  $pr_e : M^{H(G_0)} \rightarrow M^{e} \cong M^2$ be associated  projection.
Since $\delta_e(pr_e(Z_{G_0}))$ is close in the $C^0$-topology to $\delta_e Z_{G_0}$ and $\delta_e Z_{G_0} + \partial Z_{G_1} + \partial Z_{G_2}=0$, there exists 
$Z_{G_1}'$ close in the $C^0$-topology to $Z_{G_1}$ such that 
$$\delta_e(pr_e Z_{G_0}) + \partial Z_{G_1}' + \partial Z_{G_2}=0.$$ 
Hence 
\begin{equation}  \label{Ann-Euler-map0}
(pr_e(Z_{G_0}), Z_{G_1}' + Z_{G_2}) \in \mathfrak{Eul}(M)^{\blacktriangle}.
\end{equation}
%$$(Z_{G_0}, Z_{G_1}, Z_{G_2}) \mapsto (pr(Z_{G_0}), Z_{G_1}' + Z_{G_2})$$

%----

%We obtain a map
%\begin{equation} \label{Ann-Euler-map}
% \mathcal{Z}_{Ann0}/isotopy \rightarrow   \mathfrak{Eul}(M)^{\blacktriangle}.
% \end{equation}
%which is an isomorphism since both the sides are torsors on $H_1(M,\Z)$.

\end{proof}

Now consider triples  
\begin{equation} \label{triple-euler}
    (Z_0 , Z_1, Z_2)
\end{equation}
where
\begin{itemize}
\item
 $Z_0 $ and $Z_2$ are three chains on $M \times M$ close in the $C^0$ topology to $Diag_M$ and transversal to $Diag_M$;
\item $Z_1$ is a one-chain on $M$, with
$$\partial Z_1 + (Z_0 - Z_2) \cap Diag_M=0.$$
%where $Z_0 \cap Diag_M$ is considered as a chain on $Diag_M \approx M$.  
\end{itemize}
There is an obvious notion of isotopy for the pairs $(Z_0,Z_1, Z_2)$.  

%To a triple $(Z_0 , Z_1, Z_2)$ we can associate an homology class $[(Z_0 , Z_1, Z_2)] \in H_1(M,\Z)$ which depends only on its isotopy class.
%To an isotopy classes agree with $H_1(M,\Z)$. 
%Actually, there is a unique map from the class of isotopies to  $H_1(M,\Z)$ such that
%$$(Z_0,Z_1, Z_2)  \rightarrow  [Z_1]  \text{   if    } Z_0 = Z_1 .$$
%$[(Z_0 , Z_1, Z_2)]$ is characterized by the following properties:
To a triple $(Z_0 , Z_1, Z_2)$ we can associate an homology class $[(Z_0 , Z_1, Z_2)] \in H_1(M,\Z)$
characterized by the following properties:
\begin{itemize}
    \item it is invariant by isotopy;
    \item $[(Z_0,Z_1, Z_2)]  =  [Z_1]$    if    $ Z_0 = Z_1 $.
\end{itemize}
%$$[(Z_0,Z_1, Z_2)]  =  [Z_1]  \text{   if    } Z_0 = Z_2 .$$

\subsubsection{Perturbation of the moduli space}

Let  $(\overline{\mathcal{M}}_{G})_{G \in Ann0}$ be the collection of Kuranishi spaces associated to the graphs $Ann0$.
They come with  the following identification of Kuranishi spaces
\begin{equation} \label{identification-kuranishi-ann}
 \partial \overline{\mathcal{M}}_{G_1} \cong  \delta_e( \overline{\mathcal{M}}_{G_0})  \sqcup   \delta_{e'}( \overline{\mathcal{M}}_{G_2}) 
 \end{equation}
where we denote by $e $ the only internal edge of $G_0$ and by $e'$ the only element of $D(G_2)$.

Let  
\begin{equation} \label{perturbation-ann0}
(\mathfrak{s}_{G}^+)_{G \in Ann0}
\end{equation}
 be a collection of perturbations of $(\overline{\mathcal{M}}_{G})_{G \in Ann0}$ compatible with the identification of Kuranishi spaces  (\ref{identification-kuranishi-ann}). 
(\ref{perturbation-ann0}) is obtained from the inductive argument of \cite{OGW3}, which 
%s in \cite{OGW3}, let $ \overline{\mathcal{M}}_{G}^K$ be the moduli space of associated multi-curves. 
%We need to consider perturbations of   $ \{  \overline{\mathcal{M}}_{G}^K \}_{G \in Ann0}$ satisfying the constrains defined in \cite{OGW3}. These perturbations are constructed in \cite{OGW3} using an inductive argument that  
in the particular case of $Ann0$ reduces to the following two steps:
\begin{itemize}
\item choice a transversal perturbation on $\overline{\mathcal{M}}_{G_0}$ and $\overline{\mathcal{M}}_{G_2}$, which are transversal also when restricted to the RHS of (\ref{identification-kuranishi-ann}) ;
\item choice an extension to $\overline{\mathcal{M}}_{G_1}$ of the perturbation defined on its boundary by (\ref{identification-kuranishi-ann}) . 
\end{itemize}

%Assume also that $\mathfrak{s}_{G_2}^+$ is transversal when restricted to $\text{ev}^{-1}(L \times L)$.

 %Let $G^{eu}$ be the decorated graph given by one only component that is an area zero disk with one internal and one boundary marked point
%Thus $Z_{G^{eu}}$ is a closed three chain on $X \times L$ that is close in $C^0$-topology to the diagonal of $L$ $diag_L \subset L \times L \subset X \times L$.

%\subsubsection{Boundary condition}

%Use $\mathbf{fr}$ to give a local $X \times L \cong \R^3$ we can consider the condition on  $\Sigma \times \R$.

%Recall that the obstruction bundle $E_u$ on a pseudo-holomorphic $(\Sigma,u)$ of the Kuranishi Structure consists of a suitable sub-vector space  of $C^{\infty}(\Sigma, u^*(TX) \otimes \Lambda^{0,1}(\Sigma))$. If $u$ is a constant map $L^p(\Sigma, u^*(TX) \otimes \Omega^{0,1}(\Sigma))= T_{Im(u)}X \otimes L^p(\Sigma,  \Omega^{0,1}(\Sigma))$
%where $Im(u) \in X$ is the point image of $u$.  

In the not compact case we need to consider a boundary condition at infinity. 
We are going to define a particular type of boundary condition depending on the choice of  a frame $\mathbf{fr} \in \mathfrak{Fr}(\R \times \Sigma)^{\R}$. 
Using $\mathbf{fr}$ and a complex structure, we can identify $T_{p}X$ in (\ref{obstruction-bundle0}) with $\R^3 \times \R^3$.
% for different $u$ belonging to $\R^+ \times \Sigma$. 
We require  that on $ \R^+ \times \Sigma $
\begin{itemize}
\item the obstruction bundle $E_{\mathbf{p}}$ does not depend on $p$,
\item  the perturbation $\mathfrak{s}_{\mathbf{p}}$  does not depend on $p$,
\item $K^*$ is compatible with $\mathbf{fr}$.
\end{itemize}

%The associate chains $(Z_{G_0}, Z_{G_1}, Z_{G_2}^+)$ are  small perturbations of the diagonal, which is constant, when we identify a small neighbourhood with  $Diag \times \R^{H(G)}$, up to an error which is of magnitude higher, and therefore does not mutter up to isotopy.

%$(Z_{G_0}, Z_{G_1}, Z_{G_2}^+)$ are costant

%ssume that (\ref{perturbation-ann0}) is constant. Thus 
% $(Z_{G_0}, Z_{G_1}, Z_{G_2}^+)$ is constant. 
 
%ASSUMPTION: $Z_{G^{eu}}$ is invariant by translation on $\Sigma \times \R$ 
  
%\begin{equation}  \label{boundary-condition-euler}
%Z_{G_2}^+= fr \ast (0,1,0).
%\end{equation}  
 
%---

%Denote by $(Z_{G_0}, Z_{G_1}, Z_{G_2}^+)$ the corresponding chains. 

\subsubsection{Euler Structure from area zero Annulus}

Let  $(Z_{G_0}, Z_{G_1}, Z_{G_2}^+)$ be  the collection of chains associated to the perturbation defined above.
In particular $Z_{G_2}^+ \in C_3(L \times X)$ is transversal to $L \times L$ and it is a close in the $C^0$-topology to $Diag_L \subset L \times X$.
%The collection of chain $ (Z_{G_0}, Z_{G_1}, Z_{G_2}^+) $ satisfies
We have the relation
 $$\delta_e Z_{G_0} + \partial Z_{G_1} + \delta_{e'} Z_{G_2}^+=0$$
where $\delta_{e'} Z_{G_2}^+ = pr_1 ( Z_{G_2}^+ \cap (L \times L))$ with $pr_1: L \times X \rightarrow L$ is the projection on the first factor (we label the internal puncture of $G_2$ with $e'$).

To $ (Z_{G_0}, Z_{G_1}, Z_{G_2}^+) $ we can associate a triple $(Z_0 , Z_1, Z_2)$ in the sense of (\ref{triple-euler})
\begin{equation} \label{triple-associated}
    (Z_{G_0}, Z_{G_1}, Z_{G_2}^+) \leadsto (Z_0 , Z_1, Z_2).
\end{equation}
 To define $(Z_0 , Z_1, Z_2)$  we use a tubular neighborhood of $L$ inside $X$ in order to identify $Z_{G_2}^+$ with the graph of a vector field. We then perturb $Z_{G_1}$ as in the proof of Lemma \ref{isotopy-classes-Ann} to achieve the cyclic condition. 

\begin{lemma} \label{triple-zero}
    The homology class of the triple $(Z_0 , Z_1, Z_2)$ associated to $ (Z_{G_0}, Z_{G_1}, Z_{G_2}^+) $ is trivial.
\end{lemma}
\begin{proof}

The Lemma can be considered as the analogous for annulus of the well know computation of the orbifold Euler Class of the obstruction bundle for area zeros torus with one marked point. In the case of torus  we are interested to the orbibundle over $\overline{\mathcal{M}}_{0 ,1} \times X$ whose fiber at the point $((T^2, J_{T^2}),x)$ is the tensor product  $H^{0,1}(T^2, J_{T^2}) \otimes TX$. The first factor is the Hodge Bundle of the torus.

-In our case,  the Hodge bundle of the annulus is trivial and the Euler class of $TL$ vanishes.
$(\overline{\mathcal{M}}_{(0,2) ,(0, (0,1))}) \times L$

To prove the Lemma we consider a particular type of perturbation of the moduli space $(\overline{\mathcal{M}}_{G})_{G \in Ann0}$.

Pick a one dimensional vector sub bundle $E_{\Sigma}'$ of $\Lambda^{0,1} (\Sigma)$ over the moduli space of $(\Sigma,J_{\Sigma})$ which is complementary to $Im( \overline(\partial))$. 
We constrain the obstruction bundle (\ref{obstruction-bundle0}) requiring that $E_{\mathbf{p}}$ on the point  $\mathbf{p}= (u, \Sigma)$ contains the tensor product
\begin{equation} \label{obstruction-bundle01}
E_{\mathbf{p}} \supseteq    T_{p}X \otimes E_{\Sigma}' .
\end{equation}
%for some $E_{\Sigma}' \subset \Lambda^{0,1} (\Sigma)$ depending only on $ \Sigma $ and   independent on $u$. We assume that $E_{\Sigma}'$ define a vector bundle over $\Sigma_{G_1}$ TRANSVERSAL.

Pick a not vanishing vector field $V$ on $L$ and a transversal section  $\mathfrak{s}'$ of   $E_{\Sigma}'$. 
On $\overline{\mathcal{M}}_{G_1}$ consider the perturbation of the Kuranishi structure given by
\begin{equation} \label{tensor-product-perturbation}
    V \otimes \mathfrak{s}'
\end{equation}
%with $\mathfrak{s}'$ independent on $u$.

Requiring the compatibility with  (\ref{identification-kuranishi-ann}), (\ref{tensor-product-perturbation}) defines a perturbation of the Kuranishi structure on $\delta_e \overline{\mathcal{M}}_{G_0}$. From $\mathfrak{s}_{G_0}^{-1}(0)$ we can define an Euler structure $V_0$, up to small isotopy.  The homology class of $V_0$ agrees with the homology class $ V$ or its opposite $opp(V)$. 
In the same way, from $\mathfrak{s}_{G_2}^{-1}(0)$  we obtain  an Euler structure $V_2$  whose homology classes agrees with $ V$ or $opp(V)$.
The homology class $[(Z_0 , Z_1, Z_2)]$ of the triple (\ref{triple-associated}) is given by 
$$[(Z_0 , Z_1, Z_2)] = [V_0]-[V_2] .$$
Since $[(Z_0 , Z_1, Z_2)]$ has to be independent of the perturbation (and hence of $V$) we conclude that $V_0= V_2$ (note that if $V_0$ is opposite to $V_2$, $V_0-V_2$  depends on $V$). 
\end{proof}

Set 
%$Z_{G_2} =Z_{G_2}^+ \cap K $.
\begin{equation} \label{intersect-K}
Z_{G_2}=  Z_{G_2}^+ \times_X K.
\end{equation}
The dependence of $Z_{Ann_0}:=(Z_{G_0}, Z_{G_1}, Z_{G_2}) \in \mathcal{Z}_{Ann0}$ on $K$ is encoded in  (\ref{intersect-K}). 

%In the not compact case we require that $K^*$ is compatible with $\mathbf{fr}$. Then generically this implies that $Z_{G_2}=0$.
In the not compact case,  on the region $\R^+  \times \Sigma$, from the construction above we obtain that $Z_{G_1}= Z_{G_2}=0$ and $Z_{G_0}$ is equal to
$ (\R^+  \times \Sigma ) \times v$ for some $v \in \R^3$ with small norm, up to an higher order error which can be canceled up to isotopy .

The four chain $K$ defines an element $ [U_K] \in \mathfrak{Eul}(M)^{\blacklozenge}$  as follows. Assuming transversality  between $L$ and $K$ ,  we can define a four chain $\hat{K} $ on $\hat{X} $ and  set 
$$  U_K=  \partial  \hat{K}.$$

In the not compact case we  define   $ [U_K] \in \mathfrak{FrEul}(M)^{\blacklozenge}$ in an analogous way.

\begin{proposition}
%Let $Z_{Ann_0} \in   \mathcal{Z}_{Ann_0} $ defined using $Z_{G_i}= \mathfrak{s}^{-1}_{G_i}(0)$, for $i=0,1,2$. 
The isotopy class  of $Z_{Ann0}$ coincides with the homology class   $[U_K]$ using (\ref{isotopy-classes-Ann}) .
\end{proposition}
\begin{proof}

First note that, up to isotopy,  $Z_{Ann0}$  depends on $K$ only throughout $[U_K]$.
%$[K^{\dagger}]:= \partial K  \in \mathfrak{Eul}(M)^{\blacklozenge} $. 
Moreover we can define   $Z_{Ann0}$ for any  $U \in \mathfrak{Eul}(M)^{\blacklozenge}$ not necessarily  associated to a four chain $K$.   
%EULER STURTURE OR HOMOLOGY CLASS OF K

Let $\Delta$ be the difference between  (\ref{Ann-Euler-map0}) and $[U_K]$ as Euler Structures.

If we shift $[U_K]$ by an element of $H_1(L,\Z)$ the corresponding element (\ref{Ann-Euler-map0}) shifts by the same element. Hence $\Delta$ does not  depend on $U_K$. 

It is easy to check  that $\Delta$ agrees with the homology class of the triple (\ref{triple-euler}) associated to $(Z_{G_0}, Z_{G_1}, Z_{G_2}^+)$.

%---

%In the not-compact case, $\Delta$ does not depend on the frame $\mathbf{fr}$ which we have used to define the boundary condition. Assume that, on $\R^+ \times \Sigma$, $ K^*$ coincide with $graph(fr_1)$, so that $K^*$ is compatible with $\gamma \cdot \mathbf{fr}$ for each $\gamma \in H_1(\Sigma, \Z)$  ( we have used something similar in (\ref{chern-shift})). We can define an isotopy between $Z_{Ann_0}$ associated to $\mathbf{fr}$ and $Z_{Ann_0}$ associate to $\gamma \cdot \mathbf{fr}$ which does not intersect the graph of $-\mathbf{fr}_1$.  This implies that $\Delta - \gamma \cdot \Delta =0$.Hence $\Delta$ depends only on the topology of the obstruction bundle. Therefore $\Delta$ needs to be computable in terms of characteristic classes of $TL$ (this step is analogous to the well known case  of area zero torus  ). Since the tangent bundle of $L$ is trivial, $\Delta$ vanishes.

%---

%ALTERNATIVE PROOF:

% Projecting on $X$ and using a small tubular neighborhood  of $L$ in $X$ , from  $Z_{G_2}^+$ we obtain a chain $trans(Z_{G_2}^+)$ on $L$ which is close in the $C^0$-topology to $Diag_L$.  Moreover
%$ \delta_v Z_{G_2}^+ $ is close in the $C^0$-topology to $-\delta_e(trans(Z_{G_2}^+))$. Thus, there exists $Z_{G_1}'$ close in the $C^0$-topology to  $Z_{G_1}$ such that
%$$  pr_e Z_{G_0} + \partial Z_{G_1}' -\delta_e(trans(Z_{G_2}^+))=0 .$$

%NOTE the difference of sign follows because $N(L\times L \subset L \times X)$ and $N(L \subset X)$ have opposite orientation.

%Note that $ Z_{G_2}^+$ defines an Euler structure.

\end{proof}

%Using formula (\ref{difference-euler}) , $Z_{G_2}=  Z_{G_2}^+ \cap \partial K$ is the difference between the Euler structure    $ trans(Z_{G_2}^+)$  and  $\partial K$.  Hence the difference between the Euler Structures $Z_{Ann_0}$ and $\partial K$     is given by
%$$ pr(Z_{G_0}) -   trans(Z_{G_2}^+) + Z_{G_1}'  \in H_1(L, \Z) .$$
%Note that this element of $H_1(L, \Z)$ is independent on $K$ and depends only on the three manifold $L$. Since the tangent bundle of $L$ is trivial, it has to be necessarily zero.

%The Proposition follows.

%---

%Let $Z_{Ann_0}^{\blacklozenge}$ the element of  $\mathfrak{Eul}(M)^{\blacklozenge} $ corresponding to $Z_{Ann_0}$. Thus we have
%$$ Z_{Ann_0}^{\blacklozenge} = graph(pr_e(Z_{G_e}))^* +  (\pi)^*(Z_{G_1}')+ (\pi)^*(Z_{G_2} ) .$$

%According formula  (\ref{difference-euler}) , the difference between $Z_{Ann_0}^{\blacklozenge}$ and $K^*$ is given by
%\begin{multline}
% \pi_*(  Z_{Ann_0}^{\blacklozenge} \cap K^*) = \pi_* (graph(pr_e(Z_{G_e}))^* \cap K^* ) + Z_{G_1}' + \pi_*(Z_{G_2}^{+*} \cap K^*)= \\
%(\pi^* ( graph(pr_e(Z_{G_e}))^* )  + \pi^*( Z_{G_1}' )+  \pi^*( \pi_*(Z_{G_2}^{+*} )) )  \cap K^*
%\end{multline}

%It is easy to check that the RHS does not depend on $K^*$. Hence it depends only on $L$ and therefore it is equal to zero since the tangent bundle of $L$ is trivial.

\subsubsection{Abelianization and Chern Class}

In \cite{OGW3} we considered the abelianization of a $MC$-cycle. 
%As in \cite{OGW3} from $Z$ we can construct the $MC$-cycle $Z^{ab}$ using the abelianization map, which
For $Ann_0$ it reduces to 
$$Z^{ab}_{G_0}= \frac{1}{2}(Z_{G_0} -  opp(Z_{G_0})), Z^{ab}_{G_1}=Z_{G_1} , Z^{ab}_{G_2}= Z_{G_2} ,$$
where $opp(Z_{G_0})$ is the chain obtained switching the half edges of $e$. 

Using an isotopy to contract $Z^{ab}_{G_0}$ to zero, we can consider  $Z^{ab}$, up to isotopy, as an  element   $[Z_{G_1}+ Z_{G_2}] \in H_1(L, \Z)$.

%\begin{lemma}
%With the boundary condition (\ref{boundary-condition-euler} ) $2 \cdot Z_{Ann0}$ is a representative of the Chern class of the framed Euler structure $(\partial K,fr)$
%$$2 \cdot Z_{Ann0} = \mathfrak{c}(\partial K, fr) .$$
%\end{lemma}
%\begin{proof}

%Use an isotopy to contract $Z^{ab}_{G_0}$ to zero. Hence to $Z^{ab}$ it is associate the element  $[Z_{G_1}+ Z_{G_2}] \in H_1(L, \Z)$.  

%\end{proof}

\begin{lemma}
$2 Z^{ab}_{Ann0}$ is
the Chern class of the Euler structure associate to $Z_{Ann0} $. 
\end{lemma}
\begin{proof}
Given $Z \in   \mathcal{Z}_{Ann0} $,  we can define the opposite cycle  $Z^{\ominus} \in  \mathcal{Z}_{Ann0}$ as
$$Z^{\ominus}_{G_0}= opp(Z_{G_0}),  Z^{\ominus}_{G_1}=-Z_{G_1} , Z^{\ominus}_{G_2}=- Z_{G_2} .$$ 
%Then
%$$  \frak{1}{2}( Z+Z^{\ominus} ) $$
%is the $MC$-cycles .
The Euler Structure corresponding to $Z^{\ominus}$ is the opposite of the Euler Structure corresponding to $Z$. 
Since $2 Z^{ab} = Z -   Z^{\ominus}  $ the Lemma follows.
\end{proof}

\section{Kontesevich-Soibelman algebra from Multi-Curve-Homology}

%$H_1(\Sigma, \Z) $ is a lattice with endowed with skew-symmetric pairing  $ \langle ,  \rangle$ provided by the intersection pairing. 

%We start remarking that it is possible to define the product of two coherent cycles  if their support is disjoint. 

%\begin{remark} Given two coherent cycles $\text{exp}((w_1 , w_1'  ))$ and $\text{exp}(( w_2 , w_2' ))$, if the support of   
%$w_1 $ and $ w_2  $ are disjoint, we define the product as 
%\begin{equation} \label{product-coherent}
%\text{exp}((w_1 , w_1'  )) \times \text{exp}(( w_2 , w_2' ))  = \text{exp}((w_1+ w_2 , w_1'  +w_2' )).
% \end{equation}
%\end{remark}

\subsection{The case $\Sigma \times \R$}

%Let 
%$$  MCH(M, \mathfrak{c})=   \bigoplus_{\gamma \in  H_1(M,\Z) } MCH(M, \gamma |  \mathfrak{c}  )$$

%\begin{lemma}
%For each $\gamma$ and $ \mathfrak{c}$,
%$MCH([0,1] \times \Sigma, \gamma, \mathfrak{c})^{\diamond, \diamond}$ is free module of rank one.
%\end{lemma}

Particularly important is the case of the cylinder $M = \Sigma \times [0,1] $. In this case
the vector space $MCH(\R \times \Sigma )$ has a natural product: 
$$MCH([0,1] \times \Sigma, \gamma_1, \mathfrak{c}_1)^{\diamond}  \times  MCH([0,1] \times \Sigma, \gamma_2, \mathfrak{c}_2)^{\diamond} \rightarrow MCH([0,1] \times \Sigma, \gamma_1 + \gamma_2, \mathfrak{c}_1+ \mathfrak{c}_2)^{\diamond} $$
%the product of two elements of $MCH(M)^{\diamond}$ 
 obtained  from    (\ref{product-coherent}) and the isomorphism 
$$(\Sigma \times [0,1]) \sqcup_{\Sigma} (\Sigma \times [0,1]) \cong \Sigma \times [0,1] $$
obtained by gluing the boundary component $\Sigma \times \{ 0 \}$ of the first cylinder to the boundary component  $\Sigma \times \{ 1 \}$ of the second cylinder, i.e., we place the first factor on the top of the second. 

In particular the terms with $ \mathfrak{c}  =0$ define an algebra
$$  MCH([0,1] \times \Sigma, \mathfrak{c}=0)=   \bigoplus_{\gamma \in  H_1(\Sigma,\Z) } MCH([0,1] \times \Sigma, \gamma |  \mathfrak{c}=0  ).$$

Consider the Kontsevich-Soibelman algebra associated to 
the homology group  $ H_1(\Sigma,\Z) $ endowed with the standard intersection  form $\langle \bullet , \bullet \rangle$:
%$$  \mathfrak{g} = \bigoplus_{\gamma \in  H_1(\Sigma,\Z) } \Q [q^{\frac{1}{2}}, q^{-\frac{1}{2}}] \hat{e}_{\gamma} $$
$$ \hat{ \mathfrak{g} } = \bigoplus_{\gamma \in  H_1(\Sigma,\Z) } \Q [[g_s]] \hat{e}_{\gamma} $$
\begin{equation} \label{KS-algebra-product}
    \hat{e}_{\gamma_1}  \hat{e}_{\gamma_2} = 
    %(-1)^{\langle \gamma_1 , \gamma_2 \rangle}
    q^{ \frac{\langle \gamma_1 , \gamma_2 \rangle}{2}}  \hat{e}_{\gamma_1 + \gamma_2} \text{   for    }  \gamma_1, \gamma_2 \in H_1(\Sigma,\Z) 
\end{equation}
where
%We show that there exists a natural representation on $\mathfrak{H}$ of the quantum KS Kontsevich-Soibelman algebra $\mathfrak{g}$ with
$$     q^{\frac{1}{2}} = - e^{ \frac{g_s}{2}}       .$$

Fix a spin structure on $\Sigma$, and let $\sigma$ its quadratic differential.
Consider the linear map 
$$ \hat{ \mathfrak{g} } \rightarrow  MCH([0,1] \times \Sigma, \mathfrak{c}=0)   $$
\begin{equation} \label{KS-alebra-coh-map}
    \hat{e}_{\gamma} \mapsto  \sigma(\gamma) \text{exp}(( w_{\gamma} , w_{\gamma}' )) 
\end{equation}
 % $\gamma \in H_1(\Sigma,\Z)$ we can associate the coherent cycle $\hat{e}_{\gamma} \in MCH([0,1] \times \Sigma, \mathfrak{c}=0)$ defined as
where $w_{\gamma}$ is any representative of the homology class $\gamma$ with support on $\{t  \} \times \Sigma$ for some $t \in (0,1)$, and $w_{\gamma}'$ is a translation of $w_{\gamma}$ to $\{t'  \} \times \Sigma$ for some $t' \neq t$ close to $t$. Observe that in this case $w^{ann}=0$.

\begin{proposition} \label{KS-algebra-iso-proposition}
The linear map above defines an isomorphism of algebras
$$  MCH(\Sigma \times \R)^{\diamond, \diamond} =\hat{\mathfrak{g}} .$$
\end{proposition}
\begin{proof}
Let $w_1$, $w_1'$, $w_2$, $w_2'$  with support on  $\{ t_1\} \times  \Sigma, \{ t_1'\} \times  \Sigma,  \{ t_2 \} \times  \Sigma,  \{ t_2'\} \times  \Sigma$ respectively, for some  $t_1,t_1',t_2, t_2' \in [0,1]$.  By definition we have
$$
\text{exp}((w_1+ w_2 , w_1'  +w_2' ))=
\begin{cases}
  \hat{e}'_{\gamma_2} \circ \hat{e}'_{\gamma_1} \text{  if  } \text{max} \{ t_2,t_2' \} < \text{min} \{ t_1,t_1'  \} \\
  \hat{e}'_{\gamma_2 + \gamma_1} \text{  if  } \text{max} \{ t_1',t_2' \} < \text{min} \{ t_1,t_2  \} .
\end{cases} 
$$
Thus $\hat{e}'_{\gamma_2} \circ \hat{e}'_{\gamma_1}$ and $ \hat{e}'_{\gamma_2 + \gamma_1}$ are related by an isotopy which interchanges the order of $t_1'$ and $t_2$. 
%The result follows from the definition of the $\Z$-action  of Lemma \ref{action}. 
%$P_{\gamma + \gamma_1 + \gamma_2}$
The Proposition follows from the fact that under this isotopy we have
$$(w_1+ w_2 , w_1'  +w_2' ) \leadsto (w_1+ w_2 , w_1'  +w_2' ) + g_s \frac{\langle \gamma_1 , \gamma_2 \rangle}{2}.$$

\end{proof}

There is another version $KS$-algebra 
%of $\mathfrak{g}$ 
which plays an important rule in this paper. Set
$$  MCH([0,1] \times \Sigma)^{ann}=   \bigoplus_{\mathfrak{c} \in  H_1(\Sigma,\Z) } MCH([0,1] \times \Sigma, 0 |  \mathfrak{c}  )^{\diamond, \diamond},$$
which is equipped with an algebra structure inducted by the gluing as before.  We denote with $ \hat{\mathfrak{g}}^{ann}$ this algebra. 

Denote with $\hat{\mathfrak{g}}'$ the $KS$-algebra without the sign $ (-1)^{\langle \gamma_1 , \gamma_2 \rangle}$ and $\gamma \in \frac{1}{2} H_1(\Sigma, \Z)$. As in Proposition (\ref{KS-algebra-iso-proposition}) we have
$$  \hat{\mathfrak{g}}' =  MCH([0,1] \times \Sigma)^{ann}   .$$
%$$ \hat{e}_{\gamma}' \mapsto  \text{exp}(( w_{\gamma} , w_{\gamma}' ))  $$
%where $ w_{\gamma} , w_{\gamma}' $ are defined as before, but now $w_{\gamma}^{ann}  = w_{\gamma}$. 

\subsection{Kontsevich-Soibelman algebra associated to $(X,L)$}

$MCH(L, \mathfrak{c})$ is a module over $MCH(\Sigma \times [0,1])$. The module structure is induced by the identification 
$$( \Sigma \times [0,1] ) \sqcup_{\Sigma} M \cong M $$
obtained by gluing the boundary component $\Sigma \times \{ 0 \}$ to the boundary of $M$. 
From proposition above we obtain an action of $ \hat{\mathfrak{g}}$ on $MCH(L)$. We will keep the notation $ \hat{e}_{\gamma}$ for the operator on 
$MCH(L| \mathfrak{c})^{\diamond, \diamond}$ associated to $ \hat{e}_{\gamma} \in   \mathfrak{g} $.

The same considerations apply to obtain an action of $\hat{\mathfrak{g}}^{ann}$ on 
$MCH(L| \mathfrak{c})^{\diamond, \diamond}$.

As in \cite{OGW3}, the space of topological charges $\Gamma$ is defined as the relative homology:
$$  \Gamma = H_2(M,L, \Z)  .$$
%endowed with the canonical pairing $\langle \bullet , \bullet \rangle : \Gamma \times \Gamma \rightarrow \Z$.

%Define the Kontsevich-Soibelman algebra associated to the geometry $(M,L)$ as follows.
Denote by
$$Q: H_1(\Sigma,\Z) \rightarrow H_1(L, \Z)  $$
the map in homology induced by the inclusion $\Sigma \rightarrow L$. Consider the abelian group
\begin{equation} \label{KS-lattice}
%{\Gamma}_{KS}= \{  (\gamma, \beta) \in H_1(\Sigma, \Z) \times \Gamma  | \text{  } \partial \beta =  Q(\gamma) \} 
{\Gamma}^{\spadesuit}= \{  (\gamma, \beta) \in H_1(\Sigma, \Z) \times \Gamma  | \text{  } \partial \beta =  Q(\gamma) \}
\end{equation}
%$${\Gamma}_{KS}=  H_1(\Sigma, \Z) \times_{H_1(L, \Z)} \Gamma  $$
endowed with the skew-symmetric pairing
$$\langle \bullet , \bullet \rangle : {\Gamma}^{\spadesuit} \times  {\Gamma}^{\spadesuit} \rightarrow \Z$$ 
$$ \langle (\gamma_1, \beta_1 ), (\gamma_2,\beta_2) \rangle = \langle \gamma_1 , \gamma_2   \rangle ,$$
where in the right side we have used the usual intersection pairing of $H_1(\Sigma, \Z)$.

The quantum Kontsevich-Soibelman algebra associated to $(X,L)$ is defined by
\begin{equation} \label{KS-algebra} 
\mathfrak{g}^{\spadesuit}  = \bigoplus_{(\gamma, \beta ) \in \Gamma^{\spadesuit} }\Q[[g_s]]    \hat{e}_{\gamma, \beta}
 %  \mathfrak{g}^{\spadesuit} = \bigoplus_{(\gamma, \beta ) \in \Gamma_{KS}} \Q[q^{\pm \frac{1}{2}}]    \hat{e}_{\gamma, \beta}
\end{equation}
%$$  \mathfrak{g} = \bigoplus_{(\gamma, \beta ) \in \Gamma_{KS}} \Q[[g_s]]    \hat{e}_{\gamma, \beta}$$
equipped with the product structure defined by 
$$\hat{e}_{\gamma_1,\beta_1}  \hat{e}_{\gamma_2, \beta_2} =   q^{ \frac{\langle \gamma_1 , \gamma_2 \rangle}{2}}  \hat{e}_{\gamma_1 + \gamma_2, \beta_1 + \beta_2} .$$
%In our paper the parameter $q$ is a formal parameter is related to the string coupling $g_s$ by
%$$     q^{\frac{1}{2}} = - e^{ \frac{g_s}{2}}       .$$

\subsection{Boundary States}

Define the subgroup of $ \Gamma$   of \emph{flavor} charges as  
%$$   \Gamma_{flavor} = \text{Im} \{   H_1(\Sigma,\Z) \rightarrow H_1(L, \Z)   \}  .$$ 
%$$   \Gamma_{flavor} = \text{Im} \{   \Gamma_{\partial} \rightarrow \Gamma   \}  ,$$
$$   \Gamma_{flavor} =  \{  \beta | \partial \beta \in   \text{Im} (Q)   \}  .$$
Observe that $ \Gamma_{flavor} $ is the image of $\Gamma^{\spadesuit}$ in $\Gamma$ by the projection
$ (\gamma,\beta) \mapsto \beta  .$

Fix a positive real number $C^{supp}>0$. 
Fix a representative $w^{ann}$ of the Chern Class. 
A  \emph{boundary state} $Z    $   consists in an array $Z= (Z_{\beta})_{\beta }$ with $Z_{\beta} \in   MCH( \partial \beta | w^{ann})$ and 
%$ Z_{\beta} = 0$ if $\omega(beta) < 0$.
$$Z_{\beta} \neq 0 \Rightarrow || \beta || \leq C^{supp} \omega (\beta).$$
Denote by $\mathfrak{H}_{\mathfrak{c}}$ the vector space of boundary states with chern class $\mathfrak{c}$.

%$$  \mathcal{H} =   \prod_{  \beta \in  {\Gamma}_{flavor}}   \mathcal{H}_{\beta}  .$$
%$$  \mathcal{H} =   \{ \{ v_{\beta} \}_{  \beta \in  {\Gamma}_{flavor}} |   v_{\beta} \in   \mathcal{H}_{\beta} ,  \{ v_{\beta} \neq 0 \} \cap \Gamma^- \text { is finite } \}.$$
%$$  \mathfrak{H} =   \{ \{ v_{\beta} \}_{  \beta \in  {\Gamma}_{flavor}} |  \quad v_{\beta} \in   \mathcal{H}_{\beta} , \quad v_{\beta} = 0 \text{   for } \omega(\beta) \ll 0 \}.$$

%$$---  \mathfrak{H} =   \{ \{ Z_{\beta} \}_{  \beta \in  {\Gamma}_{flavor}} |  \quad Z_{\beta} \in   MCH( \partial \beta) , \quad Z_{\beta} = 0 \text{   for } \omega(\beta) \ll 0 \}.$$

There is an obvious action of $\mathfrak{g}^{\spadesuit}$  on $\mathfrak{H}$ defined by
$$ \hat{e}_{\gamma_0, \beta_0} (   ( Z_{\beta} )_{  \beta \in  {\Gamma}_{flavor}} ) =  ( \hat{e}_{\gamma_0}  Z_{\beta - \beta_0} )_{  \beta  \in  {\Gamma}_{flavor}} $$
%$$\mathcal{Z}_{\mathfrak{c}} \rightarrow  \mathcal{Z}_{\mathfrak{c}} .$$

Let $\mathfrak{g}^{ \bigstar}$ be  the set of formal sums 
$$ \sum_{ (\gamma, \beta ) \in \Gamma_{KS} } a_{\gamma, \beta}    \hat{e}_{\gamma, \beta} $$
with the properties
$$a_{\gamma, \beta} \neq 0 \Rightarrow || \beta || \leq C \omega (\beta),$$
%$$   \# \{   (\gamma, \beta) | \quad a_{\gamma, \beta} \neq 0 , \quad \omega(\beta) < E \}  < \infty \quad \quad  \forall E \in \R .$$
%for each positive real number $E $.
$$   \# \{   \gamma | \quad a_{\gamma, \beta} \neq 0  \}  < \infty \quad \quad  \forall \beta \in \Gamma  .$$
$\mathfrak{g}^{ \bigstar}$ has an obvius algebra structure. $\mathfrak{H}$ is a modulo on $\mathfrak{g}^{ \bigstar}$.

Let $\mathfrak{J}$ the left ideal of $\mathfrak{g}^{ \bigstar}$ generated by $ \{   \hat{e}_{\gamma, 0 } -1 \}_{\gamma \in \text{Ker}(Q) } $. We have an isomorphism of moduli 
$$ \mathfrak{H} =  \mathfrak{g}^{ \bigstar}/\mathfrak{J} .  $$

%The algebra structure of $   \mathfrak{g}_{KS} $ induces an algebra structure on $\hat{  \mathfrak{g}}_{KS}$. With this algebra structure  $\hat{  \mathfrak{g}}_{KS}$
%acts on $\mathfrak{H}$, extending the action of $   \mathfrak{g}_{KS} $.
%$$ \mathfrak{H} =    \mathfrak{g}_{KS}/ phase .  $$

%$$ \Gamma_{KS}=  \{ \{ a_{\gamma, \beta} \}_{  \beta \in  {\Gamma}_{flavor}} |  \quad v_{\beta} \in   MCH( \partial \beta) , \quad v_{\beta} = 0 \text{   for } \omega(\beta) \ll 0 \}$$

We can also consider the action of $ MCH([0,1] \times \Sigma)^{ann}$:
$$\hat{e}^{ann0}_{\gamma}   (   ( Z_{\beta} )_{  \beta \in  {\Gamma}_{flavor}} ) =   ( \hat{e}_{\gamma}^{ann0}  Z_{\beta} )_{  \beta  \in  {\Gamma}_{flavor}}$$
$$\mathcal{Z}_{\mathfrak{c}} \rightarrow  \mathcal{Z}_{\mathfrak{c}+ 2 Q(\gamma)} .$$

\subsubsection{Novikov Ring}

Denote by $ \mathfrak{g}(\Lambda_0)$  the Novikov Ring with coeffients on $ \mathfrak{g}$. Hence an element of    $ \mathfrak{g}(\Lambda_0)$   is defined by a formal sum $\sum_i a_i T^{\lambda_i}  $ with $a_i \in \mathfrak{g}$, $\lim_{i \rightarrow \infty} \lambda_i = + \infty$, $\lambda_i \geq 0$.
%We consider the algebra $ \mathfrak{g}((T))$ with coeffiecients on the Novikov ring.

$ \mathfrak{g}(\Lambda_0)$ is equipped with an obvious structure of  algebra.
We have an homomorphism of algebras
\begin{equation} \label{KS-algebra-nov-map}
     \hat{\mathfrak{g}}^{\bigstar}  \rightarrow   \mathfrak{g}(\Lambda_0)
\end{equation}
defined by
$$     \hat{e}_{(\gamma, \beta)} \mapsto     T^{ \omega(\beta)}  \hat{e}_{\gamma}.$$

We  also define $ \mathfrak{H}^{nov}$ whose elements are formal power series 
$\sum_i Z_i T^{\lambda_i}  $ with $Z_i \in MCH(L)$, $\lambda_i \geq 0$, $\lambda_i \rightarrow \infty$.
$ \mathfrak{H}^{nov}$ is a modulo on $ \mathfrak{g}(\Lambda_0)$.

We have an homomorphism  
$$  \mathfrak{H} \rightarrow   \mathfrak{H}^{nov}$$
defined by
$$  ( Z_{\beta} )_{  \beta \in  {\Gamma}_{flavor}} \mapsto \sum_{\beta} T^{ \omega(\beta)}  Z_{\beta}.$$
This map is compatible with the moduli structure $  \mathfrak{H} $  and $  \mathfrak{H}^{nov}$ and the homomorphism of algebras (\ref{KS-algebra-nov-map}).

\subsection{KS-algebra including the degrees} \label{KS-algebra including components section}

We now extended the objects considered above to the  context of section \ref{nice-MC2-section}.
%We need to extend the operators $ \hat{e}_{\gamma, \beta} $ in the context of section \ref{nice-MC2-section}.
%(\ref{nice-cycles2}).
In this context, coherent $MC$-cycles are defined  in the same way.

\subsubsection{The case $\Sigma \times \R$}

$MCH(\R \times \Sigma )^{\diamond , \diamond }$ is an algebra for the same reason above. 
However, it is  different from the one considered above. 

The map (\ref{KS-alebra-coh-map}) is not anymore an homomorphism of algebras. However we can define an homomorphism of Lie Algebras as follows.

%Denote by $ \hat{e}_{\gamma, \beta} $ the same element of (\ref{KS-alebra-coh-map})

%$ \hat{e}_{\gamma, \beta} $ acts on the generators (\ref{generator2}) adding a new vertex $v_0$ with data , $\beta_{v_0} = \beta$, $w_{v_0,h} = w \times \{ t_h \}$, where $w$ is any one dimensional cycle representing the homology class $\gamma$, and $(t_h)_h$ are real numbers with $t_h \in (T, \infty)$ and $t_h \neq t_{h'}$  if $h \neq h' $.   

%Relation (\ref{KS-algebra-product}) does not hold anymore on the context of (\ref{nice-cycles2}), that is, the operators $ (\hat{e}_{(\gamma, \beta)})_{(\beta,\gamma)} $  do not define an algebra. 
%However we have the following

%\begin{lemma}
Set 
\begin{equation} \label{KS-lie}
    \hat{e}_{\gamma, \beta}^{lie} = \frac{\hat{e}_{\gamma, \beta}}{q^{\frac{1}{2}}- q^{- \frac{1}{2}}} .
\end{equation}

%In analogy with  \cite{KS}, 
Set 
$$ \mathfrak{g}^{qu}= \bigoplus_{(\gamma, \beta) \in \Gamma^{\spadesuit} } \Q[[g_s]]   \hat{e}_{\gamma, \beta }^{lie} $$
equipped with the Lie bracket 
%\begin{equation} \label{KS-quantum-braket}
% [  \hat{e}_{\gamma_1}^{lie}  ,   \hat{e}_{\gamma_2}^{lie} ]  = (-1)^{\langle \gamma_1 , \gamma_2 \rangle}    \frac{q^{ \langle \gamma_1 , \gamma_2 \rangle  } - q^{ - \langle \gamma_1 , \gamma_2 \rangle  } }{q^{\frac{1}{2}}- q^{- \frac{1}{2}}}   \hat{e}_{\gamma_1 + \gamma_2}^{lie}
% \end{equation}
\begin{equation} \label{KS-quantum-braket}
 [  \hat{e}_{\gamma_1, \beta_1}^{lie}  ,   \hat{e}_{\gamma_2, \beta_2}^{lie} ]  = (-1)^{\langle \gamma_1 , \gamma_2 \rangle}    \frac{q^{ \langle \gamma_1 , \gamma_2 \rangle  } - q^{ - \langle \gamma_1 , \gamma_2 \rangle  } }{q^{\frac{1}{2}}- q^{- \frac{1}{2}}}   \hat{e}_{\gamma_1 + \gamma_2, \beta_1+ \beta_2}^{lie}
 \end{equation}

%\end{lemma}

\begin{lemma} \label{lie-algebra-homorphism-lemma}
    The map (\ref{KS-alebra-coh-map}) defines an homomorphism of Lie Algebras. 
\end{lemma}
\begin{proof}
The reason why in the context of section \ref{nice-MC2-section} we do not get the RHS of  (\ref{KS-algebra-product}) is that in the isotopy  used in the  proof of (\ref{KS-algebra-iso-proposition}) we generate   $MC$-cycles  which have two vertices instead of one. 

We claim that  these elements  are the same when  the order of the product on the LHS is switched and hence they cancel on the bracket. 

The claim follows from the fact that in the isotopy  we never cross two cycles associated to the same vertex and hence the elements generated by the RHS of  relation (\ref{jump2}) have always one vertex. 
%Hence the elements with two vertices are the one   not generated by crossing. 
%It follows that the elements generated by the isotopy with two vertices are the same when  the order of the product on the LHS is switched. Hence they cancel on the bracket. 
\end{proof}

\begin{proposition}
We have an isomorphism of algebras 
$$MCH(\R \times \Sigma )^{\diamond , \diamond } = U(\mathfrak{g}^{qu})  .$$
\end{proposition}
\begin{proof} 
It is immediate from the definitions and the Lemma (\ref{lie-algebra-homorphism-lemma}).
\end{proof}

We define the vector space $\mathfrak{H}^{\flat}$ in a way analogous to $\mathfrak{H}$ using  $MC$-cycles with components of section \ref{nice-MC2-section} instead of  \ref{generator-nice}.  
We obtain an action of $\mathfrak{g}^{\spadesuit, qu}$ on $\mathfrak{H}^{\flat}$.

%\subsubsection{-KS-algebra associated to $(X,L)$}

%Now we can proceed as before to define $\mathfrak{g}^{\spadesuit, qu}$ 

%Set 
%$$ \mathfrak{g}^{\spadesuit, qu}= \bigoplus_{(\gamma, \beta) \in \Gamma^{\spadesuit} } \Q[[g_s]]   \hat{e}_{\gamma, \beta }^{lie} $$
%equipped with the Lie bracket 
%\begin{equation} \label{KS-quantum-braket}
% [  \hat{e}_{\gamma_1, \beta_1}^{lie}  ,   \hat{e}_{\gamma_2, \beta_2}^{lie} ]  = (-1)^{\langle \gamma_1 , \gamma_2 \rangle}    \frac{q^{ \langle \gamma_1 , \gamma_2 \rangle  } - q^{ - \langle \gamma_1 , \gamma_2 \rangle  } }{q^{\frac{1}{2}}- q^{- \frac{1}{2}}}   \hat{e}_{\gamma_1 + \gamma_2, \beta_1+ \beta_2}^{lie}
% \end{equation}

\section{Compactification Frames}

So far we have considered the classes $\beta \in \Gamma$ independently. To define the wavefunction as in \cite{CCV2} we need to fix frame of $L$.

\begin{definition}
%A frame of $L $ is a maximal isotropic sublattice $F$ of $H_1(\Sigma, \Z)$. 
A frame of $L$ is a Lagrangian subvector space $F$ of $H_1(\Sigma, \R)$, such that $F \cap H_1(\Sigma, \Z)$ generate $F$ as vector space.
\end{definition}

\begin{remark}
The lattice $F$ is not necessarily transverse to $K$, that is we do \emph{not} assume that $F \rightarrow H_1(L, \Z)$ is surjective. This is in accord with the definition in three-dimensional gauge theory (\cite{CCV2}). 
\end{remark}

To a frame $F$ it is associated a compactification $L_F$ of $L$ as follows. Let $H_F$ be the handlebody with $\partial H_F = \Sigma$ and $ F= \text{Ker} \{  H_1(\Sigma) \rightarrow H_1(H_F) \} $, i.e., $H_F$ is obtained from $\Sigma$ filling the cycles of $F$. Set  
$$ L_F = L \sqcup_{\Sigma} H_F  $$
where we identify the boundary $\partial L = \partial H_F = \Sigma$.
Roughly , $L_F$ is the compactification of $L$ for which $F$ contractible at infinity.

\begin{remark}
In the physical lecturature, to define three dimensional gauge theories for not compact geometry it is necessary to fix  boundary conditions at infinity (see for example \cite{CCV}, \cite{CCV2}). These boundary conditions are defined in geometric terms of a compactification  of $L$ defining the set of cycles of $L$ contractible at infinity (see \cite{CCV}). The last data is specified by a choice of a frame $F$. This procedure  can be seen mirror to our definitions above, according to the topological string/3d gauge theory correspondence.
\end{remark}
% reducing the problem of define rational invariants to the compact context of \cite{OGW3}. This is the usual approach used  in  physics in the  context of three dimensional gauge theories (see for example \cite{CCV}, \cite{CCV2}). 

From the Mayer-Vietoris sequence 
$$... \rightarrow  H_1(\Sigma) \rightarrow H_1(L) \oplus H_1(H_F) \rightarrow H_1(L_F) \rightarrow H_0(\Sigma) \rightarrow ... $$
we obtain the isomorphism
\begin{equation} \label{homology}
 H_1({L}_F, \Z) =  H_1(L, \Z) / Q(F) .
\end{equation}

Fix a $F \subset H_1(M . \Z)$  transversal to $\text{Ker}(Q)$.
For each $\gamma \in H_1(L, \Z)$ there exists a unique $F_{\gamma} \in F$ such that $Q(F_{\gamma})= \gamma$,
%For each $\beta \in \Gamma_{flavor}$ there exists a unique $F_{\beta} \in F$ such that $(\beta, F_{\beta}  ) \in  \Gamma_{KS}$.
%This defines an isomorphism 
%$$ \Gamma_{KS} \cong  \Gamma_{flavor} \oplus \text{Ker}(Q) .$$
To simplify notation denote
$$  \hat{e}_{\beta}=    \hat{e}_{(\beta,  F_{\partial \beta} )} \text{        } \forall \beta \in \Gamma_{flavor} ,$$
$$  \hat{e}_{\gamma}=    \hat{e}_{(0,  \gamma )} \text{        } \forall \gamma \in \text{Ker}(Q)  .$$
We stress that $\hat{e}_{\beta}$ depends on the choice of the frame $F$.

%For each $\beta \in \Gamma_{flavor}$ set
%$$   \hat{e}_{\beta}=    \hat{e}_{(\beta,  \gamma_{\beta} )}   $$
%where $\gamma_{\beta}$ is the only element of $F$ with $\partial \beta =Q( \gamma_{\beta})$.

% there exist a unique $f_{\beta} \in F$ such that $(\beta, f_{\beta}  ) \in  \Gamma_{KS}$.
%This defines an isomorphism 
%$$ \Gamma_{KS} \cong  \Gamma_{flavor} \oplus \text{Ker}(Q) .$$
%In order to simplify the notation define
%$$  \hat{e}_{\beta}=    \hat{e}_{(\beta,  f_{\beta} )} \text{        } \forall \beta \in \Gamma_{flavor} ,$$
%$$  \hat{e}_{\gamma}=    \hat{e}_{(0,  \gamma )} \text{        } \forall \gamma \in \text{Ker}(Q)  .$$
%We stress that $\hat{e}_{\beta}$ depends on the choice of the frame $F$.

The  factorization property \cite{QME} implies that there exists an array   $( a(\beta) )_{\beta \in \Gamma_{flavor}}$ (depending on $F$)  with $a(\beta) \in \frac{1}{g_s} \Q[[g_s]]$  such that
\begin{equation} \label{factorization-state}
      Z =   a_0 \hat{e}^{ann0}_{F_{\mathfrak{c}}}  \prod_{\beta \in \Gamma_{flavor}} \text{exp} (a(\beta)  \hat{e}_{\beta} ) Z_0 .
\end{equation}

\subsection{Integrality}

We say that $\{ a(\beta) \}_{\beta \in \Gamma_{flavor}}$ of (\ref{factorization-state}) satisfies the Ouguri-Vafa integrality if there exists integer numbers  $  N_{\beta, j} \in \Z   $, $\beta \in H_2(X,L)$, $j \in \frac{1}{2}\Z$ with 
$$ |\{  N_{\beta, j} \neq 0, |\omega(\beta)| <E  \}  |   < \infty \text{    for   each  } E>0 $$
and
\begin{equation} \label{Ooguri-Vafa}
  a(\beta) = \sum_{n,j} N_{\frac{\beta}{n},j} \frac{q^{nj}}{n(q^{\frac{n}{2}} -  q^{-\frac{n}{2}} )}  .
\end{equation}
In the last  sum $n$ run over the integers numbers such that $\frac{\beta}{n} \in \Gamma_{flavor}$. 
%Observe that the integrality (\ref{Ooguri-Vafa})  is compatible with the change of frame.

Observe that if we twist the frame by an element $\gamma_0 \in H_1(\Sigma, \Z)$ the integers changes according to the formula
$$ N_{\beta, j} \mapsto   N_{\beta, j+ \langle \gamma_0 , \partial \beta \rangle}.$$
In particular this means that if  the integrality property holds in a frame holds in each frame.

Using (\ref{Ooguri-Vafa}) we can rewrite (\ref{factorization-state}) as 
$$  Z = a_0 \hat{e}^{ann0}_{F_{\mathfrak{c}}} \prod_{\beta \in \Gamma, j}   (  S_{\beta, j} )^{   N_{\beta, j}   }  Z_0 .$$
where
$$  S_{\beta, j} = \text{exp} (   \sum_{n}  \frac{q^{n(j+ \frac{1}{2})}}{n(q^n -  1 )} \hat{e}_{ n \beta}   ) .$$

$ S_{\beta, j}$ are generators of Cluster trasformations: 
\begin{equation} \label{cluster-trasformation}
S_{\beta, j}\hat{e}_{ \gamma} S_{\beta,j}^{-1} =
\begin{cases}
   \hat{e}_{ \gamma} \left(  \prod_{k=1}^{ |\langle \gamma , F_{\partial \beta} \rangle| } (1+ q^{j- \frac{1}{2} +k  }    \hat{e}_{ \beta} )  \right)  \text{      if    }   \langle F_{ \partial \beta} , \gamma \rangle >0 :\\
 \hat{e}_{ \gamma}    \text{      if    }   \langle F_{ \partial \beta} , \gamma \rangle =0; \\
   \hat{e}_{ \gamma}  \left(  \prod_{k=1}^{ |\langle \gamma , F_{\partial \beta} \rangle| } (1 + q^{j+\frac{1}{2} - k  }    \hat{e}_{ \gamma} )  \right)^{-1}  \text{      if    }     \langle F_{ \partial \beta} , \gamma \rangle <0.
\end{cases}
\end{equation}

%\begin{lemma} \label{polinomio} 
If the Ooguri-Vafa integrality (\ref{Ooguri-Vafa})  holds the vector $Z$ satisfy the quantum equations:
\begin{equation} \label{quantum-equation}
  \mathcal{A}_{\alpha}  Z =0  
  \end{equation}
% \text{   for  }  0 \leq i \leq g  
with
$$  \mathcal{A}_{\alpha}  \in  \prod_{(\gamma, \beta)} \Z[q^{\pm \frac{1}{2}}]   \hat{e}_{\gamma, \beta}  \text{   for  }  0 \leq \alpha \leq g  .$$
%$$  \mathcal{A}_{\alpha} (\hat{e}_{\gamma}) \Psi =0  \text{    for    }  0 \leq    \alpha \leq \frac{1}{2} \text{dim}( H_1(\Sigma, \C) )$$%\end{lemma}
%\begin{proof}

Equation (\ref{quantum-equation}) is obtained applying the sequence of cluster transformations (\ref{cluster-trasformation})  to the equations 
$$  ( \hat{e}_{\gamma} -1 ) Z_0 = 0 \text{   for each   } \gamma \in \text{Ker}(Q) .$$
The last equations should be considered as the quantization of the classical equations
\begin{equation} \label{classical}
x_{\gamma}=0 \text{ for each } \gamma \in \text{Ker}(Q)
\end{equation}
where for each $\gamma \in H_1(\Sigma,  \Z)$, $x_{\gamma}: H^1(\Sigma,  \C)\rightarrow \C $ denotes the linear map defined by $\gamma$.

Equations (\ref{classical}) define the sub-vector space 
of $H^1(\Sigma,  \C)$ 
\begin{equation} \label{vector-lagrangian}
  \frac{H^1(L, \C)}{H^1_{comp}(L, \C)}   <  H^1(\Sigma,  \C).  
\end{equation}
%The quantization due to holomorphic curves   corrects the solution of  $\Psi_0$ to the vector $\Psi$ defined in this paper.
It is a standard result, and easy to check, that (\ref{vector-lagrangian}) is a Lagrangian sub-vector space of $H^1(\Sigma,  \C)$ equipped with the symplectic form defined by the intersection pairing.

\begin{remark} \label{A-branes-remark}
An $A$-brane is defined by a lagrangian sub-manifold $L$ and an abelian flat connection of $L$. The moduli space of $A$-branes is locally modelled on $H^1(L, \C)$, where the real part can be viewed coming from the moduli of the deformations of the Lagrangian and the imaginary part defines a flat $U(1)$ connection on $L$.    

The vector space (\ref{vector-lagrangian}) is the tangent space of the classical moduli space of $A$-branes up to compact deformations.  
Hence the vector $Z$  can be considered as the quantum  corrected moduli obtained including  the contribution of holomorphic curves to the classical moduli of $A$-branes  . It is    
a subspace of  the deformation quantization of the complex symplectic vector space   $( H^1(\Sigma,  \C), \langle \bullet , \bullet \rangle )$, which is the quantum torus with algebra of functions defined by:
$$ \hat{e}_{\gamma_1}  \hat{e}_{\gamma_2} =   q^{ \langle \gamma_1 , \gamma_2 \rangle}  \hat{e}_{\gamma_2}  \hat{e}_{\gamma_1} \text{   for     } \gamma_1, \gamma_2 \in  H_1(\Sigma,  \Z). $$ 
\end{remark}

%A standard computation shows that  
%$$  \text{Ad} (\prod_{n} \text{exp} (a(n \gamma)  \hat{e}_{n \gamma} )) (\hat{e}_{\gamma}) = ()   $$
%$$  S_{\gamma_h} X_{\gamma} S_{\gamma_h}^{-1} = X_{\gamma } (\prod_{m= M_{\gamma_h}}^{M_{\gamma_h}} \Phi_{\langle \gamma, \gamma_h  \rangle} ((-1)^m q^m X_{\gamma_h})^{\epsilon N_{m, \gamma_h}} ) . $$

\subsection{Multi-link Homomorphism}

%For not compact $L$, the definition of Open Gromov-Witten invariants as rational numbers involves the choice of a frame of $L$ ( see \cite{frame} ). 
%Denote by $Q$ the obvious map
%$$ Q:  H_1(\Sigma,\Z) \rightarrow H_1(L,\Z) .$$
%Define
%$$K=\text{ker} \{ H_1(\Sigma, \Z) \rightarrow H_1(L,\Z) \} .$$ 
%$$ K=\text{ker}(Q)  .   $$

The inclusion $L \hookrightarrow {L}_F$ induces a link map 
$$  \text{Link}_F :  \sqcup_{\gamma \in \text{Im}(Q)}  \mathcal{P}_{\gamma} \rightarrow \Z .$$
It has the characterizing property:
$$ \text{Link}_F (w, w') =0  $$
if $(w,w') \in  \mathcal{P}_{\gamma}$ for $\gamma \in F$,  $w$ has support on $\{t \} \times \Sigma $, $w'$ has support on $\{t' \} \times \Sigma $ with $t \neq t' $.

%The inclusion $L \hookrightarrow {L}_F$ as open subset
%induces a map 
%$MCH_0(L) \rightarrow MCH_0({L}_F) .$
%between the multi curve homology of $L$ and the multi curve homology of ${L}_F$. 
%In this way we are reconducted to the compact case, thus we can get rational invariants using the multi-link homomorphism defined in \cite{OGW3}.

We can now extend the Multi-Link homomorphism
$$  \text{Multi-Link}_F : \mathcal{Z}_{\gamma} \rightarrow \Q[[g_s]] .$$
defined by
$$ (H ,  \{ w_h \}_h ) \mapsto \sum_{E} \frac{g_s^{E^{in}}}{Aut(E)} \prod_{e \in E^{in}} \text{Link}_F (w_e,w_e')  \prod_{e \in E^{ex}}  \langle x , w_e   \rangle  $$
%\text{exp}(  \langle x , w_e   \rangle)   $$
where $ (w_e,w_e')  = (w_h,w_{h'})$ for $e = \{ h,h'\} \in E^{in}$, $w_e = w_h$ for $e = \{ h \} \in E^{ex}$.

\subsection{Wave Function}

Introduce a set of formal variables $(x_{\beta})_{ \beta  \in \Gamma_{flavor} } $ with relation $x_{\beta_1} \cdot x_{\beta_2} = x_{\beta_1+ \beta_2}$.

Set
$$  \Psi_{\beta} = \text{Link}_F (Z_{\beta})$$
$$ \Psi(x) =  \sum_{ \beta  \in \Gamma_{flavor} } \Psi_{\beta}  x_{\beta} $$

The action of $   \mathfrak{g}^{\bigstar} $  on $\mathfrak{H}$ induces an action on $\Psi $ which can be obtain from the action  on the formal variables given by
\begin{equation} \label{action-x}
 \hat{e}_{\gamma_0, \beta_0}: x_{\beta} \mapsto
\begin{cases}
 x_{\beta + \beta_0}  \text{    if   } \gamma_0 \in F \\
 \text{exp}(g_s  \langle \gamma_0 , \partial \beta   \rangle)  x_{\beta }  \text{    if   } \gamma_0 \in \text{Ker}(Q)
 \end{cases}
 \end{equation}

\subsubsection{Novikov Ring}

Assuming $\mathfrak{c} \in 2 H_1(M,\Z)$ we make the following substation: 
$$x_{\beta } \mapsto  T^{ \omega(\beta)} \text{exp}( \langle(x, \partial \beta + \frac{1}{2}\mathfrak{c} \rangle ) $$
where $x$ is a formal variable with value on $H^1(L, \R)/ 2 \pi H^1(L, \Z)$.

Hence the wave function with coefficients on the Novikov Rings is given by
$$ \Psi(T, x) =  \text{exp}( \langle x,  \frac{1}{2}\mathfrak{c} \rangle )  \sum_{ \beta  \in \Gamma_{flavor} } \Psi_{\beta} T^{ \omega(\beta)}  \text{exp}( \langle x, \partial \beta \rangle ) .$$

The action of $   \mathfrak{g}(\Lambda_0) $  on $\mathfrak{H}^{nov}$ induces an action on $\Psi $ which can be obtain from the action  on the formal variables given by
%The action (\ref{action-x}) becomes
\begin{equation} \label{action-z}
 \hat{e}_{\gamma_0}: \langle x, \gamma \rangle  \mapsto
\begin{cases}
\langle x, \gamma + \gamma_0 \rangle \text{    if   } \gamma_0 \in F \\
 \langle x + g_s PD( \gamma_0 ), \gamma \rangle  \text{    if   } \gamma_0 \in \text{Ker}(Q)
\end{cases}
 \end{equation}

%\begin{equation} \label{action-z}
%- \hat{e}_{\gamma_0}: z^{\gamma} \mapsto
%\begin{cases}
%z^{\gamma + \gamma_0}  \text{    if   } \gamma_0 \in F \\
% \text{exp}(g_s  \langle PD(\gamma_0 ) , \gamma  \rangle)  z^{\gamma }  \text{    if   } \gamma_0 \in \text{Ker}(Q)
%\end{cases}
% \end{equation}
 
 $\hat{e}_{\gamma_0}^{ann}$ acts on the pair $(x ,\mathfrak{c}) $ as 
 \begin{equation} \label{action-z-ann}
 (x ,\mathfrak{c}) \mapsto
\begin{cases}
  (x ,   \mathfrak{c} + 2 \gamma_0)   \text{    if   } \gamma_0 \in F \\
  ( x + g_s \gamma_0 ,  \mathfrak{c}  ) \text{    if   } \gamma_0 \in \text{Ker}(Q)
\end{cases}
 \end{equation}

\section{Semiclassical limit}

We now define a representation of the semiclassical Kontsevich-Soibelman Lie algebra $\mathfrak{g}^{scl}$  in terms of multi-disk-homology introduced in \cite{KSWCF} and we show that this can be considered in the natural way as  the semiclassical limit of the representation defined in Section \ref{KS-algebra including components section}.

\subsection{Multi-Disk Homology}

%We shall define the  (nice) multi-disk homology as a formal $\Q$-vector space $\mathcal{H}^0$  modulo a suitable equivalence relation.

We now define the vector space of nice multi disk cycles $\mathcal{Z}^{disk}$.

Consider objects as 
\begin{equation} \label{nice-chains0}
 (V, ( w_v )_v , ( \beta_v )_v)  
\end{equation}
where
\begin{itemize}
\item A finite set $V$, called the set of vertices. 
\item For each $v \in V$ an element $\beta_v \in H_2(X,L)$ called the homology class of $v$, 
\item for each $v \in V$ a curve of real dimension one $\gamma_v$, with $[\gamma_{v}] = [\partial \beta_v] \in H_1(L, \Q)$. 
\end{itemize}
We assume the support property $|| \beta_v || \leq C^{supp} \omega(\beta_v)$ for each $v \in V$.

Let $\mathfrak{Gen^0}(\beta)$ be the set of objects  (\ref{nice-chains0}) modulo the obvious equivalence, with $\sum_v \beta_v = \beta$.

$\mathcal{Z}^{disk}$ is defined as the formal vector space generated by $\mathfrak{Gen^0}(\beta)$ .

%A Multi-Disk-cycle $W \in MDH$ is a formal sum of the objects (\ref{nice-chains0}).

To define isotopies of Multi-Disk cycles we consider one parameter family of the objects (\ref{nice-chains0}):
\begin{equation} \label{nice-chains1-scl}
   ( V, (\tilde{w}_v)_v  ,(\beta_v)_v, [a,b])   
   %\}_{v \in V,  t \in [a,b]}.  
\end{equation}
where $\tilde{w}_v$ are one parameter family of one cycles parametrized by $[a,b]$.

An isotopy $\tilde{Z}^{disk}$ of $MD$-cycles is a formal sum of elements (\ref{nice-chains1-scl}) with the following proprieties for a pair $v_1 \neq v_2 \in V$, $Z^{disk,t}$ jumps according to the formula 
\begin{multline} \label{MD-jump} 
 Z^{disk,t_0^+} - Z^{disk, t_0^-} =  \pm (  (V \setminus \{  v_1 ,v_2 \} ,  (\tilde{w}_v (t_0))_v ,  (\beta_v)_v   )        \sqcup   (v_0,  \tilde{w}_{v_1}(t_0) +\tilde{w}_{v_2}(t_0), \beta_{v_1} + \beta_{v_2}  )  ) 
\end{multline}
% \{  ( \tilde{\gamma}_v ,\beta_v)     (1)      \}_{v \in V}     - \{   ( \tilde{\gamma}_v , \beta_v)    (0)      \}_{v \in V}  =  \\
%  =  \pm ( \{  (  \tilde{\gamma}_v (t),\beta_v)        \}_{v \in V \setminus \{  v_1 ,v_2 \} } \cup   \{   (\tilde{\gamma}_{v_1}(t) +\tilde{\gamma}_{v_2}(t), \beta_{v_1} + \beta_{v_2})   \} ) 
where the sign is provided by the sign of the intersection $\tilde{w}_{v_1}(t_0) \cap \tilde{w}_{v_2}(t_0)$.

\subsubsection{Semiclassical Konsevich-Soibelman Algebra}

To the pair $(X,L)$  we associate the semiclassical Kontsevich-Soibelman Lie algebra 
$$\mathfrak{g}^{scl}= \bigoplus_{(\gamma, \beta) \in \Gamma^{\spadesuit}} \Q e_{(\gamma,\beta)}$$ 
with Lie bracket
$$ [ e_{(\gamma_1, \beta_1)} , e_{(\gamma_2,\beta_2)} ] = (-1)^{\langle \gamma_1 , \gamma_2 \rangle} \langle \gamma_1 , \gamma_2 \rangle e_{ \gamma_1+ \gamma_2,  \beta_1 + \beta_2} . $$
%We now define a representation of the Lie algebra $\mathfrak{g}^{scl}$  in terms of multi-disk-homology introduced in \cite{KSWCF} and we show that this can be considered in the natural way as  the semiclassical limit of the representation defined in Proposition \ref{quantum-representation}.

We fix a spin structure on $\Sigma$ and denote by $\sigma$ the associated  quadratic differential.

For a  $(\gamma,\beta) \in \Gamma^{\spadesuit}$ consider the linear map 
%$$   e_{\beta} : \mathcal{H}_k = \mathcal{H}_{k+1}  ,$$
%\begin{equation} \label{action0}
$$ e_{(\gamma, \beta)} : \mathcal{Z}^{disk} \rightarrow \mathcal{Z}^{disk}  $$
%\end{equation}
defined on the generators (\ref{nice-chains1-scl}) as
\begin{equation} \label{action0}
 e_{(\gamma, \beta)} :   (V, ( w_v )_v , ( \beta_v )_v )    \mapsto  \sigma(\gamma)  (V \sqcup \{ v_0 \}, ( w_v )_{v \in V \sqcup \{ v_0 \}}, ( \beta_v )_{v \in V \sqcup \{ v_0 \}}  ) .
% e_{(\gamma, \beta)} :   (V,  (\gamma_v)_v , (\beta_v )_v          ) \mapsto    \{    ( \tilde{\gamma}_v , \beta_v )            \}_{v \in V} \sqcup \{  (w  , \beta) \}
\end{equation}
where  $\beta_{v_0}= \beta$ and $w_{v_0}$ is a representative of $\gamma$  with support in $(-\infty, -T) \times \Sigma $ for $T$ big enough such that each  $w_v$  have support in the complementary of $(-\infty, -T) \times \Sigma $.

%The operators ( \ref{action0})  define a representation of $\mathfrak{g}_{cl}$ with the wrong sign. To correct the sign we consider as before the dependence on the spin structure and define
%$$    e_{(\gamma, \beta)} = \sigma(\gamma)    e_{(\gamma, \beta)}'   .    $$ 

\subsubsection{The case $\Sigma \times \R$}

We now show that the Multi-Disk homology of $\Sigma \times \R$ is isomorphic to the envelopping algebra of  $\mathfrak{g}^{scl}$.

%We consider the map
%$$ e_{(\gamma_1)} \otimes e_{(\gamma_2)} \otimes ... \otimes e_{(\gamma_k)} \mapsto \{ ( \gamma_i \times {t_i}, \beta_i  \}_i $$
%for some $t_1 > t_2 > ... > t_k$. 

To an element $e_{(\gamma_1, \beta_1)} \otimes e_{(\gamma_2, \beta_2)} \otimes ... \otimes e_{(\gamma_k, \beta_k)} $  we can associated an element  $(V, (w_v)_v , (\beta_v)_v )  \in \mathcal{Z}^{disk}/isotopies$
where $V= \{ 1,2,...,k  \}$ and 
$w_i$ are one dimensional cycles representing the class $\gamma_i$, supported on $\Sigma \times \{ t_i \}$,  for  $t_1 > t_2 > ... > t_k$.

%The map
%$$ e_{(\gamma_1, \beta_1)} \otimes e_{(\gamma_2, \beta_2)} \otimes ... \otimes e_{(\gamma_k, \beta_k)} \mapsto  ( \{1,2,...,k \}, \{ w_i \}_{1 \leq i \leq k}  ) $$
%where $w_i$ are one dimensional chains representing the class $\gamma_i$, supported on $\Sigma \times \{ t_i \}$,  for  $t_1 > t_2 > ... > t_k$. The element on the RHS does not depends on the  choice of $w_i, t_i$ up to isotopy.

It is easy to check that 
$$  e_{(\gamma_1, \beta_1)} \otimes e_{(\gamma_2, \beta_2)} - e_{(\gamma_2,\beta_2)} \otimes e_{(\gamma_1, \beta_1)} - e_{(\gamma_1 + \gamma_2, \beta_1+ \beta_2)} \mapsto 0 .$$

Hence  the homomorphism
$$ e_{(\gamma_1, \beta_1)} \otimes e_{(\gamma_2, \beta_2)} \otimes ... \otimes e_{(\gamma_k, \beta_k)} \mapsto (\prod_i \sigma(\gamma_i) ) (V, (w_v)_v , (\beta_v)_v ) $$
defines  map from the enveloping algebra
$$ U(\mathfrak{g}^{scl}) \rightarrow \mathcal{Z}^{disk}/\sim . $$
It is easy to check that this map is an isomorphism.

\subsubsection{Semiclassical states }

%The semiclassical Kontsevich-Soibelman algebra associated to $(X,L)$ is defined by
%\begin{equation} \label{KS-algebra} 
%\mathfrak{g}^{\spadesuit}  = \bigoplus_{(\gamma, \beta ) \in \Gamma_{KS}} \Q[[g_s]]    e_{\gamma, \beta}
% \end{equation}
%equipped with the Lie bracket 
%$$ [ e_{(\gamma_1, \beta_1)} , e_{(\gamma_2,\beta_2)} ] = (-1)^{\langle \gamma_1 , \gamma_2 \rangle} \langle \gamma_1 , \gamma_2 \rangle e_{ \gamma_1+ \gamma_2,  \beta_1 + \beta_2} . $$

Fix a positive real number $C^{supp}>0$.  
A  \emph{semiclassical boundary state} $Z    $   consists in an array $Z= (Z_{\beta}^{disk})_{  \beta \in  {\Gamma}_{flavor}}$ with $Z_{\beta}^{disk} \in   \mathcal{Z}^{disk}_{\beta}$ and 
%$ Z_{\beta} = 0$ if $\omega(beta) < 0$.
$$Z_{\beta}^{disk} \neq 0 \Rightarrow || \beta || \leq C^{supp} \omega (\beta).$$
Denote by $\mathfrak{H}^{scl}$ the vector space of semiclassical boundary states.

%$$  \mathcal{H} =   \prod_{  \beta \in  {\Gamma}_{flavor}}   \mathcal{H}_{\beta}  .$$
%$$  \mathcal{H} =   \{ \{ v_{\beta} \}_{  \beta \in  {\Gamma}_{flavor}} |   v_{\beta} \in   \mathcal{H}_{\beta} ,  \{ v_{\beta} \neq 0 \} \cap \Gamma^- \text { is finite } \}.$$
%$$  \mathfrak{H} =   \{ \{ v_{\beta} \}_{  \beta \in  {\Gamma}_{flavor}} |  \quad v_{\beta} \in   \mathcal{H}_{\beta} , \quad v_{\beta} = 0 \text{   for } \omega(\beta) \ll 0 \}.$$

%$$---  \mathfrak{H} =   \{ \{ Z_{\beta} \}_{  \beta \in  {\Gamma}_{flavor}} |  \quad Z_{\beta} \in   MCH( \partial \beta) , \quad Z_{\beta} = 0 \text{   for } \omega(\beta) \ll 0 \}.$$

There is an obvious action of $\mathfrak{g}^{\spadesuit,scl}$  on $\mathfrak{H}^{scl}$ defined by
$$ \hat{e}_{\gamma_0, \beta_0} (   ( Z_{\beta}^{disk} )_{  \beta \in  {\Gamma}_{flavor}} ) =  ( \hat{e}_{\gamma_0,  \beta_0}  Z_{\beta - \beta_0}^{disk} )_{  \beta  \in  {\Gamma}_{flavor}} .$$
%$$\mathcal{Z}_{\mathfrak{c}} \rightarrow  \mathcal{Z}_{\mathfrak{c}} .$$

%---

%Let $\mathfrak{g}^{ \bigstar}$ be  the set of formal sums 
%$$ \sum_{ (\gamma, \beta ) \in \Gamma_{KS} } a_{\gamma, \beta}    \hat{e}_{\gamma, \beta} $$
%with the properties
%$$a_{\gamma, \beta} \neq 0 \Rightarrow || \beta || \leq C \omega (\beta),$$
%$$   \# \{   \gamma | \quad a_{\gamma, \beta} \neq 0  \}  < \infty \quad \quad  \forall \beta \in \Gamma  .$$
%$\mathfrak{g}^{ \bigstar}$ has an obvius algebra structure. $\mathfrak{H}$ is a modulo on $\mathfrak{g}^{ \bigstar}$.

%Let $\mathfrak{J}$ the left ideal of $\mathfrak{g}^{ \bigstar}$ generated by $ \{   \hat{e}_{\gamma, 0 } -1 \}_{\gamma \in \text{Ker}(Q) } $. We have an isomorphism of moduli 
%$$ \mathfrak{H} =  \mathfrak{g}^{ \bigstar}/\mathfrak{J} .  $$

\subsection{Semiclassical Limit}

We are now ready to define the semiclassical limit of boundary states. 

\begin{proposition} There exists a natural map 
$$\mathfrak{scl} : \mathfrak{H}^{\flat} \rightarrow \mathfrak{H}^{scl} $$
compatible with $MC$-isotopies.

The following identity holds:
\begin{equation} \label{semiclassical-compatibility}
     \mathfrak{scl} \circ \hat{e}_{\gamma, \beta}^{lie} = e_{\gamma, \beta } \circ \mathfrak{scl} ,
\end{equation}
where $\hat{e}_{\gamma, \beta}^{lie}$ is the operator defined in section \ref{KS-algebra including components section}.
%$$\mathfrak{scl} : \mathcal{Z}/g_s \mathcal{Z} \rightarrow \mathcal{Z}^{disk} .$$
%which is an isomorphism of vector spaces compatible with $MC$-isotopies.
\end{proposition}
\begin{proof}
First notice that the left side of (\ref{jump2}) belongs to $g_s \mathcal{Z}$ if $v_1=v_2$. It follows that on
$\mathcal{Z}/g_s \mathcal{Z}^{}$ we can neglect the relative position of $w_{h,v}$ and $w_{h',v}$ for each $h,h' \in H$ and $v \in V$. Hence, up to small isotopy, on the generator (\ref{generator2}) we can assume that $w_{h,v}$ is independent on  $h \in H$. Set $w_{v}=w_{h,v}$. 
%For the same reason it does not depend on $w^{ann0}$.
We have a linear map 
$$ \mathcal{Z}/g_s \mathcal{Z} \rightarrow \mathcal{Z}^{disk} $$
defined by
$$g_s^{-|V|} (V, H , ( w_{v} )_{v} , ( w^{ann0}), \beta_v) \mapsto (V,  ( w_{v} )_{v} , \beta_v) .$$
%\text{    if      .$$
%$$g_s^{b} (V, H , \{ w_{v} \}_{v} , \beta_v) \mapsto (V,  \{ w_{v} \}_{v} , \beta_v) \text{    if    } b>-|V|$$
 It is immediate to check that this map is compatible with (\ref{MD-jump}) and (\ref{jump2}). 
 %and it is a bijection on the generators.

The identity (\ref{semiclassical-compatibility}) is easy to check from the definition.

\end{proof}


\begin{thebibliography}{10}

%\bibitem{duality} M. Bullimore, T. Dimofte, D. Gaiotto, and J. Hilburn, \emph{Boundaries, Mirror Symmetry, and Symplectic Duality in 3d N = 4 Gauge Theory}, arXiv:1603.0838.


\bibitem{ff1}  M.  Aganagic,  R.  Dijkgraaf,  A.  Klemm,  M.  Marino  and  C. Vafa, \emph{Topological strings and integrable hierarchies},hep-th/0312085.





\bibitem{AV} M. Aganagic and C. Vafa, \emph{Large N Duality, Mirror Symmetry, and a Q-deformed A-polynomial for Knots} , arXiv:1204.4709.


\bibitem{CV} S. Cecotti and C. Vafa, \emph{BPS Wall Crossing and Topological Strings} , arXiv: 0910.2615.

\bibitem{CCV} S. Cecotti, C. Cordova, and C. Vafa, \emph{Braids, Walls, and Mirrors} , arXiv:1110.2115.

\bibitem{CCV2} C. Cordova, S. Espahbodi, B. Haghighat, A. Rastogi, and C. Vafa, \emph{Tangles, generalized Reidemeister moves, and thre-dimensional mirror symmetry}, arXiv:1211.3730.


%\bibitem{GMN1} D. Gaiotto, G. W. Moore and A. Neitzke, \emph{Four-dimensional wall-crossing via threedimensional field theory}, arXiv:0807.4723.

%\bibitem{GMN2} D. Gaiotto, G. W. Moore and A. Neitzke, \emph{Wall-crossing, Hitchin Systems, and the WKB Approximation}, arXiv:0907.3987.

%\bibitem{GMN3} D. Gaiotto, G. Moore, A. Neitzke, \emph{Framed BPS States} , arXiv:1006.0146.

%\bibitem{GMN4} D. Gaiotto, G. Moore, A. Neitzke, \emph{Wall-Crossing in Coupled 2d-4d Systems}, arXiv:1103.2598.

%\bibitem{FO} K. Fukaya; K. Ono, \emph{Arnold conjecture and Gromov-Witten invariant},  Topology  38  (1999),  no. 5, 933--1048.

%\bibitem{FO3} K. Fukaya, Y.-G. OH, H. Ohta, K. Ono, \emph{ Lagrangian intersection Floer theory - anomaly and obstruction}.

\bibitem{KS} M. Kontsevich and Y. Soibelman, \emph{Stability Structures, Motivic Donaldson-Thomas Invariants and Cluster Transformations} , arXiv:0811.2435


\bibitem{CS} V. Iacovino, \emph{Master Equation and Perturbative Chern-Simons theory}, arXiv:0811.2181.




\bibitem{OGW1} V. Iacovino, \emph{Open Gromov-Witten theory on Calabi-Yau three-fold I}, arXiv:0907.5225.

\bibitem{OGW2} V. Iacovino,  \emph{Open Gromov-Witten theory on Calabi-Yau three-folds II}, arXiv:0908.0393.


\bibitem{frame} V. Iacovino, \emph{Frame ambiguity in Open Gromov-Witten invariants}, arXiv:1003.4684.


%\bibitem{OGW3} V. Iacovino, \emph{Open Gromov-Witten theory without obstruction},  arXiv:1711.05302.


\bibitem{OGW3} V. Iacovino,  \emph{Open Gromov-Witten theory without Obstruction}, arXiv:1711.05302.


\bibitem{KSWCF} V. Iacovino, \emph{Kontsevich-Soibelman Wall Crossing Formula and Holomorphic Disks}, arXiv:1711.05306.


\bibitem{QME} V. Iacovino,  \emph{Quantum Master Equation and Open Gromov-Witten theory}, to appear.

%\bibitem{PSPCS} V. Iacovino,  \emph{Point Splicting Perturbative Chern Simons}, in preparation.


\bibitem{BoundaryStates2} V. Iacovino,  \emph{Open Gromov-Witten invariants and Bounrdary States2}, in preparation .





%\bibitem{floer} V. Iacovino, \emph{Holomorphic Lagrangian Floer Homology}, in preparation.


%\bibitem{refined} V. Iacovino, \emph{Refined Wall Crossing and Open Holomorphic Curves}, in preparation.


\bibitem{OV} H. Ooguri, C. Vafa, \emph{Knot invariants and topological strings}, Nuclear Phys. B 577 (2000) 419–438 MR1765411.



\bibitem{W} E. Witten, \emph{Chern-Simons Gauge Theory As A String Theory}, arXiv:hep-th/9207094.


\bibitem{Z} E. Zaslow, \emph{Wavefunctions for a Class of Branes in Three-space}, arXiv:1803.02462.

\end{thebibliography}
\end{document}